\documentclass[a4paper,11pt]{article}%
\usepackage{amssymb,amsfonts,amsmath,amssymb,indentfirst,endnotes,rotating,arydshln,amsthm}
\usepackage[nohead]{geometry}
\usepackage[singlespacing]{setspace}
\usepackage[bottom]{footmisc}
\usepackage{graphicx,array,booktabs,bbm,multicol,colortbl,tabularx}
\usepackage{caption}
\usepackage{authblk}
\usepackage{mathrsfs,mathtools}

\usepackage{physics}
\usepackage{mathrsfs}
\usepackage{wrapfig}
\usepackage[ruled,boxed]{algorithm2e}
\usepackage{subfig}
\usepackage{tikz}
\usetikzlibrary{snakes,arrows,shapes}
\usepackage{float}

\RequirePackage[OT1]{fontenc} \RequirePackage{amsthm,amsmath}
\RequirePackage[numbers]{natbib}
\RequirePackage[colorlinks,citecolor=blue,urlcolor=blue]{hyperref}

\makeatletter
\newcommand*\rel@kern[1]{\kern#1\dimexpr\macc@kerna}
\newcommand*\widebar[1]{%
  \begingroup
  \def\mathaccent##1##2{%
    \rel@kern{0.8}%
    \overline{\rel@kern{-0.8}\macc@nucleus\rel@kern{0.2}}%
    \rel@kern{-0.2}%
  }%
  \macc@depth\@ne
  \let\math@bgroup\@empty \let\math@egroup\macc@set@skewchar
  \mathsurround\z@ \frozen@everymath{\mathgroup\macc@group\relax}%
  \macc@set@skewchar\relax
  \let\mathaccentV\macc@nested@a
  \macc@nested@a\relax111{#1}%
  \endgroup
}
\makeatother

\theoremstyle{example} \theoremstyle{remark} \theoremstyle{lemma}
\theoremstyle{definition} \theoremstyle{corol}
\theoremstyle{proposition} \theoremstyle{condition}
\theoremstyle{assumption}

\newtheorem{theorem}{\n{Theorem}}

\newtheorem{lemma}{\n{Lemma}}

\newtheorem{corollary}{\n{Corollary}}

\newtheorem{proposition}{\n{Proposition}}

\newcommand{\Rmnum}[1]{\expandafter\romannumeral #1}

\font\n=cmcsc12

 
\newcommand{\argmin}{\mathop{\rm argmin~}}
\newcommand{\argmax}{\mathop{\rm argmax~}}

\newcommand{\mle}{\wht \theta_{\mx{\tiny MLE}}}

\newcommand{\wt}{\widetilde}
\newcommand{\wht}{\widehat}

\newcommand{\opnorm}[1]{|\!|\!| #1 | \! | \!|}

\newcommand \bbR{\mathbb{R}}

\def\m{\mathcal}
\def\mb{\mathbb}

\def\mx{\mbox}

\newcommand{\bmp}{\begin{minipage}}
\newcommand{\emp}{\end{minipage}}
\newcommand{\ben}{\begin{eqnarray*}}
\newcommand{\non}{\end{eqnarray*}}

\newcommand{\bc}{\begin{center}}
\newcommand{\ec}{\end{center}}

\begin{document}
\title{\Large Statistical Inference in Mean-field Variational Bayes}
\author[1]{Wei Han\thanks{weih2@illinois.edu}}
\author[1]{Yun Yang \thanks{yy84@illinois.edu}}
\affil[1]{Department of Statistics, University of Illinois Urbana-Champaign}
\date{\vspace{-2em}}

\maketitle

\begin{abstract}
We conduct non-asymptotic analysis on the mean-field variational inference for approximating posterior distributions in complex Bayesian models that may involve latent variables. We show that the mean-field approximation to the posterior can be well-approximated relative to the Kullback-Leibler divergence discrepancy measure by a normal distribution whose center is the maximum likelihood estimator (MLE). In particular, our results imply that the center of the mean-field approximation matches the MLE up to higher-order terms and there is essentially no loss of efficiency in using it as a point estimator for the parameter in any regular parametric model with latent variables. We also propose a new class of variational weighted likelihood bootstrap (VWLB) methods for quantifying the uncertainty in the mean-field variational inference. The proposed VWLB can be viewed as a new sampling scheme that produces independent samples for approximating the posterior. Comparing with traditional sampling algorithms such Markov Chain Monte Carlo, VWLB can be implemented in parallel and is free of tuning.

\end{abstract}


\section{Introduction}\label{sec:introduction}
Variational inference \cite{jordan1999introduction} is a popular computational approach for approximating complicated probability densities that often involve intractable integrals and many latent variables arising in complex Bayesian hierarchical models. In variational inference, the complicated target is approximated by a closest member relative to the Kullback-Leibler (KL) divergence in a pre-specified family of tractable densities.
In many large-scale machine learning applications including clustering problems~\cite{corduneanu2001variational,teschendorff2005variational}, image classification~\cite{pu2016variational,rezende2015variational} and topic models~\cite{miao2017discovering,blei2003latent}, variational inference can be orders of magnitude faster than the traditional sampling based approaches such as Markov Chain Monte Carlo (MCMC). In particular, by turning the integration, or sampling, problem into an optimization problem, variational inference can take advantage of modern optimization tools such as stochastic optimization techniques \cite{kushner2003stochastic,hoffman2013stochastic} and distributed optimization architecture~\cite{ahmed2012scalable,boyd2011distributed} for further improving its efficiency. 

Among various approximating schemes, mean-field approximation is the most common type of variational inference that is conceptually simple, implementation-wise easy and particularly suitable for problems involving large numbers of latent variables. The word ``mean-field" is originated from the mean-field theory in physics where despite complex interactions among many particles in a many (infinite) body system, all interactions to any one particle can be approximated by a single averaged effect from a ``mean-field". In variational inference, by restricting the approximating family of the mean-field to be all density functions that are fully factorized over (blocks of) unknown variables, the associated optimization problem of finding a closest density can be efficiently solved via the (block) coordinate ascent algorithm \cite{bishop2006pattern}.  However, the ease of computation comes at a price of poor approximation as these fully factorized densities in the mean-field family fail to capture any dependence structure among the variables. As noticed in earlier studies \cite{wang2005inadequacy,westling2015establishing}, this disregard of dependence structure may lead to undesirable consequences such as under-estimating the uncertainties if the resulting mean-field densities are blindly used for constructing credible intervals for the parameters. Therefore, variational inference, including the mean-field approximation, are primarily used for rapidly obtaining point estimates in complex Bayesian hierarchical models where traditional methods such as EM algorithms and MCMC are either mathematically intractable ($\rm E$-step in the EM) or computationally inefficient (slow mixing in the MCMC).

Despite the great empirical success achieved by variational inference over the past decades, researchers have not developed much general theory explaining why variational approximation, in particular the mean-field approximation, works so well until recently. Some earlier threads of research characterize their statistical properties in specific problems such as Bayesian linear models~\cite{you2014variational,ormerod2017variational}, Poisson mixed effect models~\cite{hall2011theory,hall2011asymptotic}, stochastic block models~\cite{celisse2012consistency,bickel2013asymptotic,zhang2017theoretical} and normal mixture models~\cite{wang2006convergence}, among others. Many of these studies prove estimation consistency and derive convergence rate of a point estimator based on the variational proxy by
explicitly analyzing the fixed point equation of the variational optimization problem, or directly analyzing the iterative algorithm for solving the optimization problem. In addition, these analyses require the strong
conjugacy assumption on the priors of their models. 

More recently, Wang and Blei~\cite{bleinewresult} prove that the KL minimizer in variational Bayes asymptotically approaches a normal limit in regular parametric models. Their proof uses the $\Gamma$ convergence technique, and is based on a crucial local asymptotic normality (LAN) assumption on the variational objective. This LAN assumption implicitly assumes the estimation consistency and may require a case-by-case verification. Three groups~\cite{alquier2017concentration,Chao2018,yang2017alpha} provide general conditions for deriving the contraction rate of variational approximation as a probability distribution towards the $\delta$-measure at the true parameter of the data generating model, which includes both regular parametric models and infinite-dimensional nonparametric models and implies estimation consistency. Specifically, \cite{alquier2017concentration,Chao2018} focus on models that contain no latent variables, while the theory in \cite{yang2017alpha} can be applied to latent variable models such as normal mixture models. All these results justify the use of variational inference as a valid approach for rapidly obtaining rate-optimal point estimators in complex Bayesian models. However, it remains unclear how good the variational point estimator is when compared to some benchmark, such as the maximum likelihood estimator (MLE) and the posterior mean in regular parametric models---at least theoretically, since the MLE (posterior mean) may be computationally expensive to calculate when the $\rm E$-step (full conditional distributions) does not admit a closed form expression when applying the EM algorithm (Gibbs sampler). In addition, there is little work on how to conduct statistical inference, such as creating credible intervals and performing hypothesis testing in variational procedures.

In this work, we develop a new framework for studying theoretical properties of the mean-field variational approximation and for conducting statistical inference on the model parameters based on the mean-field estimator in parametric models involving latent variables. First, we prove a non-asymptotic result that provides an explicit upper bound on the KL divergence between the mean-field approximation to the marginal posterior of the model parameter and its normal approximation, where the center of the normal is precisely the MLE and the covariance matrix is the diagonal of the inverse of the observed data information matrix (which is the asymptotic covariance of the MLE) plus an extra latent variable information matrix (c.f.~Section~\ref{Sec:contration_MF}). The covaraince structure of the approximating normal limit implies that due to the neglect of the dependence between model parameters and latent variables in the mean-field approximation, the uncertainty under-estimation phenomenon is more severe in models with latent variables than models without (latent variable information is zero).
As a direct consequence of the normal approximation, we show that the mean-field variational estimator, defined as the expectation of the model parameter under the mean-field approximation to the posterior distribution, matches the MLE (or posterior mean) up to a higher-order term relative to the root-$n$ convergence rate. In other words, there is no loss of efficiency (at least asymptotically) in terms of the mean squared error criterion in using the mean-field inference for point estimation.  
Second, we propose a new class of variational weighted likelihood Bootstrap (VWLB) methods for conducting statistical inference on the model parameter via perturbing with random weights the (joint) likelihood function in the mean-field inference in the same spirit as bootstrapping. Interestingly, the VWLB can also be viewed as a new sampling scheme that produces independent samples approximating the marginal posterior of the model parameter. In terms of methodology, VWLB extends the classical ideas of weighted likelihood bootstrap~\cite{wlb} and Bayesian bootstrap~\cite{rubin1981bayesian} to complex Bayesian latent variable models. In terms of computation, the VWLB does not suffer from the slowing mixing issue in MCMC due to the independence of the generated samples and is free of tuning. In addition, unlike the sequential nature of MCMC, sampling via VWLB can be conducted in an embarrassingly parallel manner that has the same time complexity as solving a single variational optimization problem via any distributed learning architecture.

A key ingredient in our proof is a relaxed ``triangle inequality" around the projection of the limiting normal approximation to the posterior for the KL divergence when restricted to the mean-field family (c.f.~Lemma~\ref{lem:KL_identity}). In particular, the mean-field family is not a convex family of distributions, and a strict triangle inequality around the projection of a distribution onto this family with leading factor one (e.g.~Theorem 11.6.1 in~\cite{cover2012elements}) is no longer true. In addition, previous results~\cite{alquier2017concentration,Chao2018} only show a slow polynomial decay on the tail probability of the variational approximation to the posterior (when away from the true parameter), while in order to control the KL-divergence between the variational approximation and its normal approximation, we prove a stronger sub-Gaussian type tail bound (square exponential decay) that uses essential structures of the mean-field family via a variational type analysis (c.f.~Proof of Lemma~\ref{lem:expDacayverify} in Section~\ref{Sec:Proof_lem:expDacayverify}). 

Overall, our results reveal that despite uncertainty under-estimation, point estimators from the mean-field variational inference have essentially no loss of efficiency as the maximum likelihood estimator and also attains the Cram\'{e}r-Rao lower bound in parametric models involving latent variables. In addition, by combining variational inference with bootstrap, the resulting VWLB has the potential of providing a principled and more efficient algorithm for sampling from the posterior in complicated Bayesian hierarchical models.

The rest of the paper is organized as follows. In Section~\ref{sec:method}, we briefly review the mean-field variational inference for approximating the posterior in a general class of Bayesian latent variable models, and present our theoretical results on the non-asymptotic properties of the mean-field approximation. Motivated by these theoretical developments, we propose in Section~\ref{Sec:Inference_MF} a new class of variational weighted likelihood Bootstrap methods for statistical inference via the mean-field approximation, and show in Section~\ref{sec:mainres} the estimation consistency in terms of approximating the target posterior of the model parameter. In Section~\ref{sec:simulation}, we provide two simulations studies, one with latent variable and one without for validating the theory and illustrating the method. Proofs of some selected results are provided in Section~\ref{sec:proofs}. Further details about the simulation and other proofs are deferred to appendices in the supplement material.

\subsection{Notation}
We begin with the notation. As a convention, random variables will be denoted by capital letters, and their realizations by small letters if not otherwise specified. In addition, for each probability measure denoted by a capital letter such as $P$ or $Q$, we use the corresponding small letter $p$ or $q$ to denote its density function as the Radon–Nikodym derivative, depending on the context, either relative to a counting measure when the relevant random variable is discrete, or to the Lebesgue measure when continuous. Depending on whether the density function $p$ is relative to the counting measure or the Lebesgue measure, the integral $\int f(x) p(x) \dd x$ either means a sum $\sum f(x) p(x)$ or the usual Lebesgue integral of $f(x)p(x)$. Throughout the paper, $\|\cdot\|_2$ denotes the usual $\ell_2$ norm in the functional space $L^2(\mathbb{R})$ relative to the Lebesgue measure. Let $D(P\, \|\, Q) = \int \dd P\log(\dd P/\dd Q)$ denote the Kullback–Leibler (KL) divergence, $H(P,\,Q)=\big(\int (\sqrt{\dd P/\dd Q}-1)^2\,\dd Q\big)^{1/2}$ the Hellinger distance, and  $d_{TV}(P,\,Q)=\int |\dd P/\dd Q -1|\,\dd Q$ the total variation distance between two probability measures $P$ and $Q$. 
We use the notation $N(\mu, \Sigma)$ to denote a (multivariate) normal distribution with mean $\mu\in\bbR^k$ and variance-covariance matrix $\Sigma\in\bbR^{k\times k}$. Also, we use $N(a;\mu, \Sigma)$ to denote the cumulative distribution function (cdf) for the above normal distribution at $a\in\bbR$ when the dimension $k$ is one, and $\phi(a;\mu, \Sigma)$ the probability density function (pdf) at $a\in \bbR^k$. We use $\text{diag}(\Gamma)$ to denote the diagonal matrix that has the same diagonal elements as the square matrix $\Gamma$. For an arbitrary square matrix $\Gamma$, $|\Gamma|$ stands for its determinant. For a vector $a\in\mb R^d$, $\|a\|$ denotes its Euclidean $\ell_2$ norm, and for a matrix $A\in\mb R^{d\times d}$, $\opnorm{A}$ denotes its matrix operator norm relative to $\|\cdot\|$.


\section{Non-asymptotic analysis of mean-field variational approximation}\label{sec:method}
In this section, we begin with a brief review on the mean-field variational inference for a class of Bayesian latent variable models. Then we provide two perspectives for explaining the mechanism behind the mean-field approximation. After that, we state our main results in this section providing non-asymptotic analysis of the mean-field variational procedure. Our results imply the estimation consistency, and 
characterize the center and shape of the variational approximation to the exact posterior. In particular, our results reveal that although the mean-field approximation fails to capture the uncertainty, a point estimator obtained as the expectation with respect to the variational distribution matches that of the exact posterior up to high-order terms. This favorable property on the variational mean provides the basis of our inference procedure proposed in the next section.


\subsection{Mean-field variantional inference for Bayesian latent variable models}\label{subsec:revVI}
Let $X^n = (X_1, X_2,\ldots, X_n)$ be composed of independent and identically distributed random variables taking values in $\mathcal X$ from a parametric family $\mathcal P=\{P_{\theta},\, \theta \in \Theta\}$, with $n$ denoting the sample size and $\Theta$ the parameter space as a subset of $\bbR^d$.  In cases such as mixture models, the joint probability distribution of the observations $X^n$ admits a simplified representation by introducing local latent variables $S^n=(S_1,S_2,\ldots,S_n)\in\mathcal S^n$, one per observation, as
\begin{wrapfigure}{R}{0.3\textwidth}
\begin{center}
\begin{tikzpicture}[>=latex',line join=bevel,thick,scale=0.6]
  \pgfsetlinewidth{1bp}
\pgfsetcolor{black}
  \draw [->] (47.326bp,144.2bp) .. controls (44.252bp,136.01bp) and (40.538bp,126.1bp)  .. (33.597bp,107.59bp);
  \draw [->] (57.747bp,144.09bp) .. controls (59.747bp,133.62bp) and (61.997bp,120.12bp)  .. (63.0bp,108.0bp) .. controls (64.32bp,92.055bp) and (64.32bp,87.945bp)  .. (63.0bp,72.0bp) .. controls (62.294bp,63.475bp) and (60.973bp,54.27bp)  .. (57.747bp,35.907bp);
  \draw [->] (33.674bp,72.202bp) .. controls (36.748bp,64.006bp) and (40.462bp,54.102bp)  .. (47.403bp,35.593bp);
\begin{scope}
  \definecolor{strokecol}{rgb}{0.0,0.0,0.0};
  \pgfsetstrokecolor{strokecol}
  \draw (54.0bp,162.0bp) ellipse (18.0bp and 18.0bp);
  \draw (54.0bp,162.0bp) node {$\theta$};
\end{scope}
\begin{scope}
  \definecolor{strokecol}{rgb}{0.0,0.0,0.0};
  \pgfsetstrokecolor{strokecol}
  \draw (27.0bp,90.0bp) ellipse (18.0bp and 18.0bp);
  \draw (27.0bp,90.0bp) node {$s_i$};
  \draw (-0.5,-0.3) rectangle (3.5,4.7); 
\end{scope}
\begin{scope}
  \definecolor{strokecol}{rgb}{0.0,0.0,0.0};
  \pgfsetstrokecolor{strokecol}
  \draw (54.0bp,18.0bp) ellipse (18.0bp and 18.0bp);
  \draw (54.0bp,18.0bp) node {$X_i$};
\end{scope}
\end{tikzpicture}
\end{center}
\caption{Graphical representation of the Bayesian latent variable model.}\label{fig:GMR}
\end{wrapfigure}
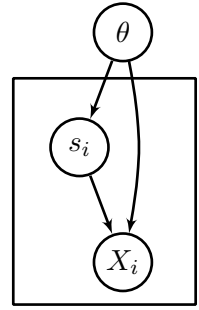
\begin{equation}\label{Eqn:lvm}
p(X^n |\, \theta) = \int_{\mathcal S^n} p(X^n |\, \theta, s^n)  p(s^n |\, \theta) \, \dd s^n,
\end{equation}
where $p(X^n|\,\theta,\, S^n)$ is the conditional density function of $X^n$ given $S^n$, and the joint density $p(X^n\, |\, \theta)$ of $S^n$ is also parametrized by (a subset of) $\theta$ under $P_\theta^n$. 
In other cases, a complex probability model, including the latent Dirichlet allocation and Bayesian hierarchical models, may itself be defined in a hierarchical fashion by first specifying the distribution of the data given latent variables and parameters, and then the latent variable distribution given parameters, as formulated in~\eqref{Eqn:lvm}.  
Due to the negative result on the inconsistency of mean-field variatioanl approximation~\cite{wang2004lack} for general state-space models with non-independent observations with non-independent latent variables, we assume the observation latent variable pair $(X_i,S_i)$ to be mutually independent, that is,
\begin{equation}\label{Eqn:BLV_model}
p(X^n |\, \theta, S^n) = \prod_{i=1}^n p(X_i | \,\theta, S_i)\quad\mbox{and}\quad p(S^n |\, \theta)=\prod_{i=1}^n p(S_i |\, \theta),
\end{equation}
where $p(S_i|\,\theta)$ denotes the marginal density function of $S_i$ parametrized by parameter $\theta$ under $P_\theta^{n}$.

In the Bayesian paradigm, we impose a prior distribution, denoted by $\Pi(\cdot)$, on the model parameter $\theta$ over $\Theta$, whose density function is denoted by $\pi(\cdot)$. 
Figure~\ref{fig:GMR} provides a graphical representation of this Bayesian latent variable model considered in the paper.
In this framework, all inference is based on the posterior probability $p(Z^n |\, X^n)$ of the collection of latent variables $Z^n=(\theta,\, S^n)$ given visible variables $X^n$. 
According to Bayes' theorem, this posterior probability has the following form:
\begin{equation}\label{eqn:pos}
\begin{aligned}
&p(S^n, \theta \,|\, X^n) = \frac{p(X^n, S^n\,|\, \theta)\pi(\theta)}{p(X^n)}, \\ &\qquad\mbox{with}\quad p(X^n) = \int_{\m S^n} p(X^n, s^n\,|\, \theta)\pi(\theta) \,\dd \theta \dd s^n,\\
&\mbox{and (joint) likelihood}\quad  p(X^n, S^n \,|\, \theta) = \prod_{i=1}^n p(X_i \,|\, S_i, \theta) p(S_i\, |\, \theta).
\end{aligned}
\end{equation}
In particular, we are interested in the marginal posterior distribution $\Pi_n$ of the parameter $\theta$ by integrating out $S^n$ in the joint posterior,
\begin{equation}\label{eqn:posterior1}
\Pi_n(A):\,=\Pi(\theta\in A\,|\,X^n) = \frac{\int_A L(\theta;\,X^n)\, \dd \Pi(\theta)}{\int_\Theta  L(\theta;\,X^n) \,\dd \Pi(\theta)}, \quad\mbox{for all measurable set $A\subset \Theta$},
\end{equation}
where the marginal likelihood $L(\theta;\,X^n)$, as a function of $\theta$, is
\begin{equation}\label{eqn:lhf}
L(\theta;\,X^n):\, =p(X^n |\, \theta) = \int_{\m S^n} p(X^n, s^n | \,\theta)\, \dd s^n.
\end{equation}
Unfortunately, in most cases $p(S^n,\theta\,|\,X^n)$ and $\Pi_n(\cdot)$ in equations~\eqref{eqn:pos} and~\eqref{eqn:posterior1} can be inconvenient to use for direct analysis due to the intractable normalization constant $p(X^n)$ involving multi-dimensional integration. Sampling based procedures such as MCMC algorithms could be computationally inefficient due to the high computational cost and slow mixing. Alternatively, variational inference turns the integration problem into an optimization problem by approximating the target distribution $p(Z^n|\,X^n)$ with a closest member $\wht q_{Z^n}$ in a pre-specified family $\Gamma$. 
Formally, the variational approximation $\wht q_{Z^n}$ to $p(Z^n|\,X^n)$ is obtained by solving the following optimization problem,
\begin{equation}\label{eqn:vbObj}
\wht{q}_{Z^n} = \underset{q_{Z^n} \in \Gamma}{\argmin}\, D\big(q_{Z^n}(\cdot) \,\|\, p(\cdot \,| X^n)\big)=\underset{q_{Z^n} \in \Gamma}{\argmin} \int_{\Theta\times\m S^n}q_{Z^n}(z^n)\log \frac{q_{Z^n}(z^n)}{p(z^n|X^n)}\,\dd z^n.
\end{equation}
In particular, we focus on the mean-field approximation where the variational family $\Gamma$ is composed of all fully factorized distributions as
\begin{equation}\label{eqn:mean-field}
q_{Z^n}(z^n) = q_{\theta}(\theta)\,q_{S^n}(s^n)= \prod_{j=1}^d q_{\theta_j}(\theta_j)\, \prod_{i=1}^n q_{S^i}(s_i), \quad z^n=(\theta,\,s^n)\in\mb R^d\times \m S^n.
\end{equation}
Alternatively, one can apply a block mean-field approximation that preserves dependence structures within some multidimensional components, such as the $d$-dim $\theta=(\theta_1,\ldots,\theta_d)$ block, in $Z^n$, and our results can be readily applied to this less stringent scheme.

\subsection{Two perspectives of the mean-field approximation}
The following decomposition of the KL divergence in~\eqref{eqn:vbObj} reveals the interplay between the parameter $\theta$ and latent variables $S^n$ pertaining to the mean-field approximation,
\begin{align}\label{Eqn:KL_decom}
D\big(q_\theta \otimes q_{S^n} \,\|\, p_{\theta,S^n}(\cdot \,| X^n)\big) = D\big(q_{\theta}\,\|\, \pi_n\big) + \int_{\Theta} D\big(q_{S^n}\,\|\, p_{S^n}(\cdot \,|\,\theta,\, X^n)\big)\, q_{\theta}(\theta)\, \dd\theta,
\end{align}
where $p_{\theta,S^n}(\cdot \,| X^n)$ stands for the joint posterior density of $(\theta,\,S^n)$,  $\pi_n$ the density function induced from the marginal posterior distribution $\Pi_n$ of $\theta$ as in~\eqref{eqn:posterior1}, and $p_{S^n}(\cdot \,|\,\theta,\, X^n)$ the conditional posterior density of $S^n$ given $\theta$. Therefore, jointly minimizing the KL-divergence over $(q_\theta,q_{S^n})$ is equivalent to first profiling out the nuisance part $q_{S^n}$ by minimizing the second term for a fixed $q_{\theta}$, and then finding the primary quantity of interest $q_\theta$ that minimizes the resulting ``profile divergence". In particular, this semiparametric profiling perspective plays a crucial role in identifying the limiting center and shape of the variational approximation $\wht q_\theta$  in its normal approximation presented in the following subsection. The lemma below provides an explicit expression for this profile divergence, whose proof is provided in Section~\ref{Sec:proof_lem:decomposeKL}. Recall that $\pi_n$ is the marginal posterior density of $\theta$.

\begin{lemma}\label{lem:decomposeKL}
For each fixed density $q_\theta$ over $\Theta$, the profile divergence takes the following form:
\begin{align*}
& \min_{q_{S^n}=\bigotimes_{i=1}^n q_{S_i}}D\big(q_\theta \otimes q_{S^n} \,\|\, p_{\theta,S^n}(\cdot \,| X^n)\big)
= D(q_{\theta}\, \|\, \pi_n) - \sum_{i=1}^n \log \Big(\int_{\m S} \exp\{r_i(s_i)\}\,\dd s_i\Big),
\end{align*}
where $r_i(s) = \int_{\Theta} \log p(s\,|\,\theta,\,X_i)\, q_\theta(\theta)\,\dd\theta$ for all $s\in\m S$.
\end{lemma} 
\noindent Roughly speaking, as $q_\theta$ approaches the $\delta$-measure at $\theta^\ast$, $r_i(s)$ tends to $\log p(s\,|\,\theta^\ast,X_i)$ and the second term in the preceding display vanishes. Consequently, minimizing the KL-divergence between the joint distributions over $(\theta,S^n)$ boils down to minimizing $D(q_{\theta}\, \|\, \pi_n)$. However, the second term still contributes to the limiting shape of $\wht q_\theta$ as we will see in the next subsection.

The original derivation of the variance inference~\cite{jordan1999introduction} provides an alternative interpretation via Jensen's inequality. Precisely, using the concavity of $\log(x)$, we can obtain an manageable lower bound to the log normalization constant (called evidence) as,
\begin{align}\label{eq:elbo}
\log p(X^n) = \log \int_{\Theta} \frac{p(X^n,\,z^n)}{q_{Z^n}(z^n)}\, q_{Z^n}(z^n)\,\dd z^n  \geq \int_{\Theta} \log \frac{p(X^n,\,z^n)}{q_{Z^n}(z^n)}\, q_{Z^n}(z^n)\,\dd z^n :\,= L(q_{Z^n}),
\end{align}
where $L(q_{Z^n})$ is called the evidence lower bound (ELBO,~\cite{Review}). In particular, the KL divergence
\begin{align*}
D\big(q_{Z^n}(\cdot) \,\|\, p(\cdot \,| X^n)\big)=\log p(X^n)-L(q_{Z^n}),
\end{align*}
quantifies the discrepancy between the evidence $\log p(X^n)$ and its lower bound approximation $L(q_{Z^n})$. Consequently, minimizing the KL divergence in optimization problem~\eqref{eqn:vbObj} is equivalent to finding a best  $q_{Z^n}$ to maximize the ELBO. The KL minimization formulation~\eqref{eqn:vbObj} is convenient for our theoretical analysis, while the ELBO formulation leads to various computational algorithms for implementing the variational inference.

In this paper, our attention toward the model is inference on $\theta$, the model parameter, where our theory and methodology is centered on. Towards this goal, it is helpful to inspect a finer decomposition of the ELBO from~\cite{yang2017alpha},
\begin{equation}\label{Eqn:ELBO_decmop}
\begin{aligned}
L(q_{Z^n})& = \int_{\Theta} \log p(X^n|\,\theta)\, Q_\theta(\dd\theta)-\Delta_J(q_\theta, \,q_{S^n}) - D(q_\theta\,||\,\pi_\theta),\quad \mbox{with}\\
\Delta_J(q_\theta, \,q_{S^n})&= \int_{\Theta} \bigg[ \log p(X^n|\,\theta) - \underbrace{\int_{\m S^n} \log\frac{p(X^n,\,s^n |\theta)}{q_{S^n}(s^n)}\,q_{S^n}(s^n)\,\dd s^n}_{\widehat{\log p(X^n|\,\theta)}}  \bigg]\, Q_\theta(\dd\theta) \geq 0,
\end{aligned}
\end{equation}
which consists of three terms: an integrated (relative to the variational distribution of $\theta$) log-marginal likelihood, the Jensen gap $\Delta_J$ due to the mean-field decomposition on latent variables $\{S_i\}_{i=1}^n$ in approximating the marginal likelihood $p(X^n|\,\theta)$ with $\widehat{p(X^n|\,\theta)}$, and the KL divergence between the variational distribution $q_{\theta}$ and the prior $\pi(\theta)$. When there is no likelihood approximation with latent variables, the Jensen gap $\Delta_J$ term vanishes, and maximizing of the ELBO value in decomposition~\eqref{Eqn:ELBO_decmop} resembles a regularized $M$-estimation problem of minimizing an objective function composed of a goodness of fit term $\int_{\Theta} -\log p(X^n|\,\theta)\, Q_\theta(\dd\theta)$ plus a regularizing term $D(q_\theta\,||\,\pi_\theta)$ over all distributions in the variational family $\Gamma$. This perspective is useful in proving the consistency and characterizing the contraction rate of the variational approximation $\wht q_\theta$ towards the true parameter $\theta^\ast$, as describe in the next subsection.

\subsection{Contraction of mean-field variational approximation and normal approximation}\label{Sec:contration_MF}
In this subsection, we introduce our non-asymptotic results characterizing the contraction rate and shape of $\wht q_\theta$ in the mean-field variational inference. Our analysis is under the frequentist perspective by assuming the observations $\{X_i\}_{i=1}^n$ as i.i.d.~copies from a data generating model $P_{\theta^\ast}$, where $\theta^\ast\in\Theta\subset \mb R^d$ is referred to as the truth parameter, or simply truth. We make following assumptions.

\vspace{0.2em}
\noindent{\bf Assumption A1 (Prior continuity and growth):} The prior density satisfies $\pi(\theta^\ast)>0$. In addition, $\log\pi(\theta)$ is differentiable in a neighborhood of $\theta^\ast$, and satisfies
\begin{align*}
|\log \pi(\theta) - \log \pi(\theta^\ast)|\leq C(1+\|\theta-\theta^\ast\|^L),\quad\forall \theta\in\Theta\ \ \mbox{ for some constant $L>0$}.
\end{align*}

\noindent We define the following quantity of Hellinger bracketing entropy that provides a measure on the model space complexity.

\vspace{0.2em}
\noindent{\bf Definition (Hellinger  bracketing entropy):}
For a set $\m F$ of functions over $\m X$ and any $\varepsilon>0$, we call a set (of pairs of functions) $\{ (f_j^L, f_j^U, j = 1,\ldots,N)\}$ a (Hellinger) $\varepsilon$-bracketing of $\m F$, if $\int_{\m X} \big[(f_j^L)^{1/2}(x)-(f_j^U)^{1/2}(x)\big]^2\dd x\leq \varepsilon^2$ for $j=1,2,\ldots,N$ 
and for any $f\in \m F$, there is a $j$ such that $f_j^L \leq f \leq f_j^U$. The (Hellinger) $\varepsilon$-bracketing metric entropy, denote by $H_B(\varepsilon,\m F)$, is defined as the logarithm of the smallest cardinality of such an $\varepsilon$-bracketing of $\m F$.

\vspace{0.5em}
\noindent{\bf Assumption A2 (Marginal likelihood regularity):}
\begin{enumerate}
\item (Smoothness) The log-marginal likelihood function $l(\theta; x)=\log p(x\,|\,\theta)$ is thrice continuously differentiable with respect to $\theta$. 

\item (Finite moments) In a neighborhood $\m B(\theta^\ast;\delta)=\{\theta\in\mb R^d:\,\|\theta-\theta^\ast\|\leq \delta\}$ of $\theta^\ast$ the fourth moments of the derivatives at $\theta^\ast$ up to order three exist under $P_{\theta^\ast}$. Moreover, there exists a measurable function $M:\,\m X\to \mb R_+$ satisfying $E_{\theta^\ast} M^2(X) < \infty$, such that the third order derivatives of $l(\theta; x)$ satisfies 
\begin{align*}
\max_{j,k,l\in[d]} \Big|\frac{\partial^3 l(\theta; x)}{\partial \theta_j\partial \theta_k\partial \theta_l}\Big| \leq M(x)(1+\|\theta-\theta^\ast\|^L),\quad \mbox{for all $x\in\m X$ and $\theta\in\Theta$}.
\end{align*}

\item (Information matrix non-degeneracy) In addition, the order of taking expectation with respect to $P_{\theta^\ast}$ and differentiation at $\theta^\ast$ is valid so that 
\begin{align*}
E_{\theta^\ast} \big[\nabla l(\theta^\ast;\, X)\nabla l(\theta^\ast;\, X)^T\big] = - E_{\theta^\ast} \nabla^2 l(\theta^\ast;\, X).
\end{align*}
The $d\times d$ Fisher information matrix in this display, denoted by $I(\theta^\ast)$, is positive definite. 

\item (Euclidean metric equivalence) The squared Hellinger distance satisfies that for some constants $(c_1,c_2)$,
\begin{align*}
c_1\,\|\theta-\theta^\ast\|^2 \leq H^2\big(P_\theta,\,P_{\theta^\ast}\big)\leq c_2\,\|\theta-\theta^\ast\|^2,\quad\mbox{for all }\theta\in\Theta.
\end{align*}

\item (Local metric entropy growth) There exists a constant $c_3$, such that the Hellinger entropy satisfies
\begin{equation*}
H_B\big(u, \big\{p_\theta,\,\theta\in\Theta:\, H(P_\theta,\, P_{\theta^\ast})\leq s\big\} \big) \leq c_3 \log \Big(\frac{s}{u}\Big),  \quad \forall u \in (0, s)\ \ \mbox{and} \ \ s\in[0,1].
\end{equation*}
\end{enumerate}
The first three assumptions in A2 are standard regularity conditions for parametric models (c.f.~Chapter 1.4 in~\cite{ghoshbayesian}).
The Euclidean metric equivalence assumption is also made in Theorem 5.1 in~\cite{ghosal2000convergence} as one of their sufficient conditions for proving posterior contraction in parametric models. The last assumption on the local metric entropy assumption is adopted from~\cite{RatesPosterior}, and often holds for parametric models (c.f.~\cite{RatesPosterior} for examples).

\vspace{0.2em}
\noindent{\bf Assumption A3 (Latent conditional density regularity):} The log-conditional density $l_s(\theta; s,x)=\log p(s\,|\,\theta,\,x)$ of the latent variable $S$ given $X$ is thrice differentiable 
with respect to $\theta$ in the neighborhood $\m B(\theta^\ast;\delta)$ of $\theta^\ast$, and the fourth moments of the derivatives at $\theta^\ast$ up to order three exist under $P_{\theta^\ast}$.  Moreover, there exists a measurable function $M_s:\,(\m S,\m X)\to \mb R_+$ satisfying $E_{\theta^\ast} M_s^2(S,X) < \infty$ such that the third order derivatives of $l_s(\theta; s,x)$ satisfies 
\begin{align*}
\max_{j,k,l\in[d]} \Big|\frac{\partial^3 l_s(\theta; s,x)}{\partial \theta_j\partial \theta_k\partial \theta_l}\Big| \leq M(s,x),\quad \mbox{for all $s\in\m S$, $x\in\m X$ and $\theta\in\m B(\theta^\ast;\delta)$.}
\end{align*}
In addition, the $d\times d$ latent (variable) information matrix $I_s(\theta)=E_{\theta^\ast}[\nabla^2 l_s(\theta; S,X)]$ is locally Lipschitz in the neighborhood $\m B(\theta^\ast;\delta)$ of $\theta^\ast$, that is,
\begin{align*}
\opnorm{I_s(\theta)-I_s(\theta^\ast)} \leq L_s \|\theta-\theta^\ast\|\quad \mbox{for all $\theta\in\m B(\theta^\ast;\delta)$.}
\end{align*}

\vspace{0.2em}
\noindent Assumption A3 includes regularity conditions on the conditional distribution of latent variables.
In particular, by viewing the latent variable $S$ as missing data, and interpreting $I_s(\theta^\ast)$ as the missing data information matrix and $I(\theta^\ast)$ as the observed data information matrix~\cite{woodbury1970missing}, we can define the complete data information matrix as $I_c(\theta^\ast) = I(\theta^\ast)+I_s(\theta^\ast)$. As we will see, the inverse of ${\rm diag}(I_c(\theta^\ast))$ characterizes the limiting shape of the variational approximation~$\wht q_\theta$, where the second term $I_s(\theta^\ast)$ causes the extra variance reduction due to the neglect of the posterior dependence between $\theta$ and $S^n$. 

The following lemma shows that with high probability, the marginal distribution $\widehat{Q}_\theta$ obtained from the mean-field variational approximation~\eqref{eqn:vbObj} has a sub-Gaussian tail probability outside an $\varepsilon$-ball centered at the truth $\theta^\ast$, for all $\varepsilon\geq \sqrt{\log n/n}$. This exponentially decaying tail behavior is essential for controlling the tail integrals in the proof of our next result that approximates $\widehat Q_\theta$ with a normal distribution. A proof is provided in Section~\ref{Sec:Proof_lem:expDacayverify}.

\begin{lemma}\label{lem:expDacayverify}
Under Assumptions A1, A2 and A3, there exist constants $(C_0,C_1,C_2,C_3)$ such that for any $M\geq 1$ and $\varepsilon_n=C_0 \sqrt{\frac{\log n}{n}}$ it holds with probability at least $1-C_1M^{-2}$ that the mean-field approximation $\widehat Q_{\theta}$ satisfies
\begin{align*}
\widehat Q_{\theta} \big( \|\theta - \theta^\ast\|\geq C_2 \varepsilon\big) \leq e^{-C_3n\varepsilon^2},\quad\mbox{for all } \ \varepsilon\geq M\varepsilon_n.
\end{align*}
\end{lemma}
\noindent Although we focus on the finite-dimensional parametric model $P_\theta$, the proof of this lemma is based on a general treatment under a similar setting as~\cite{RatesPosterior,ghosal2000convergence,yang2017alpha}. Therefore, this result can also be extended to mean-field approximations for infinite-dimensional models, for which the same sub-Gaussian tail bound holds for all $\varepsilon$ greater than a benchmark contraction rate $\varepsilon_n$ slower than the parametric root-$n$ rate pertaining to the model by making certain assumptions (c.f.~conditions~(2.2)--(2.4) in~\cite{ghosal2000convergence}) on the prior thickness and the complexity of model space. In comparison, earlier results on the consistency and convergence rates of variational approximations, such as~\cite{yang2017alpha,Chao2018}, only show a polynomially decay $C_3\varepsilon_n^2/\varepsilon^2$ on the tail probability $\widehat Q_{\theta} \big( \|\theta - \theta^\ast\|\geq C_2 \varepsilon\big)$. The proof of our exponentially decaying bound utilizes the factorization structure~\eqref{eqn:mean-field} of the mean-field approximation, and it is still an interesting open problem whether similar sub-Gaussian type tail bounds hold for a broader class of variational approximations beyond the mean-field.

Let $\mle$ denote the maximum likelihood estimator of $\theta$,
\begin{align*}
\mle = \argmax_{\theta\in\Theta} \sum_{i=1}^n\log p(X_i\,|\,\theta).
\end{align*}
The classical Bernstein von-Mises (BvM) theorem \cite{VaartTextbook,ghoshbayesian} states that the marginal posterior distribution $\Pi_n$ of $\theta$ approaches in the total variation metric to $N\big(\mle,\, [nI(\theta^\ast)]^{-1}\big)$ as $n\to\infty$ (our Lemma~\ref{lem:KL_normal_approx} with $Q_\theta=\Pi_n$ gives a stronger KL divergence version of the BvM theorem). Our next theorem shows that the marginal variational distribution $\widehat Q_\theta$ can also be approximated by a normal distribution with the same center as $\mle$, but a different variance-covariance matrix, under the stronger KL-divergence. A proof is deferred to Section~\ref{Sec:proof_thm1}. Recall that $I_c(\theta^\ast)=I(\theta^\ast) + I_s(\theta^\ast)$ is the complete data information matrix. 

\begin{theorem}\label{thm:1}
Under Assumptions A1, A2 and A3, there exist constants $(C_4,C_5)$ such that for any $M\geq 1$ it holds with probability at least $1-C_4M^{-2}$ that 
\begin{equation*}
D(\widehat{Q}_\theta \,\| \, Q^\ast_{VB}) \leq  \frac{C_5M^3(\log n)^{d + 3}}{\sqrt{n}},
\end{equation*}
where $Q^\ast_{VB} = N\big(\mle,\, [nI_{VB}]^{-1}\big)$, and $I_{VB} =\mx{Diag}\big(I_c(\theta^\ast)\big)$.
\end{theorem}
\noindent Theorem~\ref{thm:1} is a non-asymptotic result that applies to any sample size $n\geq 1$. The complementary probability $C_4M^{-2}$ decays polynomially in $M$ because we simply apply the Markov inequality with the second order moment assumption on the derivatives of the log-likelihood function. If we instead make a sub-Gaussian type assumption as in~\cite{spokoiny2012parametric}, then this remainder probability will be exponentially small in $M^2$ as $\exp\{-C_4 M^2\}$. 
As a special when there is no latent variables ($I_s(\theta^\ast)=0$), Theorem~\ref{thm:1} shows that 
the mean-field approximation $\widehat Q_\theta$ tends to the normal distribution $N\big(\mle,\, [nI_{VB}(\theta^\ast)]^{-1}\big)$ whose covariance matrix simply removes all off-diagonal components in $I(\theta^\ast)$, which is consistent with earlier results such as~\cite{bleinewresult} and explains the 
overly small variances exhibited by the mean-field approximation~\cite{wang2005inadequacy,westling2015establishing} due to the neglect of the dependence among components of $\theta$. In the general case of Bayesian latent variable models, Theorem~\ref{thm:1} shows that the overly small variances phenomenon is even more severe due to the neglect of the dependence between $\theta$ and $S^n$. 

In practice, although this mismatch on the covariance structures between the variational approximation and the exact posterior is not a serious issue when doing point estimation, erroneous characterization of uncertainty can be produced. As a consequence, variational inference is widely used for rapidly obtaining a point estimator for the model parameter $\theta$. 
Let $\wht\theta_{VB}$ to denote the variational posterior mean $\wht\theta_{VB}=\int_{\Theta} \theta \,\wht q_\theta(\theta)\,\dd\theta$. 
The following corollary, as a direct consequence of the normal approximation in Theorem~\ref{thm:1}, shows that although the shape of the exact posterior $\Pi_n$ is not properly captured by $\widehat Q_\theta$, their centers $\mle$ and $\wht\theta_{VB}$ match up to $\m O(n^{-3/4})$.
\begin{corollary}\label{cor:VB_center}
Under the conditions and the high-probability event of Theorem~\ref{thm:1}, there exists a constant $C_6$ such that
\begin{align*}
\|\sqrt{n}\,( \wht\theta_{VB}-\mle)\| \leq \frac{C_6M^{3/2}(\log n)^{d/2 + 3/2}}{n^{1/4}}.
\end{align*}
\end{corollary}
\noindent This corollary implies that there is essentially no loss of efficiency in using the mean-field approximation as a fast approach for obtaining a point estimator in low-dimensional parametric models. Moreover, this interesting finding suggests that we can also conduct statistical inference in mean-field approximation by Bootstrapping the point estimator $\wht\theta_{VB}$.


\section{Statistical inference in mean-field approximation}\label{Sec:Inference_MF}
Motivated by results in the previous section, we propose an inferential framework for mean-field variational Bayes in Bayesian models with latent variables by borrowing the classical idea of weighted likelihood Bootstrap (WLB,~\cite{wlb}) for approximate Bayesian computation.
We begin this section with a brief review on the original WLB as a way to simulate approximately from a posterior distribution when there is no latent variables. After that, we extend the WLB to incorporate latent variables, which further leads to our variational weighted likelihood Bootstrap (VWLB) for approximating the marginal posterior distribution of $\theta$ in the mean-field variational Bayes. 

\subsection{Weighted likelihood Bootstrap}
WLB is an extension of the Bayesian Bootstrap~\cite{rubin1981bayesian} from nonparametric models to parametric and semiparametric models by approximating the exact posterior via a random sample of parameter values, each maximizing a weighted likelihood function with random weights.
The original WLB proposed in~\cite{wlb} directly operates on the marginal density function $p(\cdot\,|\,\theta)$ of the i.i.d.~observations $\{X_i\}_{i=1}^n$ without introducing the latent variables $\{S_i\}_{i=1}^n$. More specifically, in WLB the $b$th random sample $\wt\theta^{(b)}$ of the parameter $\theta$, for $b=1,2,\ldots,B$, is produced by maximizing the following weighted likelihood function obtained from tilting the likelihood function $L(\theta;\,X^n)$:
\begin{align}\label{Eqn:wmlilelihood}
\widetilde L^{(b)}(\theta;\, X^n)=\prod_{i=1}^n \big[p(X_i\,|\,\theta)\big]^{W_i^{(b)}},\quad\mx{for all }\theta\in\Theta,
\end{align}
where the weights $W^{n,(b)}=(W_1^{(b)},W_2^{(b)}\ldots,W_n^{(b)})$ satisfy the following assumption.

\vspace{0.2em}
\noindent{\bf Assumption W (Weight randomness):} The weights $\{W_i^{(b)}:\,i=1,2,\ldots,n,\, b=1,2,\ldots,B\}$ are i.i.d.~copies of a nonnegative random variable $W$ with $E[W]=\mx{Var}(W)=1$. In addition, $W$ is sub-exponential, that is, there exist some constants $(c_0,c_1)$ such that $E[e^{\lambda(W-1)}]\leq e^{c_0\lambda^2/2}$ holds for all $|\lambda|\leq c_1$.
\vspace{0.2em}

Unlike other commonly used sampling schemes such as MCMC, the random samples $\{\wt\theta^{(b)}\}_{b=1}^B$ from WLB are conditionally independent given the data $X^n$, where the extra randomness in $\wt\theta^{(b)}$ is induced by the distribution of the random weights. Under some mild conditions on the model, one can show that the conditional distribution of $\wt\theta^{(b)}$ given data $X^n$ approaches the exact posterior distribution $p(\theta\,|\,X^n)$ of $\theta$ as $n\to\infty$ (readers may refer to~\cite{wlb} for more details on the WLB and its accompanied theory). 
To accommodate this WLB idea to variational inference, it is helpful to also associate the weighted likelihood function $\widetilde L^{(b)}(\theta;\, X^n)$ with a weighted posterior distribution 
 \begin{equation}\label{eq:weightedPo}
\widetilde{\Pi}^{(b)}_n(A):\, = \wt \Pi^{(b)}(A\,|\,X^n) = \frac{\int_A \widetilde L^{(b)}(\theta;\, X^n)\, \dd \Pi(\theta)}{\int_\Theta \widetilde L^{(b)}(\theta;\, X^n)\, \dd \Pi(\theta)},\quad\mx{for all measureable set }A\subset \Theta.
\end{equation}
This weighted posterior can be viewed as a generalization of the fractional posterior~\cite{o1995fractional,bhattacharya2019bayesian}
by raising the probability density $p(X_i\,|\,\theta)$ of $X_i$ in the likelihood to a sample specific power $W_i^{(b)}$. Similar to the classical BvM theorem, the following proposition shows that the mean $\wt \theta^{(b)}_{B}$ of the weighted posterior distribution $\widetilde{\Pi}^{(b)}_n$ based on the marginal likelihood matches the maximum weighted likelihood estimator $\wt \theta^{(b)}$ up to a higher-order remainder term. 

\begin{proposition}\label{prop:wposterior_mode}
Under Assumptions A1, A2 and W, there exist constants $(C_7,C_8)$, such that for any $M\geq 1$, it holds with probability at least $1-C_7M^{-2}$ that
\begin{align*}
\|\sqrt{n}\,(\wt\theta^{(b)}_{B}-\wt\theta^{(b)})\| \leq C_8 \frac{M \log n}{\sqrt{n}}.
\end{align*}
\end{proposition}

\subsection{Variational weighted likelihood Bootstrap}
In this subsection, we propose variational weighed likelihood Bootstrap (VWLB) as a variational approximation method for simulating random samples from the marginal posterior distribution $\Pi_n$ of $\theta$ in the Bayesian latent variable model~\eqref{Eqn:BLV_model}, thereby facilitating statistical inference on parameter $\theta$. Motivated by the MLB method described in the previous subsection, we define the weighted joint likelihood function $\wt p(X^n, S^n | \,\theta)$ and weighted joint posterior density $\wt p(S^n,\theta\,|\,X^n)$ that incorporate latent variables $S^n$ as,
\begin{align}
\wt p^{(b)}(X^n, S^n | \,\theta) &= \prod_{i=1}^n \big[p(X_i\, |\, S_i, \theta) \,p(S_i \,|\, \theta)\big]^{W_i^{(b)}},\label{eqn:wlhf}\\
\wt p^{(b)}(S^n,\theta\,|\,X^n) &= \frac{\wt p^{(b)}(X^n, S^n | \,\theta)\,\pi(\theta)}{\int_\Theta \wt L^{(b)}(\theta;\, X^n)\, \pi(\theta)\,\dd \theta},\quad\mx{for all }\theta\in\Theta,\label{eqn:wposterior}
\end{align}
where recall that $\wt L^{(b)}(\theta;\, X^n)$ is the weighted (marginal) likelihood function defined in~\eqref{Eqn:wmlilelihood}. Note that the ``marginalization" of $S^i$ in the denominator of~\eqref{eqn:wposterior} is before raising to the power $W^{(b)}_i$. Therefore, the weighted joint posterior density is not properly normalized and strictly speaking, not a real density function. Similar to the optimization problem~\eqref{eqn:vbObj} of variational approximation, we define the weighted variational approximation $\wt q^{(b)}_{Z^n}=\wt q_\theta^{(b)} \otimes\bigotimes_{i=1}^n \wt q^{(b)}_{S_i}\in\Gamma$ to $\wt p^{(b)}(Z^n|\,X^n)$ as
\begin{equation}\label{eqn:vwlbObj}
\begin{aligned}
\wt{q}^{(b)}_{Z^n} =&\, \underset{q_{Z^n} \in \Gamma}{\argmin}\wt D^{(b)}\big(q_{Z^n}(\cdot)\,||\,\wt p^{(b)}(\cdot\,|\,X^n)\big)\\
&\quad\quad:\,= \underset{q_{Z^n} \in \Gamma}{\argmin} \int_{\Theta\times\m S^n}  q_{Z^n}(z^n)\log \frac{q_{\theta}(\theta)\prod_{i=1}^n\big[q_{S_i}(s_i)\big]^{W_i^{(b)}} }{\wt p^{(b)}(z^n|X^n)}\,\dd z^n,
\end{aligned}
\end{equation}
where $\wt D^{(b)}\big(q_{Z^n}(\cdot)\,||\,\wt p^{(b)}(\cdot\,|\,X^n)\big)$ calculates the expectation with respect to $q_{Z^n}$ of the log-ratio between the tilted $q_{Z^n}$ and $\wt p^{(b)}(Z^n\,|\,X^n)$. It is worthy noticing that in practical implementations, the denominator $\wt p^{(b)}(z^n|X^n)$ in the preceding display simply contributes a constant independent of $q_{Z^n}$ in above objective function, and there is no need to explicitly compute this quantity. We include this term mainly for theoretical purposes---the $\wt D^{(b)}(\cdot\,||\,\cdot)$ with this term reduces to the usual KL divergence when $W_i^{(b)}\equiv 1$. 
In addition, the weighted variational approximation leads to the following decomposition that generalizes the unweighted KL decomposition formula~\eqref{Eqn:KL_decom} and plays an essential role in our theoretical analysis,
\begin{align}\label{Eqn:WKL_decom}
\wt D^{(b)}\big(q_{Z^n}(\cdot)\,||\,\wt p^{(b)}(\cdot\,|\,X^n)\big) = D\big(q_{\theta}\,\|\, \wt \pi^{(b)}_n\big) + \sum_{i=1}^n W_i^{(b)}\int_{\Theta} D\big(q_{S_i}\,\|\, p_{S_i}(\cdot \,|\,\theta,\, X^n)\big)\, q_{\theta}(\theta)\, \dd\theta.
\end{align}
\noindent This identity implies the weighted variational objective function to be nonnegative. Accordingly, this decomposition formula enables us to divide the joint minimization problem~\eqref{eqn:vwlbObj} into two steps---first profiling out the nuisance part $q_{S^n}$ by minimizing the second term for a fixed $q_\theta$, and then minimizing the resulting ``weighted profile divergence" as a function of $q_\theta$. In particular, the first step of minimizing over $q_{S^n}$ admits a convenient closed form expression for the theoretical analysis, as summarized in the following lemma.
\begin{lemma}\label{lem:decomposeWKL}
For each fixed density $q_\theta$ over $\Theta$, the weighted profile divergence takes the form as:
\begin{align*}
& \min_{q_{S^n}=\bigotimes_{i=1}^n q_{S_i}}\wt D^{(b)}\big(q_{Z^n}(\cdot)\,||\,\wt p^{(b)}(\cdot\,|\,X^n)\big) 
= D(q_{\theta}\, \|\, \wt \pi^{(b)}_n) - \sum_{i=1}^n W_i^{(b)}\log \Big(\int_{\m S} \exp\{r_i(s_i)\}\,\dd s_i\Big),
\end{align*}
where functions $\{r_i(\cdot)\}_{i=1}^n$ are defined in Lemma~\ref{lem:decomposeKL}.
\end{lemma} 
\noindent The same remark after Lemma~\ref{lem:decomposeKL} regarding the limiting behavior of $r_i$ and its implications as $q_\theta$ approaches to the $\delta$-measure at $\theta^\ast$ applies to the weighted case. 
Note that the decomposition formula in the lemma is convenient for deriving a normal approximation to $\wt q^{(b)}_\theta$ and cannot be directly used for practical computation since the weighted marginal posterior $\wt \pi^{(b)}_n$in the first term is computationally intractable. In addition, it is worthy mentioning that a similar weighted ELBO decomposition generalizing~\eqref{Eqn:ELBO_decmop} holds by replacing the first goodness of fit term and the second Jensen gap term with their weighted counterparts. This weighed ELBO decomposition would be useful for deriving the contraction rate of the weighted variational approximation $\wt q^{(b)}_\theta$ beyond parametric models by adapting the proof techniques in~\cite{yang2017alpha}.

Having obtained the weighted variational approximation $\wt q^{(b)}_\theta$ to the weighted posterior $\wt\pi_n$ of $\theta$ with the $b$th random weights $W^{n,(b)}$, we can then proceed as in the usual variational inference by using the variational mean $\wt\theta^{(b)}_{VB}=\int_\Theta \theta \,\wt q^{(b)}_\theta(\theta)\,\dd\theta$ as the $b$th random sample approximately drawn from the marginal posterior $\pi_n$, for $b=1,2,\ldots,B$. Unlike MCMC sampling algorithms, the random samples $\{\wt\theta^{(b)}_{VB}\}_{b=1}^B$ from the variational WLB are i.i.d.~draws approximately from $\pi_n$ given data $X^n$. These random samples can be used for statistical inference, such as constructing credible sets and conducting hypothesis testing.
Algorithm~\ref{algo:vbA} below summarizes the pseudo-code for implementing the variational WLB.
\begin{figure}[ht]
  \centering
  \begin{minipage}{.9\linewidth}
\begin{algorithm}[H]
\SetAlgoLined
\KwIn{Number $B$ of random draws}
\KwData{$X_1,X_2, \ldots, X_n$}
\For{$b \leftarrow 1$ \KwTo $B$}{
  \emph{Generate random weights} $W_1^{(b)},W_2^{(b)}, \ldots, W_n^{(b)}$ \;
  \emph{Solve optimization problem~\eqref{eqn:vwlbObj} to obtain variational density $\wt q^{(b)}_\theta$, for example, via the weighted CAVI in Algorithm~\ref{algo:vbB}}\;
  \emph{Compute weighted variational mean $\wt\theta^{(b)}_{VB}=\int_\Theta \theta \,\wt q^{(b)}_\theta(\theta)\,\dd\theta$}\;
}

\KwOut{Random samples $\{\wt\theta^{(b)}_{VB}\}_{b=1}^B$ for approximating posterior $\pi_n$ of $\theta$}
\caption{Variational Weighted Likelihood Bootstrap}
\label{algo:vbA}
\end{algorithm}
\end{minipage}
\end{figure}

\subsection{Coordinate ascent algorithm for computation}
In this subsection, we discuss computational aspects of the optimization problem~\eqref{eqn:vwlbObj} in the inner loop of Algorithm~\ref{algo:vbA} via a weighted variant of the coordinate ascent.
Coordinate ascent variational inference (CAVI,~\cite{bishop2006pattern}) is a popular optimization algorithm tailored for solving~\eqref{eqn:vbObj} in the usual mean-field approximation that is scalable to large datasets.
CAVI, as an optimization counterpart of Gibbs sampling, utilizes the special structure of the mean-field solution $\wht q_{Z^n}$ to~\eqref{eqn:vbObj} that based upon the optimality, each factor in the decomposition~\eqref{eqn:mean-field} should be proportional to the exponential of the expected log of the joint posterior with respect to the rest factors, which under our notation, simplifies to
\begin{equation}\label{Eqn:CAVI}
\begin{aligned}
\wht q_{\theta_j}(\theta_j)& \propto \exp\Big\{E_{\wht{q}_{-\theta_j}}\big[\log \pi(\theta) + \log p(S^n, X^n|\,\theta)\big]\Big\},\quad\theta\in\Theta, \quad \mx{for $j=1,2,\ldots,d$,\quad and}\\
\wht q_{S_i}(s_i)& \propto \exp\Big\{E_{\wht{q}_{\theta}}\big[\log p(s_i\,|\,\theta,X_i)\big]\Big\},\quad s_i\in\m S,\quad \mx{for $i=1,2,\ldots,n$,}
\end{aligned}
\end{equation}
where the notation $E_{\wht{q}_{-\theta_j}}$ stands for taking expectation with respect to all factors in~\eqref{eqn:mean-field} except for $\wht q_{\theta_j}$ (similar notation convention applies to the weighted case below). CAVI iteratively updates each factor $q_{\theta_j}$ or $q_{Z_i}$ until convergence.  Since the ELBO value is non-decreasing along the iterations, CAVI is guaranteed to converge to a local minimum. In practice, the convergence of the CAVI can be assessed by monitoring the ELBO value, and multiple random initializations can be deployed for finding the global minimum by picking one that yields the highest ELBO value.

Now we generalize the CAVI to its weighted counterpart for solving optimization~\eqref{eqn:vwlbObj}. More specifically, the optimality condition of the optimization problem~\eqref{eqn:vwlbObj} with random weights $W^{n,(b)}=(W^{(b)}_i)$ is
\begin{align}
\wt q^{(b)}_{\theta_j}(\theta_j)& \propto \exp\Big\{E_{\wt{q}_{-\theta_j}^{(b)}}\big[\log \pi(\theta) + \log\wt p(S^n, X^n|\,\theta)\big]\Big\},\quad\theta\in\Theta, \quad \mx{for $j=1,2,\ldots,d$,\quad and}\label{eq:wcavig} \\
\wt q_{S_i}^{(b)}(s_i)& \propto \exp\Big\{E_{\wt{q}_{\theta}^{(b)}}\big[\log p(s_i\,|\,\theta,X_i)\big]\Big\},\quad s_i\in\m S,\quad \mx{for $i=1,2,\ldots,n$,}\label{eq:wcavil}
\end{align}
which only differs from the optimality condition~\eqref{Eqn:CAVI} of the usual mean-field optimization in replacing the joint likelihood function $p(X^n,S^n|\,\theta)$ with its weighted version $\wt p(X^n,S^n|\,\theta)$. Similarly to the CAVI, we can repeatedly update each $q_{\theta_j}$ and $q_{S_i}$ until convergence, and a stopping criterion can be based on the change in the weighted ELBO. Again, the weighted CAVI reduces the the CAVI when are weights are identically one. Due to the similarity between~\eqref{Eqn:CAVI} and the preceding display, only minor changes are needed in order to implement the weighted CAVI based on existing statistical softwares for the CAVI. Algorithm~\ref{algo:vbB} below summarizes the pseudo-code for the weighted CAVI.
\begin{figure}[ht]
  \centering
  \begin{minipage}{.9\linewidth}
\begin{algorithm}[H]
\SetAlgoLined
\KwIn{Random weights $W^{(b)}_1,W^{(b)}_2,\ldots,W^{(b)}_n$}
\KwData{$X_1,X_2, \ldots, X_n$}
\BlankLine
\emph{Initialize $q_{Z^n}=\bigotimes_{j=1}^d q_{\theta_j}\otimes \bigotimes_{i=1}^n q_{S_i}$}\;
\While{not converged}{
\For{$j \leftarrow 1$ \KwTo $d$}{
\emph{Update $q_{\theta_j}$ as in equation~\eqref{eq:wcavig}} \;
}
\For{$i \leftarrow 1$ \KwTo $n$}{
\emph{Update $q_{S_i}$ as in Equation~\eqref{eq:wcavil}} \;
}
}
\KwOut{Variational approximation $\wt q_{\theta}^{(b)}=\bigotimes_{j=1}^d q_{\theta_j}$ for the weighted posterior $\wt \pi_n^{(b)}$ of $\theta$}
\caption{Weighted CAVI Algorithm}
\label{algo:vbB}
\end{algorithm}
\end{minipage}
\end{figure}
\\


\section{Non-asymptotic analysis of variational weighted likelihood Bootstrap}\label{sec:mainres}
In this section, we develop theoretical justifications for the statistical inference procedure developed in Section~\ref{Sec:Inference_MF}. Our results show that unlike the mean-field variational approximation $\wht q_\theta$ to the marginal posterior $\pi_n$ of $\theta$ that generally underestimates the variance (c.f.~Theorem~\ref{thm:1}), the variational weighted likelihood Bootstrap generates independent random samples of $\theta$ given the data whose distribution approaches $\pi_n$ as $n\to\infty$. As a consequence, credible intervals constructed from these random samples of $\theta$ has frequentist coverage approaching to their nominal levels as $n\to\infty$. Our analysis is non-asymptotic and leads to explicit high probability error bounds on the discrepancies.

\subsection{Contraction of weighted posterior distribution and its normal approximation}\label{Sec:contration_WMF}
In this subsection, we investigate the theoretical properties of the weighted variational approximation $\wt q^{(b)}_\theta$ as the optimum of the weighted variational optimization~\eqref{eqn:vwlbObj}. 
Similar to the study of the mean-field approximation $\wht q_{\theta}$ in  Section~\ref{Sec:contration_MF}, we begin with he contraction property of the weighted posterior distribution $\wt \pi_n^{(b)}$ defined in~\eqref{eq:weightedPo} that leads to the contraction of $\wt q^{(b)}_\theta$ with a sub-Gaussian type tail bound. It proof is provided in Section~\ref{app:Proof_thm_WMF_contraction}.
\begin{theorem}\label{thm:WMF_contraction}
Under Assumptions A1, A2 and W, there exist constants $(C'_0,C'_1,C'_2,C'_3)$ such that for any $M\geq 1$ and $\wt \varepsilon_n=C'_0 \frac{\log n}{\sqrt{n}}$ it holds with probability at least $1-C'_1M^{-2}$ that the weighted posterior $\wt\Pi^{(b)}_n$ satisfies
\begin{align}\label{Eqn:contraction_WP}
\wt \Pi^{(b)}_n \big( \|\theta - \theta^\ast\|\geq C'_2 \varepsilon\big) \leq e^{-C'_3n\varepsilon^2},\quad\mbox{for all } \ \varepsilon\geq M\wt\varepsilon_n.
\end{align}
In addition, if Assumption A3 is also true, then there exist some constants $(C''_2,C''_3)$ such that under the same high probability event, the mean-field approximation $\widehat Q_{\theta}$ satisfies
\begin{align*}
\wt Q_{\theta}^{(b)} \big( \|\theta - \theta^\ast\|\geq C''_2 \varepsilon\big) \leq e^{-C''_3n\varepsilon^2},\quad\mbox{for all } \ \varepsilon\geq M\wt\varepsilon_n.
\end{align*}
\end{theorem}
\noindent The proof of~\eqref{Eqn:contraction_WP} involves a uniform control on the weighted likelihood ratio via the bracket entropy, and uses the proof technique of Theorem 2 in~\cite{RatesPosterior}. 
Similar to Lemma~\ref{lem:expDacayverify}, the proof of Theorem~\ref{thm:WMF_contraction} is not specific to parametric models where $\wt\varepsilon_n\asymp \sqrt{\log n/n}$ and can be extended to general cases such as infinite-dimensional models as long as the bracket entropy of the model space is properly controlled. 

Recall that $\wt\theta^{(b)}$ is the maximizer of the weighted likelihood function,
\begin{align*}
\wt\theta^{(b)}=\argmax_{\theta\in\Theta}\widetilde L^{(b)}(\theta;\, X^n)=\argmax_{\theta\in\Theta}\sum_{i=1}^n W_i^{(b)}\log p(X_i\,|\,\theta).
\end{align*}
Let $\wt Q^{\ast(b)}$ denote the normal distribution $N\big(\wt\theta^{(b)},\,[nI(\theta^\ast)]^{-1}\big)$.
The next theorem extends the classical BvM theorem to the normal approximation $\wt Q^{\ast(b)}$ of the weighted posterior distribution $\wt\Pi_n^{(b)}$ of $\theta$.
\begin{theorem}\label{thm:BvM_wp}
Under Assumptions A1, A2 and W, there exist constants $(C'_4, C'_5)$ such that for any $M\geq 1$, it holds with probability at least $1-C'_4M^{-2}$ that
\begin{align*}
\max\bigg\{d_{TV}\big(\wt\Pi_n^{(b)}(A), \wt Q^{\ast(b)}(A)\big),\ D\big( \wt\Pi_n^{(b)}\,||\, \wt Q^{\ast(b)}\big)\bigg\} \leq \frac{C'_5 M (\log n)^{d+4}}{\sqrt{n}}.
\end{align*}
\end{theorem}
\noindent Theorem~\ref{thm:BvM_wp} is a non-asymptotic result providing an explicit upper bound to certain discrepancy measures between $\wt\Pi_n^{(b)}$ and its normal approximation $\wt Q^{\ast(b)}$. This theorem also extends the classical BvM type results from the total variational metric to the stronger KL-divergence, due to the strong sub-Gaussian tail bound~\eqref{Eqn:contraction_WP} for controlling the expectation of the 
log-density ratio $\log(\wt\pi_n^{(b)}/\wt q^{\ast(b)})$ relative to $\wt\Pi_n^{(b)}$ outside a $\sqrt{\log n/n}$-neighborhood of $\wt \theta^{(b)}$.

\subsection{Consistency of variational weighted likelihood Bootstrap}
In this subsection, we discuss the consistency of variational weighted likelihood Bootstrap summarized in Algorithm~\ref{algo:vbA} as a new sampling scheme for approximating the marginal posterior distribution $\Pi_n$ of $\theta$.
First, we present a result on the normal approximation of the weighted mean-field approximation $\wt Q^{(b)}_\theta$. 
\begin{theorem}\label{thm:3}
Under Assumptions A1, A2, A3 and W, there exist constants $(C'_6,C'_7)$ such that for any $M\geq 1$, it holds with probability at least $1-C'_6M^{-2}$ that
\begin{equation*}
D( \wt Q^{(b)}_\theta\,||\, \wt{Q}^{\ast{(b)}}_{VB} ) \leq  \frac{C'_7M^3(\log n)^{d + 4}}{\sqrt{n}},
\end{equation*}
where $\wt{Q}^{\ast{(b)}}_{VB} = N(\wt{\theta}^{(b)}, [nI_{VB}]^{-1})$, with the matrix $I_{VB}\in\mb R^{d\times d}$ defined in Theorem~\ref{thm:1}.
\end{theorem}
\noindent This result shows that the weighted version $\wt{Q}^{\ast{(b)}}_{VB}$ shares the same covariance structure as $\wt{Q}^{\ast}_{VB}$ in Theorem~\ref{thm:1}, but the center changes from the MLE $\wht\theta$ to the weighted MLE $\wt\theta^{(b)}$. Since the diagonal matrix $I_{VB}$ not only ignores all off-diagonal entries in the information matrix $I(\theta^\ast)$ but also inflates the diagonals by an extra additive term $I_s(\theta^\ast)$ due to the mean-field approximation on the latent variables, statistical inference based on $\wt Q^{(b)}_\theta$ will be erroneous. Fortunately, as the center of $\wt Q^{(b)}_\theta$ approximates the weighted MLE $\wt\theta^{(b)}$, we may instead conduct inference based on this quantity.
Formally, recall that $\wt\theta^{(b)}_{VB} =\int_\Theta \theta\,\wt q_{\theta}^{(b)}(\theta)\,\dd\theta$ is the weighted variational mean of $Q_{\theta}^{(b)}$ in the $b$th replicate of Algorithm~\ref{algo:vbA}. 
\begin{corollary}\label{cor:WVB_center}
Under the conditions and the high-probability event of Theorem~\ref{thm:3}, there exists a constant $C'_8$ such that
\begin{align*}
\|\sqrt{n}\,( \wt\theta^{(b)}_{VB}-\wt{\theta}^{(b)})\| \leq \frac{C_6M^{3/2}(\log n)^{d/2 + 2}}{n^{1/4}}.
\end{align*}
\end{corollary}
\noindent This corollary is the weighted version of Corollary~\ref{cor:VB_center} that indicates the closeness between the weighted variational mean of $Q_{\theta}^{(b)}$ and the weighted MLE $\wt\theta^{(b)}$, whose conditional distribution given $X^n$ approximates the sampling distribution of the MLE $\mle$.

As a consequence of Corollary~\ref{cor:WVB_center}, the following theorem provides a theoretical justification of using the random samples $\{\wt\theta^{(b)}_{VB}\}_{b=1}^B$ for approximating the posterior $\Pi_n$ of $\theta$. For two vectors $u,v\in\mb R^d$, we use $u\leq v$ to mean that $u$ is element-wise less than or equal to $v$. We use $U^{1/2}$ to denote the matrix square root of a positive definite matrix $U\in\mb R^{d\times d}$.
\begin{theorem}\label{thm:VWLB_consistency}
Under the assumptions of Theorem~\ref{thm:3}, there exist constants $(C'_9,C'_{10},C'_{11})$ such that for any $M\geq 1$, it holds with probability at least $1-C'_9M^{-2}$ that
\begin{align*}
\sup_{u\in\mb R^d}\Big|P \big(\sqrt{n}\,[I(\theta^\ast)]^{1/2} \big(\wt\theta^{(b)}_{VB}-\mle\,\big) \leq u\,\big|\,X^n\big) - P(Z \leq u) \Big|& \leq \frac{C'_{10}M^{3/2}(\log n)^{d/2+2}}{n^{1/4}},\\
\quad\mx{and}\qquad \sup_{u\in\mb R^d}\Big|P \big(\wt\theta^{(b)}_{VB} \leq u\,\big|\,X^n\big) - \Pi\big(\theta \leq u\,\big|\,X^n\big) \Big|& \leq \frac{C'_{11}M^{3/2}(\log n)^{d/2+2}}{n^{1/4}},
\end{align*}
where $Z\sim N(0, I_d)$ is the $d$-variate standard normal distribution, and the randomness of conditional probability $P(\cdot\,|\,X^n)$ is on the random weights $(W^{(b)}_1,W^{(b)}_2,\ldots,W^{(b)}_n)$.
\end{theorem}
\noindent The first display in Theorem~\ref{thm:VWLB_consistency} implies that the conditional cdf of $\sqrt{n}\,[I(\theta^\ast)]^{1/2} \big(\wt\theta^{(b)}_{VB}-\mle\,\big)$ given $X^n$ uniformly converges to that of $N(0, I_d)$ as $n\to\infty$, and the second display implies the uniform convergence of the conditional cdf of $\wt\theta^{(b)}_{VB}$ to that of the marginal posterior distribution $\Pi(\cdot\,|\,X^n)$ of $\theta$. As a direct consequence of Theorem~\ref{thm:VWLB_consistency} and the classical BvM result on the posterior $\Pi(\cdot\,|\,X^n)$, we may use sample quantiles of $\{\wt\theta^{(b)}_{VB}\}_{b=1}^B$ to construct a credible interval, whose frequentist coverage is at most $\m O(\sqrt{\log n}/{n^{1/4}})$ away from its nominal level for sufficiently large $B$.

\section{Numerical study}\label{sec:simulation}
In this section, we provide two numerical examples, the Gaussian mixture model and the Bayesian linear regression, to evaluate the performance of the variational weighted likelihood Bootstrap for credible interval constructions. 
We will compare four types of credible intervals with level $\alpha=95\%$. The first interval is based on the sample quantiles of the draws from the Gibbs sampler for sampling from the posterior $\Pi_{n}$. The second interval is based on the quantile of the mean-field variational approximation $\wht Q_\theta$. The third interval is based on the sample quantile of the draws $\{\wt\theta_{VB}^{(b)}\}_{b=1}^B$ from our VWLB procedure (Algorithm~\ref{algo:vbA}). The last interval is the usual bootstrap interval based on the sample quantile of $\{2\wt\theta_{VB}^{(b)}-\wht\theta_{VB}\}_{b=1}^B$, since by Corollary~\ref{cor:VB_center} and Theorem~\ref{thm:VWLB_consistency}, the conditional distribution of $\sqrt{n}(\wt\theta_{VB}^{(b)}-\wht\theta_{VB})$ given $X^n$ provides a good approximation to the sampling distribution of $\sqrt{n}(\wht\theta_{VB}-\theta^\ast)$. Note that the last two types of intervals are asymptotically equivalent since the limiting distribution of $\sqrt{n}(\wht\theta_{VB}-\theta^\ast)$ is symmetric.

\subsection{Gaussian Mixture Model}\label{Sec:GMM}
Gaussian Mixture Model(GMM) is a classical example of latent variable models. We consider the following one-dimensional GMM in this simulation. Assume $(X_1,X_2,\ldots,X_n)$ as i.i.d.~sample from the data generating model $P_{\theta} = \sum_{k = 1}^K \frac{1}{K} N(\mu_k,1)$, with $K=3$, where the parameter is the centers $\theta =\mu=(\mu_1,\ldots,\mu_K)$. This model has a latent variable representation by associating each $X_i$ with a latent assignment (variable) $c_i$ that follows the categorical distribution:
$c_i \sim \text{Categorical}(1/K, \ldots, 1/K)$. Conditioning on $c_i$ and $\mu$, the observation $X_i$ follows the normal distribution $N(c_i^T\mu,1)$. 
We specify the prior distribution of the parameter $\mu = (\mu_1,\mu_2,\mu_3)$ as i.i.d.~$\mu_k \sim N(0,\sigma^2)$ with $\sigma=5$. The true parameter $\mu^\ast=(-\Delta,0,\Delta)$, with $\Delta$ ranging from $0$ to $8$.

We apply the mean-field approximation that approximates the joint posterior distribution $p(\{c_i\}_{i=1}^n,\mu\,|\,X^n)$ with a fully factorized distribution $q_{Z^n}=\bigotimes_{k=1}^K q_{\mu_k}\otimes \bigotimes_{i=1}^nq_{c_i}$. It turns out that the normal prior of $\mu$ is ``conjugate" in the sense that the KL minimizer $\wht q_{Z^n}=\bigotimes_{k=1}^K\wht q_{\mu_k}\otimes \bigotimes_{i=1}^n\wht q_{c_i}$ as in the optimization problem~\eqref{eqn:vbObj} must be from the same distribution family: each $\wht q_{\mu_k}$ is a normal distribution, and each $\wht q_{c_i}$ is a categorical distribution. Therefore, we can parametrize them by
\begin{align}\label{Eqn:GMM_MF}
Q_{\mu_k}=N(\mu_k; m_k, s_k^2),\ \ \mx{and}\ \ Q_{c_i} = \text{Categorical}(\phi_{i1}, \ldots, \phi_{iK}),
\end{align}
for $k=1,2,3$ and $i=1,\ldots,n$.
Solving the infinite-dimensional optimization problem~\eqref{eqn:vbObj} boils down to optimizing the objective function over these parameters $\{m_k,s_k\}_{k=1}^3$ and $\{(\phi_{i1},\ldots,\phi_{iK})\}_{i=1}^n$.
Implementation details of the algorithm is provided in Appendix~\ref{app:GMM}.

\begin{figure}[ht!]
\centering
\includegraphics[scale=0.55]{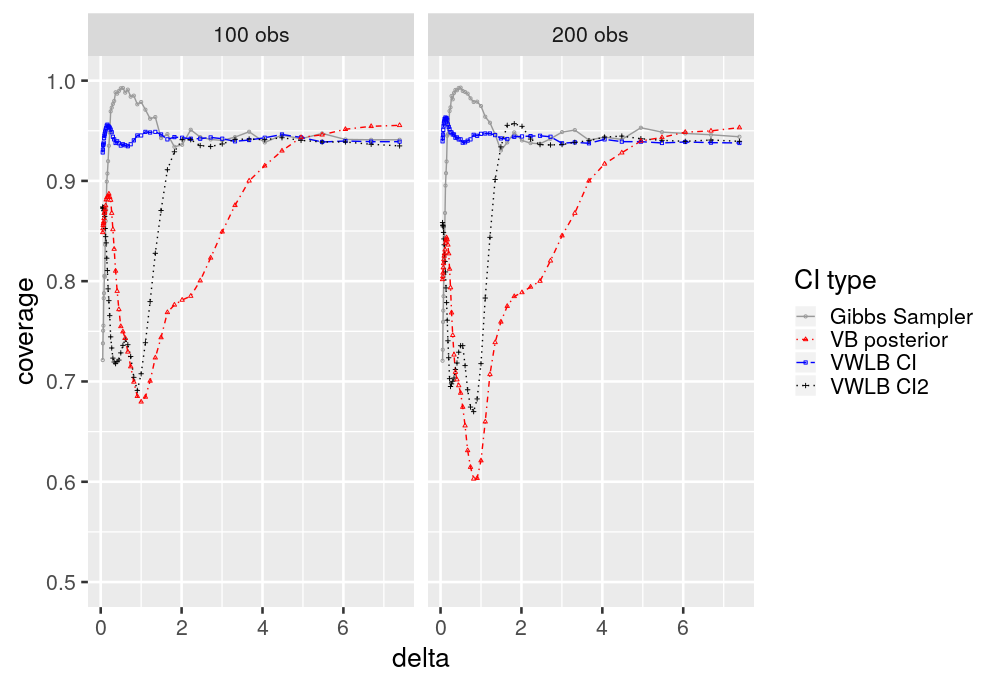} \\
\includegraphics[scale=0.55]{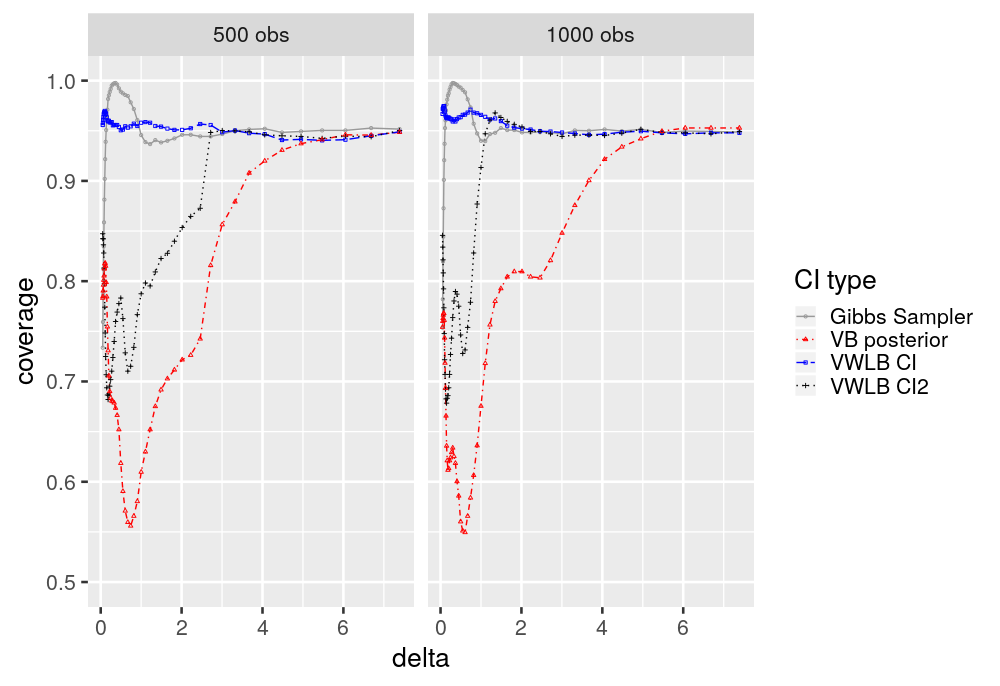} 
\caption{Coverage probabilities of four credible intervals (CI) in the GMM under $n\in\{100,200,500,1000\}$. Gibbs sampler: CI based on Gibbs sampling; VB posterior: CI based on the mean-field approximation $\wht Q_\theta$; VWLB CI: CI based on the VWLB samples $\{\wt\theta_{VB}^{(b)}\}_{b=1}^B$; VWLB CI2: bootstrap type CI based on $\{2\wt\theta_{VB}^{(b)}-\wht\theta_{VB}\}_{b=1}^B$.}\label{fig:1}
\end{figure}

We compare the numeric performance of the four credible intervals. In the Gibbs sampler, we take out the first $10,000$ iterations as the burn-in, and fetch $500$ samples every $20$ iterations to reduce the auto-correlation. In our VWLB, we set the number of draws to be $B = 500$. Figure~\ref{fig:1} reports the averaged coverage probability over three parameters $(\mu_1,\mu_2,\mu_3\}$ in each type of intervals, under sample size $n=100$, $200$, $500$ and $1000$ respectively, and Table~\ref{Table:1} reports the credible interval lengths under $\Delta\in\{0, 1, 3, 5\}$ and $n=500$ (clusters $1$ and $3$ are symmetric, so we only report clusters $1$ and $2$). In Appendix~\ref{app:GMM}, we provide more details about the individual coverage probability for each of the three normal centers, and the respectively estimated posterior density curves. 
Interestingly, when $\Delta$ is near $0$, all methods except for the VWLB have degeneracy problem, as the three clusters are no longer distinguishable (our regularity Assumptions A2 and A3 are violated). As we can see, the degeneracy window width decreases as the sample size grows. In the good region where the cluster gap $\Delta$ is sufficiently large so that the normal components are statistically distinguishable and our theory applies, the two credible intervals based on VWLB tends to attain the nominal $95\%$ level as the Gibbs sampler. Moreover, the lengths of the credible intervals based on Gibbs and VWLB are similar. In contrast, the credible intervals directly constructed from the mean-field approximation $\wht Q_\theta$ are shorter than those from the Gibbs sampling, and tend to under-estimate the uncertainty until the gap $\Delta$ exceeds 5 where the three normal components become nearly disjoint across all sample sizes. 

\begin{table}[ht!]
\resizebox{\columnwidth}{!}{
\begin{tabular}
{l|ll|ll|ll|ll}
&cluster $1$ &cluster $2$
&cluster $1$ &cluster $2$
&cluster $1$ &cluster $2$
&cluster $1$ &cluster $2$ \\\hline
Gibbs Sampler  
&$0.490(0.053)$ &$0.427(0.068)$
&$0.537(0.056)$ &$0.848(0.104)$
&$0.372(0.019)$ &$0.490(0.044)$
&$0.312(0.014)$ &$0.325(0.016)$ \\
VB Posterior 
&$0.303(7 \times 10^{-5})$ &$0.303(4 \times 10^{-5})$
&$0.303(0.001)$ &$0.304(0.001)$
&$0.303(0.007)$ &$0.303(0.007)$
&$0.303(0.009)$ &$0.303(0.009)$ \\
VWLB Sampler  
&$0.472(0.132)$ &$0.354(0.118)$
&$0.581(0.040)$ &$0.857(0.117)$
&$0.375(0.039)$ &$0.499(0.088)$
&$0.310(0.024)$ &$0.322(0.027)$ 
\end{tabular}
}
\caption{Credible interval lengths for $\mu_1$ and $\mu_2$ under $n=500$ and $\Delta\in\{0,1,3,5\}$.}\label{Table:1}
\end{table}

\begin{figure}[ht!]
\centering
\subfloat[Variational variance of $\mu_1$ (with error bars).]{{\includegraphics[scale=0.25]{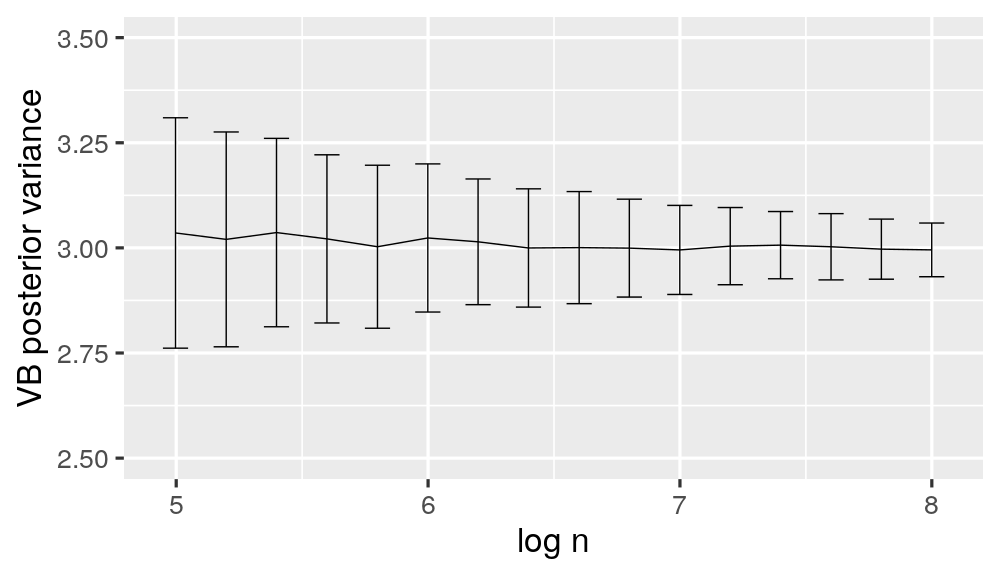} }}%
\qquad
\subfloat[Ratio of variances]{{\includegraphics[scale=0.25]{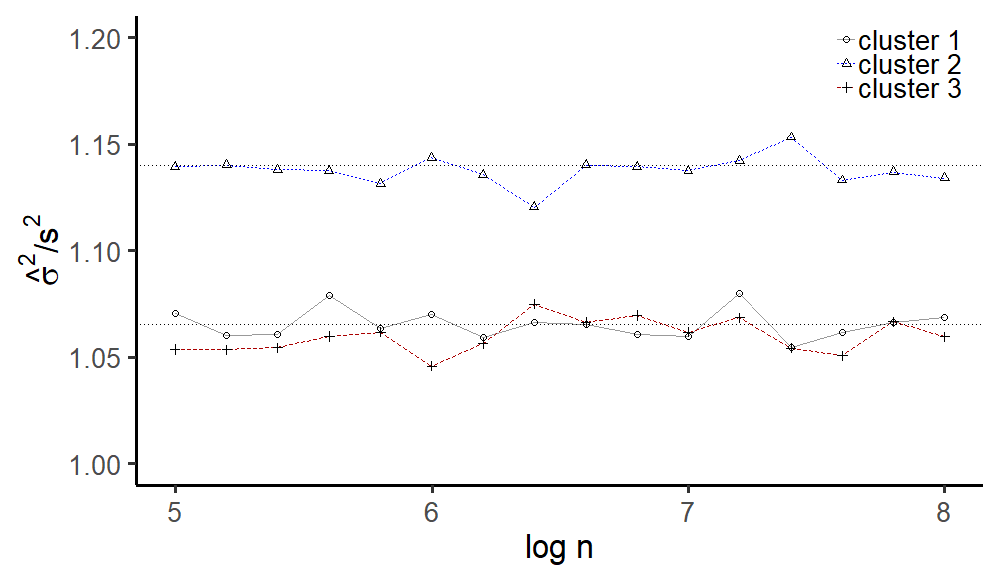} }}%
\caption{Characterizations of the limiting variances in the mean-field approximation.}%
\label{fig:gmmVBinfo}
\end{figure}

A direct calculation indicates that the mean-field information matrix $I_{\mathrm{VB}}(\theta^\ast)$ appeared in the normal approximation in Theorem~\ref{thm:1} is simply the diagonal matrix with each diagonal component being equal to $1/K$, regardless the relative positions between the components of $\mu$. Therefore, Theorem~\ref{thm:1} implies that the limiting variational variance of each $\mu_k$ is $K/n$. Figure~\ref{fig:gmmVBinfo}(a) shows the Monte Carlo estimation of the rescaled variance (multiply by $n$) of $\mu_1$ in the mean-field approximation $\wht Q_\theta$ under $\Delta=3$, which empirically justified our theoretical prediction. In addition, the second row in Table~\ref{Table:1} shows that the length of the credible interval from the variational approximation $\wht Q_\theta$ remains $0.303$ across all $\Delta$, where $0.303\approx 2\times 1.96 \times \sqrt{K/n}$ under $K=3$ and $n=500$, which is consistent with our theoretical prediction.  When $\Delta$ is near $0$, the Fisher information matrix $I(\theta^\ast) \approx I(0) = (1/K^2)1_K1_K^T$ becomes degenerate, which explains the under-covering issue associated with the mean-field credible interval. As $\Delta$ increases, $(\mu_1,\mu_2,\mu_3)$ tends to become independent in the joint posterior $\Pi_n$, and the credible intervals from the mean-field $\wht Q_\theta$ approaches the nominal level $95\%$ since $I(\theta^\ast)$ converges to $I_{\mathrm{VB}}(\theta^\ast)$ (c.f.Table~\ref{Table:1}).

We also verify our theoretical finding in Theorem~\ref{thm:1} that due to the presence of latent variables, the mean-field information matrix $I_{\mathrm{VB}}(\theta^\ast)=\mx{Diag}(I(\theta^\ast)+I_s(\theta^\ast))$ has an extra term caused by the missing data information matrix $I_s(\theta^\ast)$.  For each $k=1,2,3$, according to the classical Bernstein von-Mises theorem, we can estimate the $k$-th diagonal element of the Fisher information matrix $I(\theta^\ast)$ as $n^{-1}$ times the inverse of $\hat{\sigma}^2_k = [(\hat{\Sigma})^{-1}]_{kk}$, where $\hat{\Sigma}$ denotes the sample covariance matrix of the draws from the posterior distribution $\Pi_n$ via Gibbs sampling. The diagonal of $I_{VB}(\theta^\ast)$ can be approximated by $n^{-1}$ times the inverse of the variance $s_k^2$ in the variational approximation~\eqref{Eqn:GMM_MF}. Therefore, $\hat{\sigma}^2_k/s_k^2$ serves as an estimate of the ratio $[I_{\mathrm{VB}}(\theta^\ast)]_{kk}/[I(\theta^\ast)]_{kk}$ for $k=1,2,3$. Figure~\ref{fig:gmmVBinfo}(b) plots the $\hat{\sigma}^2_k/s_k^2$ versus the sample size for the three clusters $k=1,2,3$ under $\Delta=5$, as well as their respective theoretical value calculated from Theorem~\ref{thm:1}. As we can see, the numerical results closely conform to our theoretical prediction, where all these three ratios are greater than one due to the latent variables.

\subsection{Bayesian Linear Regression}\label{Sec:BLR}
In this example we consider the Bayesian linear regression of the following form:
\begin{equation*}
Y^n = X\beta + \varepsilon, \quad \varepsilon \sim N(0, \sigma^2 I_n).
\end{equation*}
where $Y^n$ is the response vector in $\mb R^n$, and $X=(X_1,\ldots,X_p)$ is the design matrix of dimension $n \times p$. We target at making inference on $\beta=(\beta_1,\ldots,\beta_p)\in\mb R^p$, the regression coefficient vector. To facilitate variable selection, we impose independent point mass mixture priors for $\{\beta_j\}_{j=1}^p$ by introducing $p$ latent binary variables $\{\gamma_j\}_{j=1}^p$, where $\gamma_j$ indicates whether $X_j$ is included in the model (or $\beta_j=0$ or not), that is,
\begin{align*}
\beta_j\, | \, (\gamma_j = 0,\sigma^2) \sim \delta_0 \quad\mx{and}\quad \beta_j | (\gamma_j = 1,\sigma^2) \sim N(0,v_1 \sigma^2),
\end{align*}
where $\delta_a$ denotes the point mass measure at point $a$.
The prior for $\gamma_j$'s is specified as i.i.d.~Bernoulli$(\xi)$, where $\xi$ follows the Beta prior distribution as Beta$(a_0, b_0)$, and the prior for the noise variance $\sigma^2$ is the inverse Gamma distribution as IG$(\nu/2, \nu\lambda/2)$. In our simulation, we take $v_1=2$, $a_0=1$, $b_0= 1$, $\nu= 0.002$ and $\lambda= 1$ for the hyperparameters. In this example, there is no latent variable and the parameter is $\theta=(\{(\gamma_j,\beta_j)\}_{j=1}^p,\sigma^2)$.  

We apply the following block mean-field approximation~\cite{carbonetto2012scalable} for approximating the joint posterior distribution $p(\{(\gamma_j,\beta_j)\}_{j=1}^p,\sigma^2\,|\,Y^n)$ using the blockwise-factorized family $q_{\theta}=\bigotimes_{j=1}^p q_{\gamma_j,\beta_j} \otimes q_{\sigma^2}$. It again turns out that the point mass mixture priors on $\beta$ and inverse gamma prior on $\sigma$ is the ``conjugate" prior in the sense that in the KL minimizer $\wht q_{\theta}=\bigotimes_{j=1}^p \wht q_{\gamma_j,\beta_j} \otimes \wht q_{\sigma^2}$, each $\wht q_{\gamma_j,\beta_j}$ is a point mass mixture and $\wht q_{\sigma^2}$ is an inverse Gamma. Therefore, we can parametrize them by
\begin{align}\label{Eqn:MF_block}
Q_{\gamma_j,\beta_j} = (1-\phi_j) \delta_0(\gamma_j)\otimes\delta_0(\beta_j) + \phi_j \delta_1(\gamma_j)\otimes N(\mu_j,\sigma_j^2), \ \ \mx{and} \ \ Q_{\sigma^2} = \mx{IG}(c,d),
\end{align}
for $j=1,\ldots,p$.  Therefore, finding $\wht Q_{\theta}$ amounts to optimizing the objective function~\eqref{eqn:vbObj} over these parameters $\{\phi_j,\mu_j,\sigma_j^2)\}_{j=1}^p$ and $\{c,d\}$ via the coordinate descent algorithm.
Implementation details of the algorithm is provided in Appendix~\ref{app:BLM}, where we have adopted a slightly different but more efficient algorithm~\cite{FastBLR}.

\begin{figure}[ht!]
\centering
\subfloat[$\beta_1$]{\includegraphics[scale=0.25]{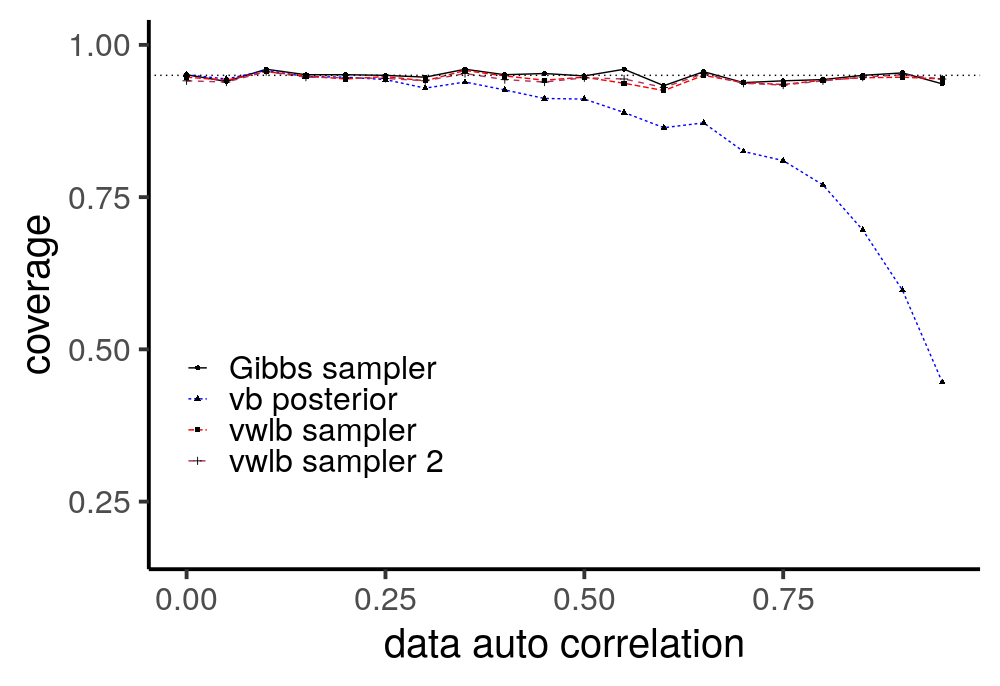}}
\qquad
\subfloat[$\beta_4$]{\includegraphics[scale=0.25]{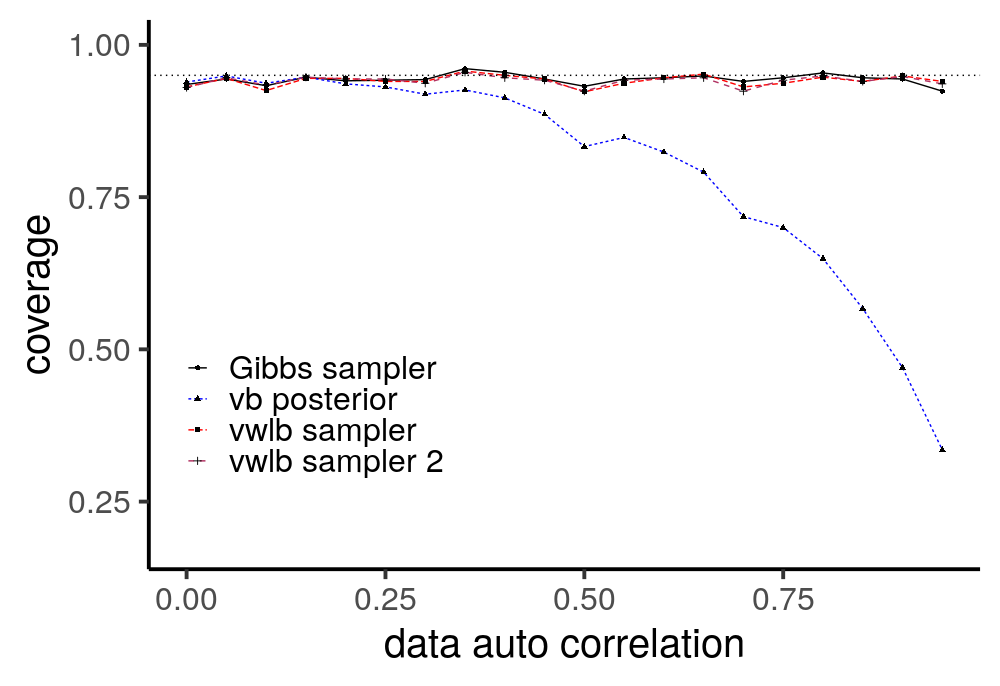}}
\caption{Coverage probabilities based on four types of credible intervals versus design correlation $\rho$.}\label{fig:BLM}
\end{figure}

We generate the data in our simulation as follows. The number of observations is $n = 1,000$ and number of covariates $p = 10$, with ground truth parameter $\beta^\ast = (2,3,2,4,1, 2,1, 0,0,2)^T$. We generate $X$ conforming to an AR(1) process associated with a unit variance white noise process as below: for each $i=1,\ldots,n$,
\begin{equation*}
X_{ij} = \rho X_{i(j-1)} + \varepsilon_j, \ \ \varepsilon_j \sim N(0,1)\ \ \mx{for $j=2,\ldots,p$, \ \ and} \ \ X_{i1}\sim N(0,(1-\rho^2)^{-1}),
\end{equation*}
where $\rho\in[0,1)$ is the auto-correlation, and each $X_{ij}$ has the same marginal distribution as $N(0,(1-\rho^2)^{-1})$.
In the simulation, we consider different settings of the auto-correlation as $\rho = 0, 0.05, \ldots, 0.95$. We draw $B = 1000$ samples from the Gibbs sampler and the VWLB sampler (Algorithm~\ref{algo:vbA}). When evaluating four credible intervals, we compute the coverage probability based on $1,000$ replicates for each setting. Figure~\ref{fig:BLM} reports the trends of coverage probabilities for two coefficients $\beta_1$ and $\beta_4$ using the four aforementioned types of credible intervals. 
As we can expect, the coverage probabilities of the credible intervals based on the mean-field approximation $\wht Q_{\theta}$ rapidly fall below the nominal level $95\%$ when the auto-correlation is large, as $\wht Q_{\theta}$ completely ignores the (high) correlations among $\{\beta_j\}_{j=1}^p$ in the joint posterior distribution $\Pi_n$. In comparison, the rest three methods exhibit similar patterns and nearly attain the nominal level (horizontal dotted line) across all $\rho$ values.
We also report the lengths of the intervals reflecting the estimated uncertainty magnitudes in Table~\ref{Table:2}.
\begin{table}[ht!]
\resizebox{\columnwidth}{!}{
\begin{tabular}{l|ll|ll|ll}
& $\beta_1$ & $\beta_4$
& $\beta_1$ & $\beta_4$
& $\beta_1$ & $\beta_4$ \\\hline
Gibbs Sampler &$0.1245(0.0066)$  &$0.1243(0.0064)$   
&$0.1243(0.0064)$ &$0.1388(0.0071)$ & $0.1222(0.0099)$ &$0.1675(0.0174)$    \\
VB Posterior  &$0.1241(0.0028)$ &$0.1241(0.0029)$
&$0.1074(0.0025)$ &$0.1075(0.0025)$ & $0.0387(0.0010)$  &$0.0388(0.0010)$ \\
VWLB Sampler  &$0.1223(0.0070)$  &$0.1226(0.0070)$
&$0.1223(0.0071)$ &$0.1368(0.0077)$  & $0.1218(0.0071)$ &$0.1680(0.0097)$ 
\end{tabular}
}
\caption{Credible interval lengths for $\beta_1$ and $\beta_4$ under $n=1,000$ and $\rho\in\{0,0.5,0.95\}$.}\label{Table:2}
\end{table}
As we can infer from this table, the degree of uncertainty underestimation in the mean-field approximation $\wht Q_\theta$ increases as the correlation among $\{\beta_j\}_{j=1}^p$ in their joint posterior increases, as the posterior covariance matrix of the regression coefficient vector $\beta$ is approximately proportional to the auto-covariance matrix of the AR(1) process with auto-correlation $\rho$. For example, when $\rho=0$, there is no visible uncertainty underestimation, while when $\rho=0.95$, the variational standard deviation of $\beta_4$ reduces to roughly one quarter of the true (marginal) posterior standard deviation. In comparison, the lengths from Gibbs sampler and our VWLB sampler are always close.


\section{Proofs of the main results}\label{sec:proofs}
In this section, we provide a selective proofs of the main results in the paper, and leave the rest and some technique results to the supplement. To simplify the presentation, we use letter $C$ to denote a generic constant whose value may change from one line to another throughout the proof.

\subsection{Proof of Lemma~\ref{lem:expDacayverify}}\label{Sec:Proof_lem:expDacayverify}
Before showing that the variational approximation $\widehat Q_\theta$ has the desired sub-Gaussian type tail bound, we first show that the marginal posterior distribution $\Pi_n$ of $\theta$ has a similar type bound.
In fact, Assumption A2 implies $E_{\theta^\ast}\log \frac{p(X\,|\,\theta^\ast)}{p(X\,|\,\theta)} \leq C \|\theta-\theta^\ast\|^2$ and $E_{\theta^\ast}\big[ \log\frac{p(X\,|\,\theta^\ast)}{p(X\,|\,\theta)} \big]^2 \leq C \|\theta-\theta^\ast\|^2$. Therefore, we can apply Theorem 5.1 in~\cite{ghosal2000convergence} to obtain that for any $M\geq 1$, it holds with probability at least $1-CM^{-2}$ that the posterior $\Pi_n$ of $\theta$ satisfies 
\begin{align}\label{Eqn:post_contraction}
\Pi_n\big(\|\theta-\theta^\ast\|\geq C\varepsilon\big) \leq e^{-Cn\varepsilon^2}, \quad\mx{for all $\varepsilon\geq M\varepsilon_n$},
\end{align}
where we slightly strengthen their result by providing the explicit posterior tail bound and by extending it from a single $\varepsilon=M\varepsilon_n$ to all $\varepsilon\geq M\varepsilon_n$ simultaneously, due to the same argument as the proof of equation~(6.10) in \cite{yang2016bayesian}. 

We will use the optimality of $\wht Q = \wht Q_{Z^n}=\wht Q_{\theta} \otimes \wht Q_{S^n}$ for optimization problem~\eqref{eqn:vbObj} to prove the desired result. Let $Q|_A$ denote the restriction of a probability measure $Q$ onto a set $A\subset \Theta\times \m S^n$, that is, $Q|_A(B)=Q(A\cap B)/Q(A)$ for all $B\subset \Theta\times \m S^n$. For a fixed $\varepsilon \geq M\varepsilon_n$ and each $j=1,2,\ldots,d$, we construct a sequence $\m Q= \{Q^\dagger_{\lambda}:\,\lambda\in[0,1]\}$ of distributions as
\begin{align*}
Q^\dagger_\lambda = (1-\lambda)\, \wht Q|_{A_{n,j}^c} + \lambda\, \wht Q|_{A_{n,j}}, \quad\mx{for all $\lambda\in[0,1]$,}
\end{align*}
where $A_{n,j}=\{\theta:\,|\theta_j-\theta^\ast_j|\geq D\varepsilon\} \times \m S^n$ for some sufficiently large constant $D\geq C$. It is easy to verify that $Q^\dagger_\lambda$ is a valid distribution belonging to the mean-field family $\Gamma$ for any $\lambda\in[0,1]$. Let $\wht\lambda =\wht Q(A_{n,j})$, so that $\wht Q = Q^\dagger_{\wht\lambda}$. Due to the optimality of $\wht Q$ for minimizing $D\big(Q\,||\,P(\cdot\,|\,X^n)\big)$ when restricting to the smaller family $\m Q$, we obtain $\wht\lambda =\argmin_{\lambda\in[0,1]} D\big(Q^\dagger_\lambda\,||\,P(\cdot\,|\,X^n)\big)$, where we can express
\begin{align*}
&D\big(Q_\lambda^\dagger\, \| \,P(\cdot\,|\,X^n)\big) = \int_{\Theta\times \m S^n} \log(\frac{Q^\dagger_{\lambda}(\dd Z^n)}{ P(\dd Z^n\,|\,X^n)})\,\big[(1 - \lambda)\,\wht{Q}|_{A_{n,j}^c}(\dd Z^n) + \lambda\, \wht{Q}|_{A_{n,j}}(\dd Z^n)\big]  \\
&\quad = (1 - \lambda)\int_{A_{n,j}^c}  \log(\frac{(1 - \lambda)\,\wht{Q}|_{A_{n,j}^c}(\dd Z^n)}{P(\dd Z^n\,|\,X^n)})\,\wht{Q}|_{A_{n,j}^c}(\dd Z^n) + \lambda \int_{A_{n,j}} \log(\frac{\lambda\,\wht{Q}|_{A_{n,j}}(\dd Z^n)}{P(\dd Z^n\,|\,X^n)})\,\wht{Q}|_{A_{n,j}}(\dd Z^n),
\end{align*}
where the second equality is due to mutual orthogonality of $\wht Q|_{A_{n,j}}$ and $\wht Q|_{A_{n,j}^c}$. 
Due to a similar decomposition $P(\dd Z^n\,|\,X^n) = P|_{A_{n,j}}(\dd Z^n\,|\,X^n)\, P(A_{n,j}\,|\,X^n) + P|_{A_{n,j}^c}(\dd Z^n\,|\,X^n)\, P(A_{n,j}^c\,|\,X^n)$ for $P(\dd Z\,|\,X^n)$, the preceding display can be further decomposed into
\begin{align}
&(1 - \lambda)\int_{A_{n,j}^c}  \log(\frac{\wht{Q}|_{A_{n,j}^c}(\dd Z^n)}{P|_{A_{n,j}^c}(\dd Z^n\,|\,X^n)})\,\wht{Q}|_{A_{n,j}^c}(\dd Z^n) + \lambda \int_{A_{n,j}} \log(\frac{\wht{Q}|_{A_{n,j}}(\dd Z^n)}{P_{A_{n,j}}(\dd Z^n\,|\,X^n)})\,\wht{Q}|_{A_{n,j}}(\dd Z^n)\notag\\
&\qquad\qquad\qquad\qquad\qquad\qquad\qquad
 + (1 - \lambda)\log(\frac{1-\lambda}{P(A_{n,j}^c\,|\,X^n)}) + \lambda\log(\frac{\lambda}{P(A_{n,j}\,|\,X^n)})  \notag \\
=&  (1 - \lambda)\, d_{n2} + \lambda  \, d_{n1}
+ D\big(\text{Ber}(\lambda)\, \| \, \text{Ber}(\beta_n)\big),\label{eqn:lambda_obj}
\end{align}
where $\beta_n = P(A_{n,j}\,|\,X^n) = \Pi_n(|\theta_j-\theta^\ast_j|\geq D\varepsilon)$, 
$\text{Ber}(\lambda)$ denotes a Bernoulli distribution with success probability $\lambda$, and two constants $(d_{n1},d_{n2})$ independent of $\lambda$ are
\begin{align*}
d_{n1} = D\big(\wht{Q}|_{A_{n,j}}\, \| P|_{A_{n,j}}(\cdot\,|\,X^n)\big) \quad \mx{and} \quad d_{n2} = D\big(\wht{Q}|_{A_{n,j}^c}\, \| P|_{A_{n,j}^c}(\cdot\,|\,X^n)\big).
\end{align*}
Since $\wht \lambda$ minimizes $D\big(Q_\lambda^\dagger\, \| \,P(\cdot\,|\,X^n)\big)$, by setting the derivative of~\eqref{eqn:lambda_obj} as a function of $\lambda$ to be zero, we obtain
\begin{align*}
\wht{\lambda} =& \,\frac{\beta_n\, \exp(-d_{n1})}{\beta_n\, \exp(-d_{n1}) + (1-\beta_n)\, \exp(-d_{n2})} 
\leq \frac{\beta_n\,\exp(-d_{n1})}{(1-\beta_n)\, \exp(-d_{n2})}
\leq \frac{\beta_n}{1-\beta_n} \exp(d_{n2}),
\end{align*}
where we have used the fact that $d_{n1}\geq 0$ in the last step.
On the other hand, from the fact that~\eqref{eqn:lambda_obj} is equal to $D\big(Q_\lambda^\dagger\, \| \,P(\cdot\,|\,X^n)\big)$ at $\lambda=\wht\lambda$ and the non-negativeness of the three terms therein, we obtain $D\big(\wht Q\,||\, P(\cdot\,|\,X^n)\big)\geq (1-\wht\lambda) \, d_{n2} + D\big(\text{Ber}(\wht \lambda)\, \| \, \text{Ber}(\beta_n)\big)$. 
Combining this with the preceding display, we can reach
\begin{align}\label{Eqn:two_bounds}
\wht{\lambda} \leq 2\beta_n \exp(\frac{D\big(\wht{Q}\, \| \,P(\cdot\,|\,X^n))}{1 - \wht{\lambda}}) \quad \mx{and}\quad  D\big(\text{Ber}(\wht \lambda)\, \| \, \text{Ber}(\beta_n)\big) \leq D\big(\wht{Q}\, \| \,P(\cdot\,|\,X^n)),
\end{align}
where we have used inequality~\eqref{Eqn:post_contraction} so that $\beta_n\leq \Pi_n(\|\theta-\theta^\ast\|\geq D\varepsilon)\leq e^{-CD^2n\varepsilon^2}\leq n^{-CD^2} \leq 1/2$ for a sufficiently large $D$. 

Now we invoke the following lemma for bounding $D\big(\wht{Q}\, \| \,\Pi_n)$ from above. Its proof is deferred to Appendix~\ref{app:lem:KL_bound_proof} in the supplement.
\begin{lemma}\label{lem:KL_bound}
Under Assumptions A1, A2 and A3, it holds with probability at least $1-CM^{-2}$ that 
\begin{align*}
D\big(\wht{Q}\, \| \,P(\cdot\,|\,X^n)) \leq CnM^2\varepsilon_n^2.
\end{align*}
\end{lemma}
Using this lemma, the second display in inequality~\eqref{Eqn:two_bounds}, the bound in inequality~\eqref{Eqn:post_contraction} on $\beta_n$, and the fact that $\lambda \log \lambda +(1-\lambda ) \log(1-\lambda) \geq -\log 2$ for any $\lambda\in[0,1]$, we can obtain that $\wht\lambda \leq C/D^2\leq 1/2$ by choosing a sufficiently large $D$. Then a combination of the same Lemma~\ref{lem:KL_bound}, the inequality~\eqref{Eqn:post_contraction} on $\beta_n$, and the first display in inequality~\eqref{Eqn:two_bounds} on $\wht\lambda$ implies
\begin{align*}
\wht Q_{\theta}(|\theta_j-\theta^\ast_j|\geq D\varepsilon) = \wht{\lambda} \leq 2\exp(-CD^2 n\,\varepsilon^2 + 2CnM^2\varepsilon_n^2)\leq \exp(-CD^2n\varepsilon^2/2),
\end{align*}
for a sufficiently large constant $D$, which yields the claimed result on the tail probability of $\wht Q_{\theta}$ via a union bound over $j=1,\ldots,d$ (since $d^{-1/2}\|\theta-\theta^\ast\|\leq \max_{j} |\theta_j-\theta_j^\ast|$).


\subsection{Proof of Theorem~\ref{thm:1}}\label{Sec:proof_thm1}
We focus on illustrating the key proof idea of Theorem~\ref{thm:1} in this subsection, and leave the proofs of technical lemmas to the appendix. In addition, the proof of Theorem~\ref{thm:3} (weighted version of Theorem~\ref{thm:1}) will follow the same strategy, and details are also deferred to Appendix~\ref{app:Proof_thm:3}.

One main difficulty of the proof lies in the fact that the KL divergence does not satisfy the triangle inequality (unless the approximate family is convex, which is not true in the mean-field case), when specialized to our problem, taking the form as
\begin{align*}
D( \wht Q_\theta \,||\, Q^\ast_{VB} ) \leq D\big( \wht Q_\theta \,||\,N(\mle,\, [n I_c(\theta^\ast)]^{-1}) \big)-D\big( Q^\ast_{VB} \,||\,N(\mle,\, [n I_c(\theta^\ast)]^{-1}) \big),
\end{align*}
where recall that $Q^\ast_{VB}$ shares the same center $\mle$ as $N(\mle,\, [n I_c(\theta^\ast)]^{-1})$, but the precision matrix (inverse of covariance matrix) of the former is the diagonal part of the latter. This triangle inequality, if true, would imply the desired bound on $D( \wht Q_\theta \,||\, Q^\ast_{VB})$. In fact, according to Lemma~\ref{lem:KL_identity} below, $Q^\ast_{VB}$ minimizes $D\big(Q_\theta\,||\, N(\mle,\, [n I_c(\theta^\ast)]^{-1} \big)$ over all $Q_\theta$ within the mean-field family ($Q_{\theta}$ factorizes into $\bigotimes_{j=1}^d Q_{\theta_j}$). Moreover, the variational optimum $\wht Q_\theta$ also approximately minimizes this divergence since according to Lemma~\ref{lem:LQA_VB}, the KL divergence $D\big(Q_\theta \,||\,N(\mle,\, [n I_c(\theta^\ast)]^{-1}) \big)$ is, for all $Q_\theta$ with a similar sub-Gaussian tail as in Lemma~\ref{lem:expDacayverify}, close to the ``profile divergence" defined in Lemma~\ref{lem:decomposeKL}, for which $\wht Q_\theta$ is the minimizer since by definition $(\wht Q_{\theta},\wht Q_{S^n})$ jointly minimizes the variational objective function~\eqref{eqn:mean-field}. In other words, the preceding triangle inequality reveals the local strongly convexity structure of the profile variational objective function (after profiling out $Q_{S^n}$) around the populational level minimizer $Q^\ast_{VB}$. Unfortunately, such a triangle inequality for KL divergence is not true in general. Therefore, our first step of the proof is to establish a similar ``triangle inequality" restricted on the mean-field family with respect to the KL-divergence around the populational level minimizer $Q^\ast_{VB}$. After that, we will formalize the above intuition that $\wht Q_\theta$ effectively minimizes $D\big(Q_\theta \,||\,N(\mle,\, [n I_c(\theta^\ast)]^{-1}) \big)$. 

\vspace{0.5em}
\noindent{\bf Step one:} In this step, we build a ``triangle inequality" for the KL-divergence restricting to the mean-field family.
We will repeatedly use the following decomposition of the KL-divergence $D(Q\,||\,P)$ when the first probability measure $Q$ belongs to the mean-field family and the second probability measure $P$ is a normal distribution. Here we use the notation $\mu_{Q_\theta}$ to denote the expectation of any probability measure $Q$ over a space $\Theta$, and for any vector $u\in\Theta$, use $Q_{u}$ to denote the translation of $Q$ whose expectation is $u$, that is, $Q(A + \mu_{Q_\theta}) = Q_{(u)}(A + u)$ for any measurable set $A\subset \Theta$. Let $\lambda_{\min}(\Gamma)$ denote the smallest eigenvalue of a positive definite matrix $\Gamma$.

\begin{lemma}\label{lem:KL_identity}
Let $Q$ be a probability measure on $\mb R^d$ that factorizes as $Q=\bigotimes_{j=1}^d Q_j$, $\mu \in\mb R^d$ be a $d$-dim vector and $\Gamma\in\mb R^{d\times d}$ a $d$-by-$d$ positive definite matrix.
Let $Q^\ast = N(\mu, (\text{diag}(\Gamma))^{-1})) = \bigotimes_{j=1}^d Q^\ast_j$, where $Q^\ast_j$ is the univariate normal distribution $N(\mu_j,\Gamma_{jj}^{-1})$ for $j=1,\ldots,d$.
Then we have the decomposition
\begin{align*}
D\big(Q\,||\,N(\mu,\,\Gamma^{-1})\big) - D\big(Q^\ast\,||\,N(\mu,\,\Gamma^{-1})\big)= \sum_{j=1}^d
D\big(Q_{(\mu_j),j}\,||\,Q^\ast_j\big) + \frac{1}{2}(\mu_{Q_\theta} - \mu)^T \Gamma\,(\mu_{Q_\theta} - \mu), 
\end{align*}
where $\mu_{Q_\theta}\in\mb R^d$ denotes the expectation of $Q$. Furthermore, we have
\begin{align*}
D(Q\, \|\, Q^\ast) \leq \frac{\max_j \Gamma_{jj}}{\lambda_{\min}(\Gamma)}\,
\bigl[ D(Q\, \|\, N(\mu, \Gamma^{-1})) - D(Q^\ast\, \| \, N(\mu, \Gamma^{-1})) \bigr].
\end{align*}
\end{lemma}
\noindent Our proof based on direct calculation is provided in Appendix~\ref{app:Proof_KL_identity}. As a direct consequence, the first identity in the lemma implies that $Q^\ast$ minimizes $D\big(Q\,||\,N(\mu,\,\Gamma^{-1})\big)$ over all $Q$ within the mean-field family. When the covariance matrix $\Gamma$ is a multiple of the identity matrix, the second inequality in the lemma has leading factor $1$ and becomes the triangle inequality for the KL divergence at $Q^\ast$. We will apply this result with $Q^\ast$ as $Q^\ast_{VB}$ and $N(\mu,\Gamma^{-1}) = N\big(\mle,[n I_c(\theta^\ast)]^{-1}\big)$ below in the proof.

\vspace{0.5em}
\noindent{\bf Step two:} In this step, we study the profile divergence defined in Lemma~\ref{lem:decomposeKL}. 
Specifically, let $F_n(Q_\theta)=\min_{Q_{S^n}=\bigotimes_{i=1}^n Q_{S_i}}D\big(Q_\theta \otimes Q_{S^n} \,\|\, P(\cdot \,|\, X^n)\big)$ be the profile divergence.
The following lemma provides an approximation formula for $F_n(Q_\theta)$ for all $Q_\theta$ with a suitable tail decay property. Here, for any probability measure $Q$, we use $\Sigma_Q$ to denote the $d$-by-$d$ covariance matrix of $Q$. In particular, $\Sigma_Q$ becomes diagonal when $Q$ belongs to the mean-field family.

\begin{lemma}\label{lem:LQA_VB}
Suppose Assumption A3 holds. Then for any $M\geq 1$, it holds with probability at least $1-CM^{-2}$ that for any probability measure $Q=\bigotimes_{j=1}^d Q_j$ satisfying the same sub-Gaussian tail decay as in Lemma~\ref{lem:expDacayverify}, we have
\begin{align*}
\Big|F_n(Q_\theta) - D\big(Q_\theta\,||\,\Pi_n\big) - \frac{n}{2}\,\text{tr}\big(\Sigma_{Q_\theta}\, I_s(\theta^\ast)\big)\Big| \leq \frac{CM^3(\log n)^{3/2}}{\sqrt{n}},
\end{align*}
where $I_s(\theta^\ast)$ is defined in Assumption A3.
\end{lemma}
\noindent The proof is deferred to Appendix~\ref{app:Proof_LQA_VB} based on Taylor expansions. When there is no latent variable, $F_n(Q_\theta)$ is exactly $D\big(Q_\theta\,||\,\Pi_n\big)$. This means that the third extra term $\frac{n}{2}\,\text{tr}\big(\Sigma_{Q_\theta}\, I_s(\theta^\ast)\big)$ is due to the mean-field approximation between the parameter $\theta$ and latent variables $S^n$. This third term is not negligible since it will contribute to the precision matrix $nI_{VB}$ of the normal approximation $Q^\ast_{VB}$ to $\wht Q_\theta$, as we will show in the following steps.

\vspace{0.5em}
\noindent{\bf Step three:} The marginal posterior distribution $\Pi_n$ of $\theta$ in the approximation of $F_n(Q_\theta)$ in Lemma~\ref{lem:LQA_VB} is not convenient for our analysis. However, from the classical Bernstein von-Mises theorem, $\Pi_n$ should be well approximated by a normal distribution $N\big(\mle, \, [n I(\theta^\ast)]^{-1}\big)$. The following lemma shows that for any $Q_\theta$ with a suitable tail decay property, $D\big(Q_\theta\,||\,N\big(\mle, \, [n I(\theta^\ast)]^{-1}\big)\big)$ provides a good approximation to $D\big(Q_\theta\,||\,\Pi_n\big)$. Note that this lemma only assumes a sub-Gaussian type tail but not the mean-field structure (therefore it implies a KL divergence version of the BvM theorem). Then we may apply the first identity in Lemma~\ref{lem:KL_identity} to analyze this KL divergence term.
\begin{lemma}\label{lem:KL_normal_approx}
Suppose Assumptions A1 and A2 hold. Then for any $M\geq 1$, it holds with probability at least $1-CM^{-2}$ that for any probability measure $Q_\theta$ satisfying the same sub-Gaussian tail decay as in Lemma~\ref{lem:expDacayverify}, we have
\begin{align*}
\Big|D\big(Q_\theta\,||\,\Pi_n\big) -D\big(Q_\theta\,||\,N\big(\mle, \, [n I(\theta^\ast)]^{-1}\big)\big)\Big| \leq \frac{CM^3(\log n)^{d+3}}{\sqrt{n}},
\end{align*}
where $I(\theta^\ast)$ is the information matrix defined in Assumption A2.
\end{lemma}
\noindent The proof of this lemma is provided in Appendix~\ref{app:proof_KL_normal_approx}. In particular, we will apply this lemma with $Q_\theta$ as $\wht Q_\theta$ and $Q^\ast_{VB}$ that both satisfy the sub-Gaussian tail decay property in Lemma~\ref{lem:expDacayverify}.

\vspace{0.5em}
\noindent{\bf Step four:} Let $\alpha_n= \frac{CM^3(\log n)^{d+3}}{\sqrt{n}}$ denote the error upper bound in Lemmas~\ref{lem:LQA_VB} and~\ref{lem:KL_normal_approx}. In this step, we will formalize the intuition that $\wht Q_\theta$ effectively minimizes $D\big(Q_\theta \,||\,N(\mle,\, [n I_c(\theta^\ast)]^{-1}) \big)$. For any $Q=\bigotimes_{j=1}^d Q_j$ satisfying the same sub-Gaussian tail decay as in Lemma~\ref{lem:expDacayverify}, we have by Lemmas~\ref{lem:LQA_VB} and~\ref{lem:KL_normal_approx} that
\begin{align}\label{Eqn:equivalent_OBJ}
\Big|\,F_n(Q_\theta) - D\big(Q_\theta\,||\,N\big(\mle, \, [n I(\theta^\ast)]^{-1}\big)\big) -\frac{n}{2}\,\text{tr}\big(\Sigma_{Q_\theta}\, I_s(\theta^\ast)\big)\Big|& \leq 2\alpha_n.
\end{align}
This inequality indicates that $\wht Q_\theta$ is effectively minimizing the (negative) sum of the second and third terms in it, since $\wht Q_\theta$ minimizes $F_n(Q_\theta)$.
Next we will relate this approximate objective function with $D\big(Q_\theta \,||\,N(\mle,\, [n I_c(\theta^\ast)]^{-1}) \big)$.

Using equation~\eqref{Eqn:KL_expression} twice in the proof of Lemma~\ref{lem:KL_identity} in Appendix~\ref{app:Proof_KL_identity} with $(\mu,\Gamma)=(\mle, n I_c(\theta^\ast))$ and with $(\mu,\Gamma)=(\mle,n I(\theta^\ast))$, respectively, we obtain
\begin{align*}
&D\big(Q_\theta\,||\,N\big(\mle, \, [n I_c(\theta^\ast)]^{-1}\big)\big)
= \sum_{j=1}^d\int_{\mb R} q_j(\theta_j)\log q_j(\theta_j)\, \dd\theta_j 
+\frac{1}{2} \log\big((2\pi)^{d} |n I_c(\theta^\ast)|^{-1}\big) \\
&\qquad\qquad\qquad\qquad\qquad\qquad\qquad\qquad
+ \frac{n}{2}\,\text{tr}\big(\Sigma_{Q_\theta}\, I_c(\theta^\ast)\big) + \frac{n}{2}(\mu_{Q_\theta}-\mle)^T I_c(\theta^\ast)\,(\mu_{Q_\theta}-\mle),\\
&D\big(Q_\theta\,||\,N\big(\mle, \, [n I(\theta^\ast)]^{-1}\big)\big)
= \sum_{j=1}^d\int_{\mb R} q_j(\theta_j)\log q_j(\theta_j)\, \dd\theta_j 
+\frac{1}{2} \log\big((2\pi)^{d} |n I(\theta^\ast)|^{-1}\big)\\
 &\qquad\qquad\qquad\qquad\qquad\qquad\qquad\qquad
+ \frac{n}{2}\,\text{tr}\big(\Sigma_{Q_\theta}\, I(\theta^\ast)\big) + \frac{n}{2}(\mu_{Q_\theta}-\mle)^T I(\theta^\ast)\,(\mu_{Q_\theta}-\mle),
\end{align*}
where we have used the fact that $\Sigma_{Q_\theta}$ is a diagonal matrix, so that we can express the sum as the trace. Taking the difference between two and using the fact that $I_c(\theta^\ast)=I(\theta^\ast)+I_s(\theta^\ast)$, we can further obtain by rearranging the terms that
\begin{equation}\label{Eqn:KL-relation}
\begin{aligned}
&D\big(Q_\theta\,||\,N\big(\mle, \, [n I(\theta^\ast)]^{-1}\big)\big) + \frac{n}{2}\,\text{tr}\big(\Sigma_{Q_\theta}\, I_s(\theta^\ast)\big) \\
=&\, D\big(Q_\theta\,||\,N\big(\mle, \, [n I_c(\theta^\ast)]^{-1}\big)\big) - \frac{n}{2}(\mu_{Q_\theta}-\mle)^T I_s(\theta^\ast)\,(\mu_{Q_\theta}-\mle) + R_c,
\end{aligned}
\end{equation}
where $R_c=\frac{1}{2} \log\big( |I_c(\theta^\ast)|\cdot|I(\theta^\ast)|^{-1}\big)$ is a constant independent of $Q_\theta$. A combination of this identity with inequality~\eqref{Eqn:equivalent_OBJ} indicates that up to a translation, $\wht Q_\theta$ minimizes $D\big(Q_\theta \,||\,N(\mle,\, [n I_c(\theta^\ast)]^{-1}) \big)$. Moreover, the first identity in Lemma~\ref{lem:KL_identity} indicates that the KL-divergence also contains a translation related term (second term) that strictly dominates $-\frac{n}{2}(\mu_{Q_\theta}-\mle)^T I_s(\theta^\ast)\,(\mu_{Q_\theta}-\mle)$, so that the center can still be captured by minimizing $F_n(Q_\theta)$, as we will show in the next step below.

\vspace{0.5em}
\noindent{\bf Step five:} This is the last step where we will use the optimality of $\wht Q$ to prove the claimed bound. More specifically, by the optimaility of $\wht Q_\theta$ and the feasibility of $Q^\ast_{VB}$ for the optimization problem $\min_{Q_\theta=\bigotimes Q_{\theta_j}} F_n(Q_\theta)$, we obtain
\begin{align}\label{Eqn:optimality_condition}
F_n(\wht Q) \leq F_n(Q^\ast_{VB}).
\end{align}
Combining this inequality with~\eqref{Eqn:equivalent_OBJ} and~\eqref{Eqn:KL-relation} and using the fact that both $\wht Q_\theta$ (by Lemma~\ref{lem:expDacayverify}) and $Q^\ast_{VB}$ satisfy the sub-Gaussian tail condition therein, we can reach
\begin{equation}\label{Eqn:basic_ineq}
\begin{aligned}
 &D\big(\wht Q_\theta\,||\,N\big(\mle, \, [n I_c(\theta^\ast)]^{-1}\big)\big) - \frac{n}{2}(\mu_{\wht Q_\theta}-\mle)^T I_s(\theta^\ast)\,(\mu_{\wht Q_\theta}-\mle) \\
 &\qquad\qquad\qquad \qquad\qquad 
 \leq  D\big(Q^\ast_{VB}\,||\,N\big(\mle, \, [n I_c(\theta^\ast)]^{-1}\big)\big) + 4\alpha_n,
\end{aligned}
\end{equation}
where we have used the fact that the mean of $Q^\ast_{VB}$ is $\mle$.
Now we combine the above with the first identity in Lemma~\ref{lem:KL_identity} with $(\mu,\Gamma)=(\mle, n I_c(\theta^\ast))$ so that $Q^\ast=Q^\ast_{VB}$ to obtain,
\begin{align}\label{Eqn:mean_error}
\frac{n}{2}(\mu_{\wht Q_\theta}-\mle)^T I(\theta^\ast)\,(\mu_{\wht Q_\theta}-\mle) \leq 4\alpha_n,
\end{align}
where we have used the nonnegativeness of KL divergence. Since $I(\theta^\ast)$ is positive definite by Assumption A2, we have from the above that for some constant $C>0$,
\begin{align*}
\frac{n}{2}(\mu_{\wht Q_\theta}-\mle)^T I_s(\theta^\ast)\,(\mu_{\wht Q_\theta}-\mle) \leq C\alpha_n,
\end{align*}
Finally, by combining the preceding display, inequality~\eqref{Eqn:basic_ineq} and the second inequality in 
Lemma~\ref{lem:KL_identity}, we obtain
\begin{align*}
D\big(\wht Q_\theta\,||\,Q^\ast_{VB}\big)& \leq C\,\Big[D\big(\wht Q_\theta\,||\,N\big(\mle, \, [n I_c(\theta^\ast)]^{-1}\big)\big)-D\big(Q^\ast_{VB}\,||\,N\big(\mle, \, [n I_c(\theta^\ast)]^{-1}\big)\big)\Big]\\
& \leq C\alpha_n + \frac{n}{2}(\mu_{\wht Q_\theta}-\mle)^T I_s(\theta^\ast)\,(\mu_{\wht Q_\theta}-\mle) \leq C'\alpha_n,
\end{align*}
which concludes the proof.


\subsection{Proof of Corollary~\ref{cor:VB_center}}
The claimed bound is a direct consequence of inequality~\eqref{Eqn:mean_error} in the proof of Theorem~\ref{thm:1}.


\bibliographystyle{plain}
\bibliography{MF}

\newpage 
\appendix

  \begin{center}
    {\LARGE\bf Supplement to Statistical Inference in Mean-field Variational Bayes}
\end{center}

\section{Computational details in the numerical study}

\subsection{Gaussian mixture model}\label{app:GMM}
\noindent {\bf Computation details:}
Using the notation in Section~\ref{Sec:GMM}, the evidence lower bound (ELBO, equation~$(21)$ of paper~\cite{Review}) has the following explicit form,
\begin{equation}\label{eq:GmObj}
\begin{aligned}
L(q_{Z^n}) =&\, E_{q_{Z^n}}\big[\log\big(p(X^n, Z^n)) - \log\big(q_{Z^n}(Z^n)\big)\big]  \\
\propto&\, E_{q_{Z_n}}\Big[-\frac{\mu^T\mu}{2\sigma^2} - \frac 12\sum_i (x_i - c_i^T \mu)^2 -\sum_i \log(c_i^T \phi_i) + \sum_k \log(s_k^2) \Big] \\
=&\, -\frac{1}{2\sigma^2}\sum_k(s_k^2 + m_k^2) - \frac12\sum_{i,k} \phi_{ik}((x_i - m_k)^2 + s_k^2) - \sum_{i,k} \phi_{ik}\log(\phi_{ik}) + \frac12\sum_k \log(s_k^2),
\end{aligned}
\end{equation} 
where in the second line we have omitted a constant term independent of the variational posterior $q_{Z^n}$.
Since $Q_{Z^n}$ is parametrized by $\{(m_k,s_k^2\}_{k=1}^K$ and $\{(\phi_{i1},\ldots,\phi_{iK})\}_{i=1}^n$,  the $t$-th step inside the while loop in Algorithm~\ref{algo:vbB} (with unit weights) can be summarized as follows,
\begin{align*}
{\phi}_{ik}^{(t)} &\propto \exp\Big\{x_i m_k^{(t-1)} - \frac12 (s_k^{(t-1)})^2 + (m_k^{(t-1)})^2\Big\}, \quad \text{(normalized over $k=1,\ldots,K$)}\\
{m}_k^{(t)} &= \frac{\sum_i x_i\phi_{ik}^{(t)}}{1/\sigma^2 + \sum_i \phi_{ik}^{(t)}}, \ \ \mx{and} \ \ 
({s}_k^{(t)})^2 = \frac{1}{1/\sigma^2 + \sum_i \phi_{ik}^{(t)}}, \ \ k=1,\ldots, K, \ \ i=1,\ldots,n,
\end{align*}
where iteration proceeds until $L(q_{Z^n})$ stabilizes.\\ 
In our VWLB methods, the weighted CAVI algorithm can be updated in a similar fashion, where in the $t$-th step, we use the following updates:
\begin{align*}
{\phi}_{ik}^{(t)} &\propto \exp\Big\{x_i m_k^{(t-1)} - \frac12 (s_k^{(t-1)})^2 + (m_k^{(t-1)})^2\Big\}, \quad \text{(normalized over $k=1,\ldots,K$)}\\
{m}_k^{(t)} &= \frac{\sum_i  w_{bi}x_i\phi_{ik}^{(t)}}{1/\sigma^2 + \sum_i  w_{bi}\phi_{ik}^{(t)}}, \ \ \mx{and} \ \ 
({s}_k^{(t)})^2 = \frac{1}{1/\sigma^2 + \sum_i  w_{bi}\phi_{ik}^{(t)}}, \ \ k=1,\ldots, K, \ \ i=1,\ldots,n.
\end{align*}

\begin{figure}[ht!]
\centering
\subfloat[Gibbs Sampler]{\includegraphics[scale=0.23]{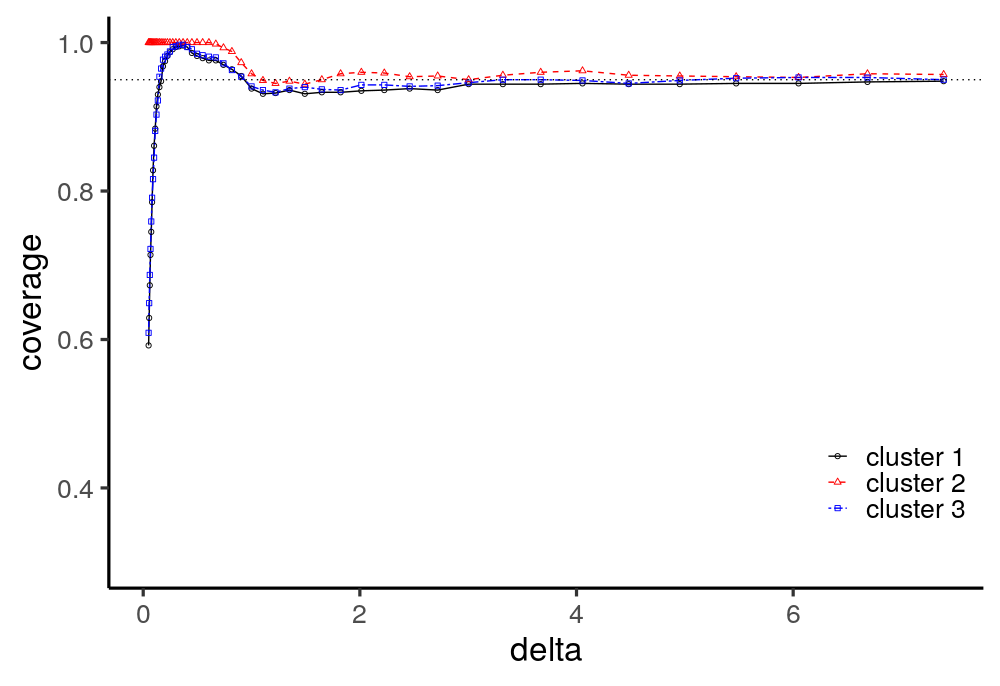}}
\qquad
\subfloat[VB posterior]{\includegraphics[scale=0.23]{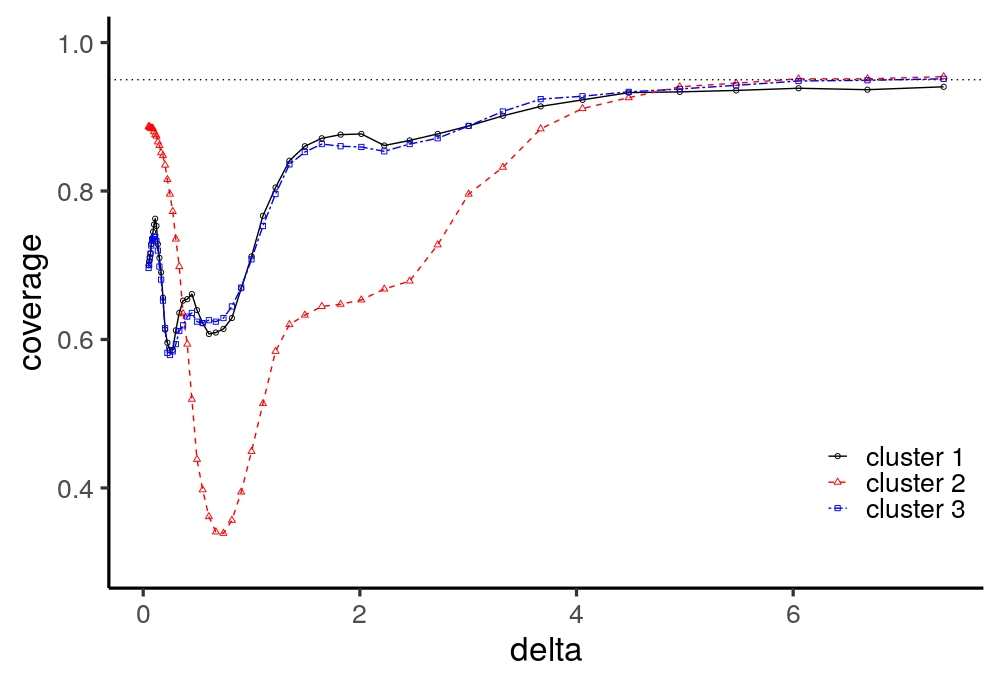}}
\\
\subfloat[VWLB CI]{\includegraphics[scale=0.23]{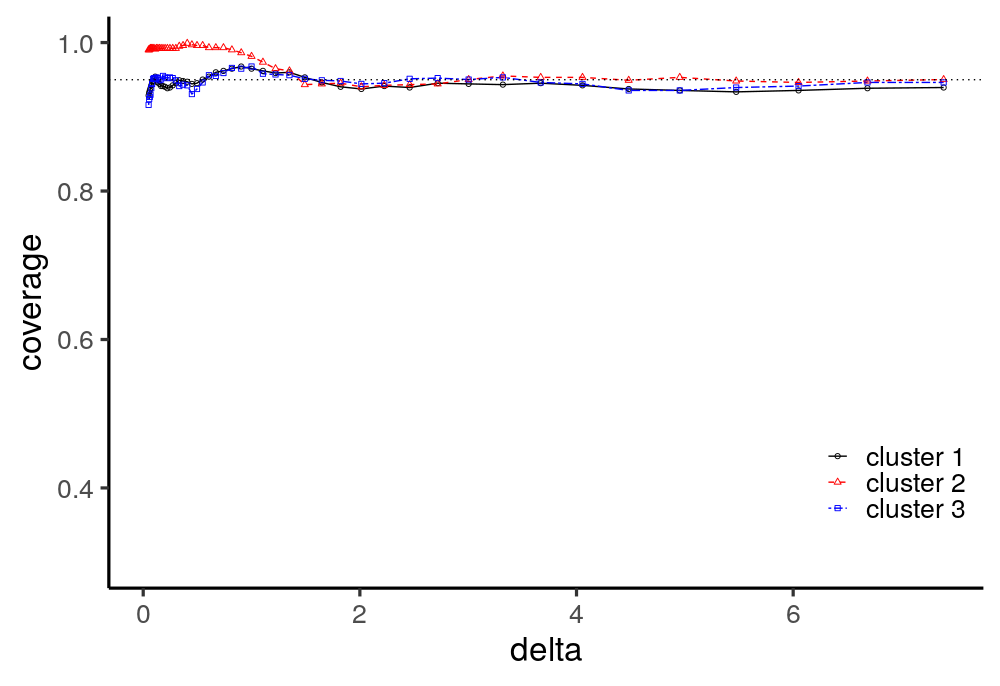}}
\qquad
\subfloat[VWLB CI2]{\includegraphics[scale=0.23]{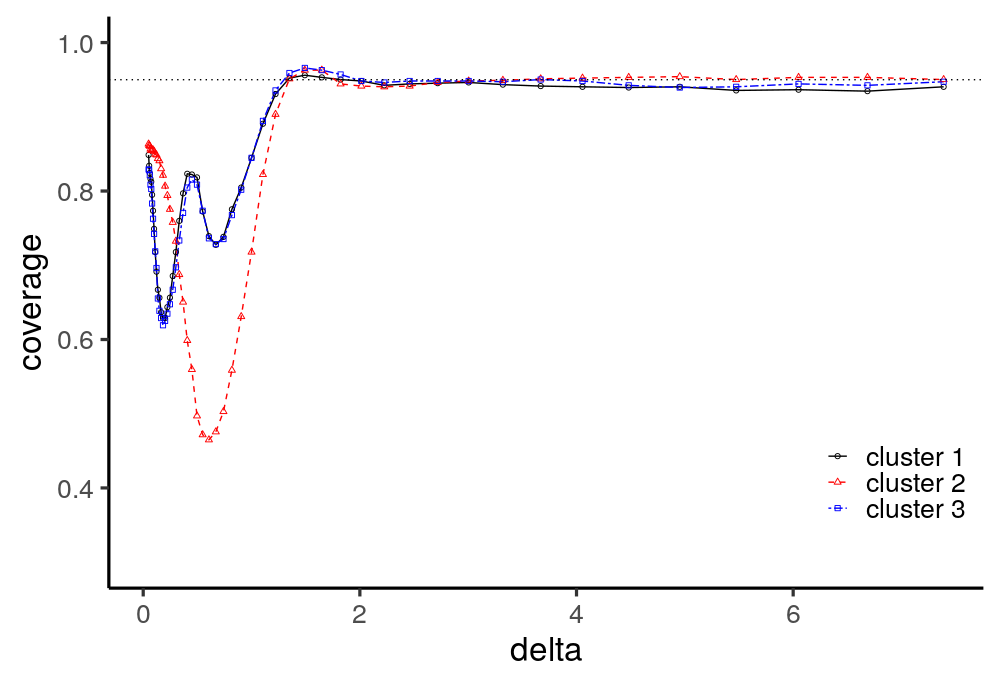}}
\caption{Coverage probabilities for each cluster center $\mu_k$, $k=1,2,3$, versus $\Delta$. }\label{fig:GMM_three}
\end{figure}

\begin{figure}[ht!]
\centering
\subfloat[$\Delta=0$]{\includegraphics[scale=0.26]{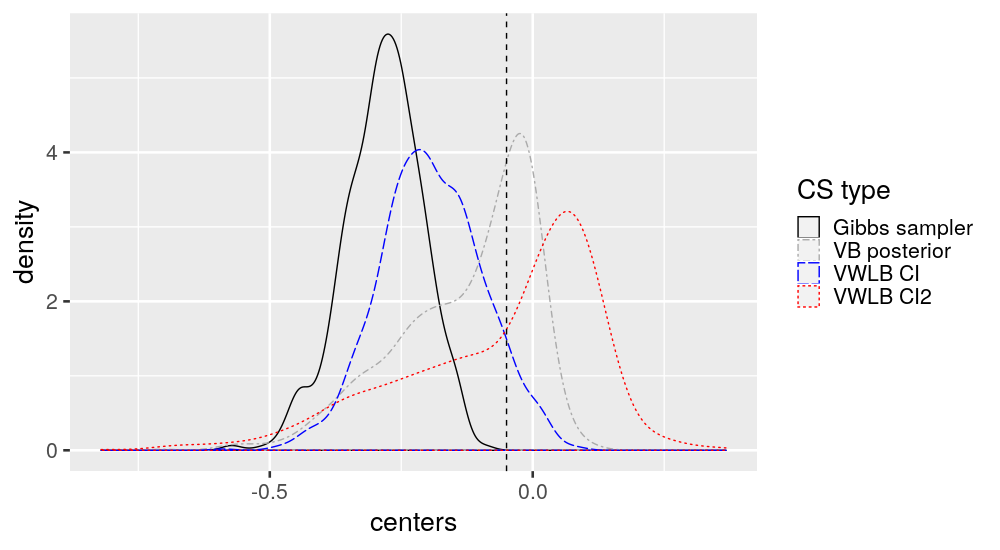}}
\qquad
\subfloat[$\Delta=1$]{\includegraphics[scale=0.26]{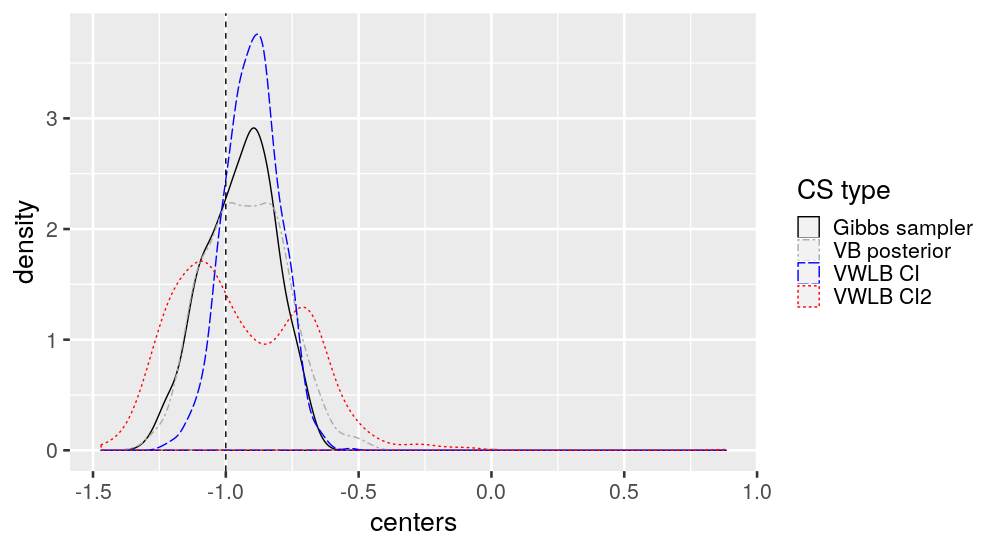}}
\\
\subfloat[$\Delta=3$]{\includegraphics[scale=0.26]{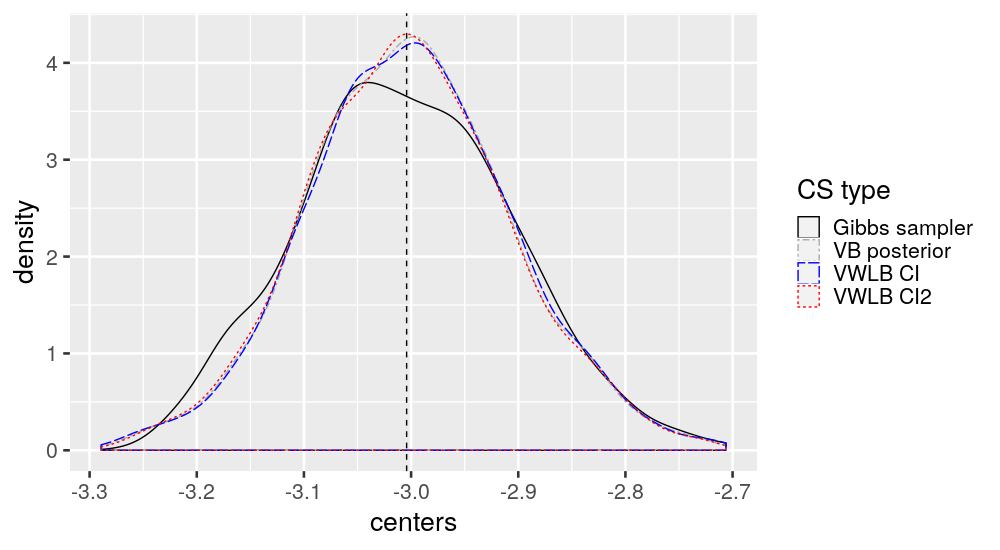}}
\qquad
\subfloat[$\Delta=5$]{\includegraphics[scale=0.26]{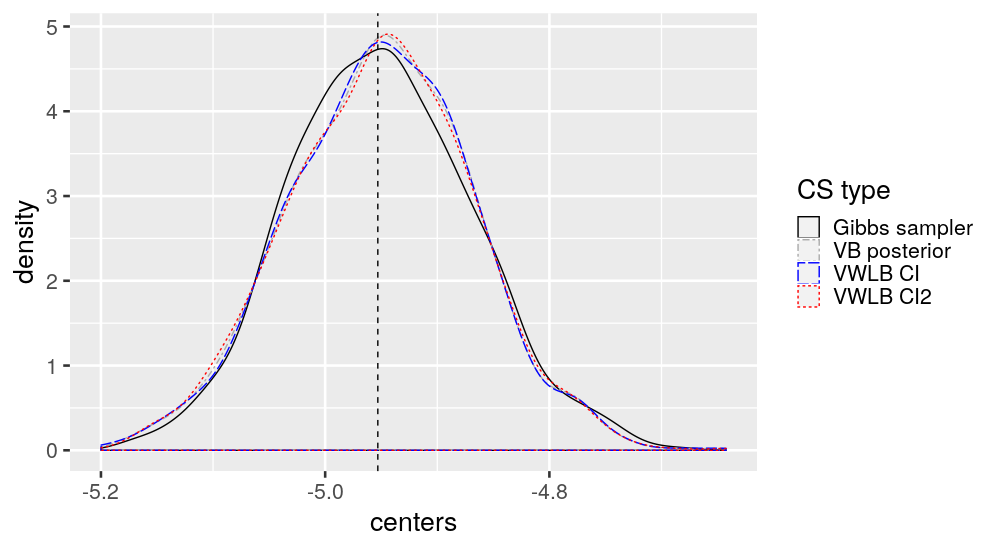}}
\caption{Empirical distribution of the centers of four credible intervals for $\mu_1$.}\label{fig:GMM_c_density}
\end{figure}

\begin{figure}[ht!]
\centering
\subfloat[$\Delta=0$]{\includegraphics[scale=0.26]{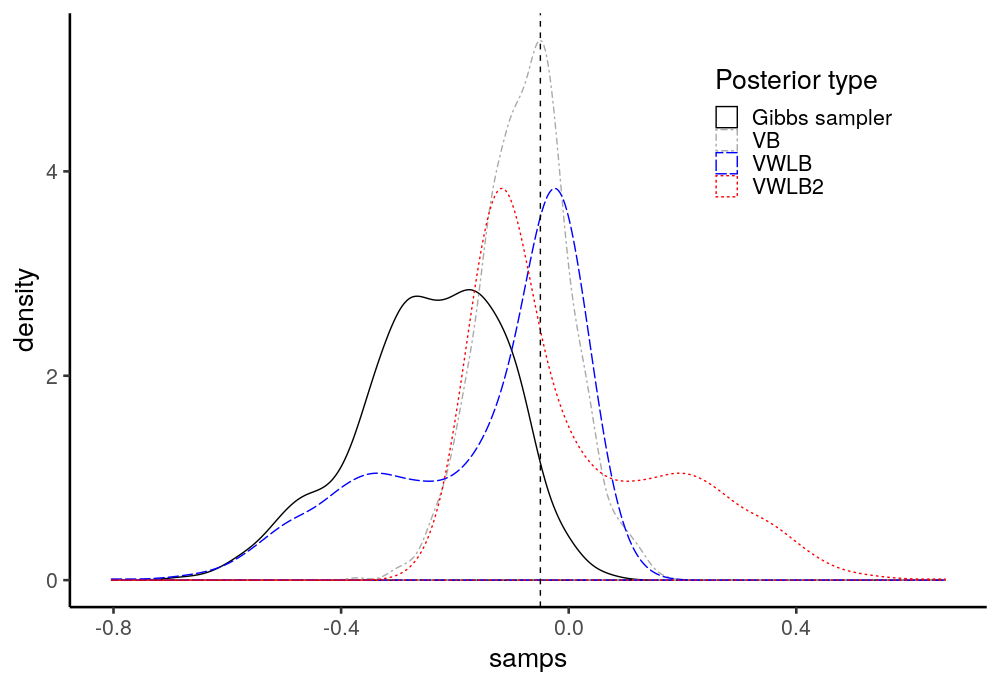}}
\qquad
\subfloat[$\Delta=1$]{\includegraphics[scale=0.26]{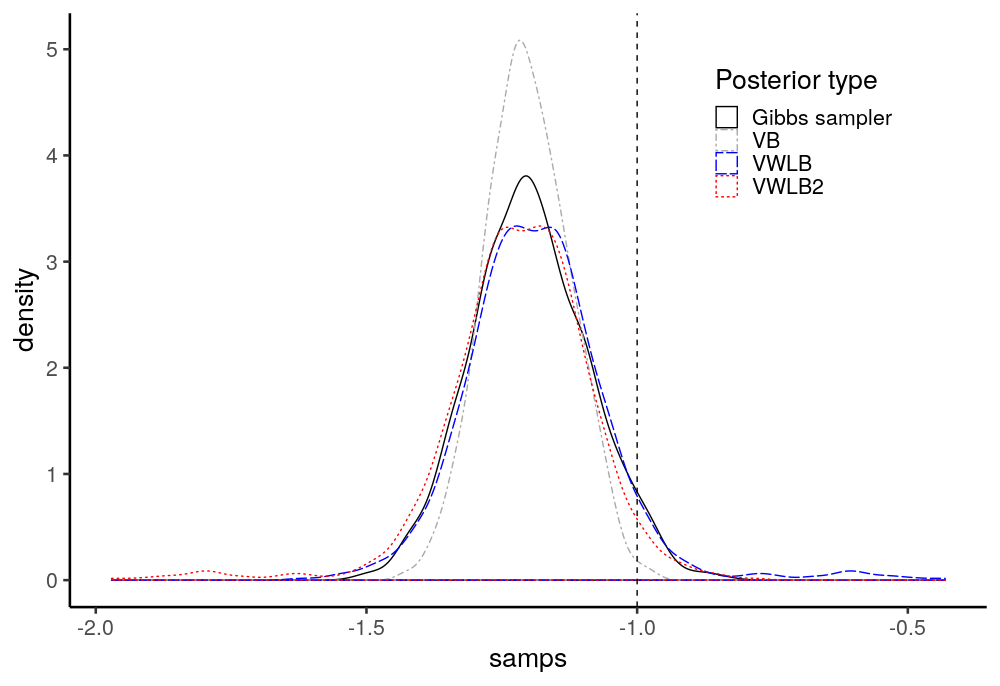}}
\\
\subfloat[$\Delta=3$]{\includegraphics[scale=0.26]{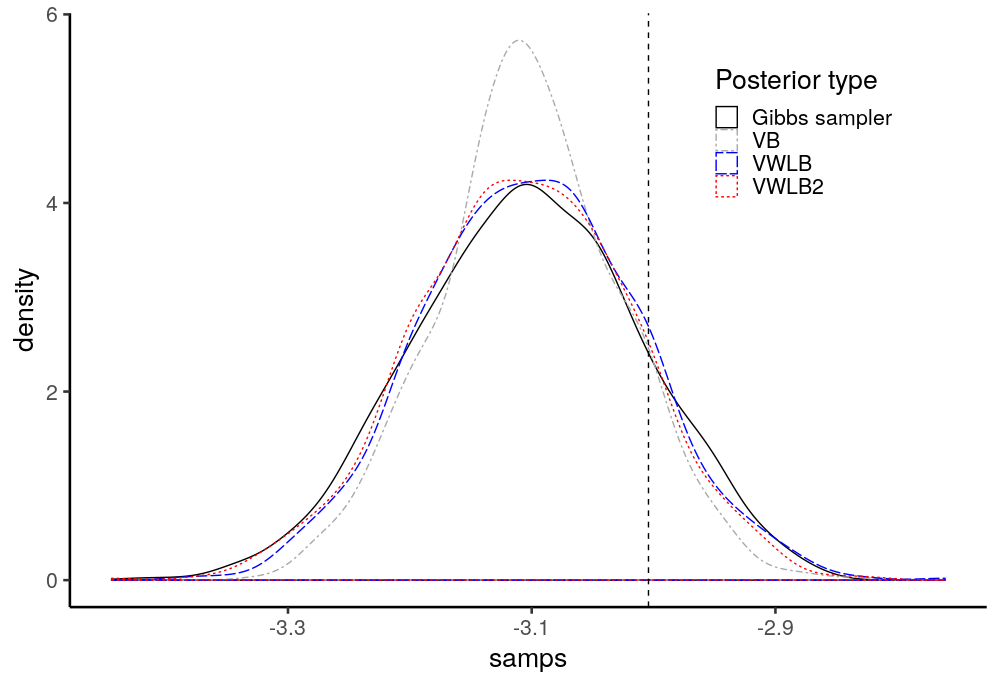}}
\qquad
\subfloat[$\Delta=5$]{\includegraphics[scale=0.26]{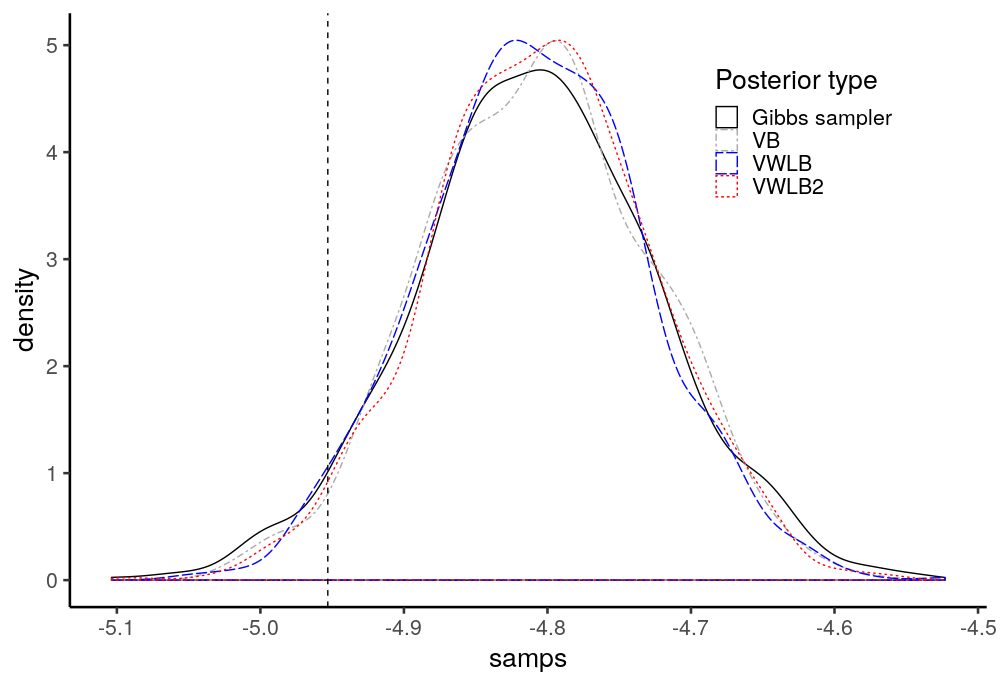}}
\caption{Approximated posteriors corresponding to four types of credible intervals.}\label{fig:GMM_density}
\end{figure}

\vspace{0.5em}
\noindent {\bf Detailed simulation results:}
We provide more detailed analysis on GMM when $n=500$, and explain different behaviors between the two credible interval construction schemes based on our VWLB methods. Firgure~\ref{fig:GMM_three} displays the individual coverage probabilities for each cluster center $\mu_k$, $k=1,2,3$, of four credible intervals versus the separation gap $\Delta$. From these results, it appears that the credible intervals based on samplers $\{\wt\theta^{(b)}_{VB}\}_{b=1}^B$ directly drawn from VWLB achieve the nominal level even when the model becomes degenerate ($\Delta\to 0_+$), while those based on Gibbs sampler (or true posterior) go from under-covering, to over-covering, and finally become stabilized at the nominal level. In contrast, then other two methods exhibit more drastic under-covering issues at small $\Delta$ values.
Figure~\ref{fig:GMM_c_density} shows the empirical distributions of the credible intervals centers for $\mu_1$ at $\Delta\in\{0,1,3,5\}$, which explains the discrepancy between the two methods based on the VWLB. From this plot, we can see that the second one (VWLB CI2) based on the reverting idea in bootstrap of using quantiles of $\{2\wt\theta_{VB}^{(b)}-\wht\theta_{VB}\}_{b=1}^B$ further suffers from the extra variability due to the variational posterior mean $\wht\theta_{VB}$ that causes the bimodal distribution for the interval centers at small $\Delta$ values. In contrast, the first one (VWLB CI) does not use $\wht\theta_{VB}$ for correcting the center, and therefore is able to avoid introducing the extra systematic bias in $\wht\theta_{VB}$ due to the model degeneracy. Finally, we provide the approximated posterior distribution corresponding to the four types of credible intervals in Figure~\ref{fig:GMM_density}. As we can expect, the mean-field approximation (VB) always underestimates the dispersion of the posterior distribution (which is well-approximated by the Gibbs sampler), especially at small $\Delta$ values.
In contrast, the posterior approximations from the two VWLB samplers become close to the true posterior as $\Delta$ exceeds value $1$.

\subsection{Bayesian linear regression}\label{app:BLM}
\noindent {\bf Computation details:} Recall that $X = (X_1, \ldots, X_p)$ is the $n\times p$ design matrix. For each $j=1,\ldots,p$, let $R_j = Y^n - \sum_{j'\neq j} X_{j'}\beta_{j'}$ denote the vector of residuals without the $j$-th covariate $X_j$. The Gibbs sampler based on the point mass mixture prior for $\beta$ consists of cycling through sampling from the following full conditionals ($\cdot\,|\,-$ stands for conditioning on the rest parameters): for $j=1,2,\ldots,p$,
\begin{equation}\label{eq:BLMGibbsUpdate1}
\begin{split}
[\gamma_j\,|\,-] &\,\sim \mathrm{Ber}(\theta_j), \ \text{with} \
\theta_j = \frac{\xi}{\xi + \frac{1 - \xi}{\sqrt{1 + \nu_1 \|X_j \|^2}} \exp{ \frac{\langle R_j, X_j \rangle^2 }{2 \sigma^2 (1/\nu_1 + \|X_j \|^2)} - \frac{\|\beta \|^2 - \beta_j^2}{2\nu_1 \sigma^2} } },\\
[\beta_j \,|\,-] & \, \sim   1(\gamma_j=0)\, \delta_0 +  1(\gamma_j=1)\, N\left( \frac{\langle R_j, X_j \rangle}{1/\nu_1 + \|X_j \|^2} , \ \frac{\sigma^2}{1/\nu_1 + \|X_j \|^2} \right),
\end{split} 
\end{equation}
\begin{equation}\label{eq:BLMGibbsUpdate2}
\begin{split}
[\sigma^2\,|\,-] &\, \sim \mathrm{IG}\Big(\nu/2 + n/2 + \sum_i (1 - \gamma_i)/2 , \,\nu\lambda/2 + \|Y^n - X\beta \|^2/2 + \|\beta \|^2/(2\nu_1)\Big) \\
[\xi\,|\,-] &\, \sim \mathrm{Beta}\Big( a_0 + \sum_i \gamma_i , b_0 + \sum_i (1 - \gamma_i) \Big)
\end{split}
\end{equation}

Since $\beta$ and $\gamma=\{\gamma_j\}_{j=1}^p$ are our primary parameters of interest, we adopt the same strategy as in~\cite{FastBLR} by estimating $\sigma^2$ and $\xi$ by their respective maximum a posterior (MAP) estimators, while still using the block mean-field approximation as in~\eqref{Eqn:MF_block} on $(\gamma,\beta)=\{(\gamma_j,\beta_j)\}_{j=1}^p$ to facilitate fast computation. According to the derivation in~\cite{FastBLR}, under this setup the evidence lower bound takes the form as
\begin{align*}
L(q_{Z^n})=& \, E_{q_{\beta,\gamma}}\Bigl[ \log\big(p(Y^n,\beta,\gamma,\sigma^2,\xi)\big) - \log(q_{\beta,\gamma}(\beta, \gamma))
\Bigr] \\
=&\, -\frac{n}{2} \log(\sigma^2) - \frac{1}{2\sigma^2} E_{q_{\beta,\gamma}}\Bigl[ (Y^n - X\beta)^T (Y^n - X\beta)\Bigr] + (a_0 - 1)\log(\xi) \\
& \,+(b_0 -1) \log(1 - \xi)-\Big(\frac{\nu}{2} + 1\Big)\log(\sigma^2) - \frac{\nu\lambda}{2\sigma^2} + \sum_{j = 1}^p E_{q_{\beta,\gamma}}\Bigl[ 
\log(\frac{\pi(\beta_j \,|\, \gamma_j)\, \pi(\gamma_j \,|\, \xi)}{q_{\beta_j,\gamma_j}(\beta_j,\gamma_j)}) 
\Bigr].
\end{align*}
The detailed steps of the coordinate ascent algorithm for optimizing $L(q_{Z^n})$ jointly over $q_{Z^n}$ and $(\sigma^2,\xi)$ can be found in paper~\cite{FastBLR}.\\
The weighted CAVI Algorithm~\ref{algo:vbB} in our VWLB sampler can proceed in a similar way, where the only difference is in replacing the sum of square $(Y^n - X\beta)^T (Y^n - X\beta)$ by its weighted version $(Y^n - X\beta)^T W^n (Y^n - X\beta)$ in the above ELBO $L(q_{Z^n})$, where $W^n$ is the $n\times n$ diagonal matrix whose $i$-th diagonal component is $W_i^{(b)}$.

\section{Proofs of results in the main paper}
In this section, we provide the remaining proofs of the results in the paper.


\subsection{Proof of Lemma~\ref{lem:decomposeKL}}\label{Sec:proof_lem:decomposeKL}
By explicitly writing out the integral in the KL-divergence, we obtain that for any $q_{Z^n}= q_\theta\otimes q_{S^n}$,
\begin{align}
D(q_{\theta}\otimes q_{S^n} \| p(\cdot\,|\,X^n)) 
&=D(q_{\theta}\, \|\, \Pi_n) + \int q_{\theta}(\theta)\, D\big(q_{S^n}\,\|\, p_{S^n}(\cdot\, | \,X^n, \theta)\big)\, \dd \theta \notag\\
&= D(q_{\theta}\, \|\, \Pi_n) + \sum_{i=1}^n \int q_{\theta}(\theta)\, D\big(q_{S_i}\,\|\, p_{S_i}(\cdot\, | \,X_i, \theta)\big)\, \dd \theta, \label{eq:fullExpansion}
\end{align}
where recall that $\Pi_n$ is the marginal posterior distribution of $\theta$ given $X^n$. 
By the definition of the $r_i(s)$ function in the lemma, we have
\begin{equation*}
q_{S_i}(s_i) = \frac{\exp(r_i(s_i))}{\int_{\mathcal{S}} \exp(r_i(s'_i))\,\dd s'_i}.
\end{equation*}
By substituting the above into the second term of equation~\eqref{eq:fullExpansion}, we obtain
\begin{align*}
& \int q_{\theta}(\theta)\, D\big(q_{S_i}\,\|\, p_{S_i}(\cdot\, | \,X_i, \theta)\big)\, \dd \theta \\
= &\,\int q_{\theta}(\theta)\, \int_{\m S} q_{S_i}(s_i)\, \Big(\log\big(q_{S_i}(s_i)\big) - \log\big(p(s_i\, |\, X_i, \theta)\big)\Big)\, \dd s_i\dd \theta \\
= &\,\int_{\m S} q_{S_i}(s_i) \,\Big(\log\big(q_{S_i}(s_i)\big) - \int q_{\theta}(\theta)\, \log\big(p(s_i\,| \, X_i, \theta)\big)\,\dd \theta\Big) \,\dd s_i \\
= &\,\int_{\m S} q_{S_i}(s_i) \,\Big[r_i(s_i) - \log \Big(\int_{\mathcal{S}} \exp\big(r_i(s'_i)\big)\,\dd s'_i\Big) - r_i(s_i) \Big] \, \dd s_i \\
= &\, - \log \Big(\int_{\mathcal{S}} \exp\big(r_i(s'_i)\big)\,\dd s'_i\Big).
\end{align*}
A combination of the above with equation~\eqref{eq:fullExpansion} leads to the claimed identity.


\subsection{Proof of Proposition~\ref{prop:wposterior_mode}}
The proof of this proposition is similar to the proof of Lemma~\ref{lem:deIntApprox} in Appendix~\ref{app:proof_deIntApprox} for the unweighted posterior. We just need to point out the difference. 

In fact, the approximation bound~\eqref{Eqn:normalization} of Lemma~\ref{lem:deIntApprox} can be rewritten as
\begin{equation}\label{Eqn:posterior_meana}
\begin{aligned}
&\int \pi(\theta)\,\exp\big\{l(\theta; X^n) - l(\mle; X^n)\big\}\,\dd\theta  \\
&\qquad
= \bigg(1 + \m O\Big( \frac{CM^3(\log n)^{d+3}}{\sqrt{n}}\Big)\bigg)\cdot \int \pi(\theta^\ast)\,  \exp\big\{-\frac{n}{2}\,(\theta-\mle)^T I(\theta^\ast)\, (\theta-\mle)\big\}\,\dd\theta.
\end{aligned}
\end{equation}
An almost same line-by-line derivation can be used to prove 
\begin{equation}\label{Eqn:posterior_meanb}
\begin{aligned}
&\int \|\theta-\mle\|\,\pi(\theta)\,\exp\big\{l(\theta; X^n) - l(\mle; X^n)\big\}\,\dd\theta  \\
&\quad= 
\bigg(1 + \m O\Big( \frac{CM^3(\log n)^{d+3}}{\sqrt{n}}\Big)\bigg)\cdot \int \|\theta-\mle\|\, \pi(\theta^\ast)\,  \exp\big\{-\frac{n}{2}\,(\theta-\mle)^T I(\theta^\ast)\, (\theta-\mle)\big\}\,\dd\theta,
\end{aligned}
\end{equation}
by instead analyzing the integration of the difference between
\begin{align*}
i'_1(\theta)& =\|\theta-\mle\|\,\pi(\theta)\,\exp\big\{l(\theta; X^n) - l(\mle; X^n)\big\},\quad\mx{and}\\
i'_2(\theta) &= \|\theta-\mle\|\, \pi(\theta^\ast)\,  \exp\big\{-\frac{n}{2}\,(\theta-\mle)^T I(\theta^\ast)\, (\theta-\mle)\big\},
\end{align*}
via the same strategy of dividing $\mb R^d$ into $A_1$ and $A_2$. The integral over $A_1$ can be handled by exactly the same way of a local Taylor expansion to the exponent in the exponential, and the integral over $A_2$ by using the sub-Gaussian tail bound of $\Pi_n$ guaranteed by Lemma~\ref{lem:expDacayverify}. Note that inequalities~\eqref{Eqn:posterior_meana} and~\eqref{Eqn:posterior_meanb} together imply
$\|\int \theta\, \pi_n(\theta)\,\dd\theta - \mle\| \leq \frac{CM^3(\log n)^{d+3}}{\sqrt{n}}$, a bound between the posterior mean and the MLE.

Return to the current weighted posterior scenario. Similarly, the claimed bound is implied by a weighted version of~\eqref{Eqn:posterior_meana} and~\eqref{Eqn:posterior_meanb} as
\begin{align*}
&\int \|\theta-\wt\theta^{(b)}\|^{k}\,\pi(\theta)\,\exp\big\{\wt l^{(b)}(\theta; X^n) - \wt l^{(b)}(\wt\theta^{(b)}; X^n)\big\}\,\dd\theta  \\
& = \bigg(1 + \m O\Big( \frac{CM^3(\log n)^{d+3}}{\sqrt{n}}\Big)\bigg)\cdot
\int\|\theta-\wt\theta^{(b)}\|^{k}\, \pi(\theta^\ast)\,  \exp\big\{-\frac{n}{2}\,(\theta-\wt\theta^{(b)})^T I(\theta^\ast)\, (\theta-\wt\theta^{(b)})\big\}\,\dd\theta,
\end{align*}
for $k=0,1$, 
where $\wt\theta^{(b)}$ maximizes the log-weighted likelihood $\wt l^{(b)}(\theta; X^n) = \log\big[ \wt L^{(b)}(\theta;\, X^n)\big]$, and plays the role of the MLE $\mle$ in the unweighted case. The integral over $A_1$ can be handled by the same way of a local Taylor expansion at $\wt\theta^{(n)}$, and the integral over $A_2$ by using the sub-Gaussian tail bound of $\wt\Pi_n^{(b)}$ guaranteed by inequality~\eqref{Eqn:contraction_WP} of Theorem~\ref{thm:WMF_contraction} for the weighted posterior distribution.

\subsection{Proof of Lemma~\ref{lem:decomposeWKL}}
The proof is almost the same as that of Lemma~\ref{lem:decomposeKL} by explicitly writing out the integral in the KL divergence decomposition~\eqref{Eqn:WKL_decom}, where the only difference is in keeping track of the weights.


\subsection{Proof of Theorem~\ref{thm:WMF_contraction}}\label{app:Proof_thm_WMF_contraction}
For the first inequality~\eqref{Eqn:contraction_WP}, we generalize the contraction result~\eqref{Eqn:post_contraction} to the weighted posterior via the adapting the proof strategy (for the usual posterior) in~\cite{RatesPosterior}.
In particular, a key ingredient of our proof is based on an extension of the probability inequality developed in~\cite{wong1995} from controlling the likelihood ratio empirical process to controlling the weighted likelihood ratio process, as in Lemma~\ref{lem:locent}. Following the notation of \cite{RatesPosterior}, we let $\wt m_n(A) = \int_A \widetilde L^{(b)}(\theta;\, X^n)/\widetilde L^{(b)}(\theta_0;\, X^n) \, \dd \Pi(\theta)$ for any measurable set $A\subset \Theta$. Then the weighted likelihood posterior~\eqref{eq:weightedPo} can be rewritten as:
\begin{equation*}
\widetilde{\Pi}^{(b)}_n(A) = \frac{\int_A \widetilde L^{(b)}(\theta;\, X^n)\, \dd \Pi(\theta)}{\int_\Theta \widetilde L^{(b)}(\theta;\, X^n)\, \dd \Pi(\theta)} 
= \frac{\wt m_n(A)}{\wt m_n(\Theta)},\quad\mx{for any measureable set $A\subset\Theta$,}
\end{equation*}
where we have used the fact that $\widetilde L^{(b)}(\theta_0;\, X^n)$ is free of $\theta$. 
Let $A_n := \{ \theta \in \Theta: \| \theta - \theta^\ast \| \leq C_2' \varepsilon \}$. Then the desired inequality~\eqref{Eqn:contraction_WP} is equivalent to $\frac{\wt m_n(A_n^c)}{\wt m_n(\Theta)}\leq e^{-C_3'n\varepsilon^2}$.
We analyze the denominator and numerator separately by using the following two lemmas, whose proofs are deferred to Appendix~\ref{app:proof_WPdenominator} and~\ref{app:proof_WPnumerator} respectively.

\begin{lemma}[Denominator lower bound]\label{lem:WPdenominator}
Under Assumptions A2, W, for every $\delta > 0$ and any $C>0$, it holds with probability at least $1 - \frac{2}{C^2n\delta^2}$ that
\begin{equation}\label{eq:WPdenominator}
\wt m_n(\Theta) \geq \Pi(\theta\in \Theta: \|\theta - \theta^\ast \|^2 \leq C_0\, \delta^2 ) e^{-(1+C)\,n\delta^2}, \quad\mx{for some $C'>0$}.
\end{equation}
\end{lemma}

\begin{lemma}\label{lem:WPnumerator}
Under Assumptions A2 and W, for any $M\geq 1$, it holds with probability at least $1 - CM^{-2}$ that
\begin{equation}\label{eq:WPnumerator}
\sup_{\| \theta - \theta^\ast \| \geq C_2' \varepsilon } \frac{\widetilde L^{(b)}(\theta;\, X^n)}{ \widetilde L^{(b)}(\theta^\ast;\, X^n)} \leq \exp(-n\varepsilon^2/24), \quad\mx{for any $\varepsilon\geq M\wt\varepsilon_n$.}
\end{equation}
\end{lemma}

\noindent The uniform bound~\eqref{eq:WPnumerator} on the weight likelihood ratio immediately implies $\wt m_n(A_n^c) \leq \exp(-n\varepsilon^2/24)$.
By combining this with the denominator bound~\eqref{eq:WPdenominator}(with $C = 1$, and $\delta = \varepsilon/16$), we obtain
\begin{equation*}
\widetilde{\Pi}^{(b)}_n(A_n^c)= \frac{\wt m_n(A_n^c)}{\wt m_n(\Theta)} \leq  \exp(-C_3''n\varepsilon^2) \quad\mx{holds for any $\varepsilon\geq M\wt \varepsilon_n$,}
\end{equation*}
with probability at least $1 - CM^{-2}$.

The second part concerning the sub-Gaussian tail of $\wt Q^{(b)}_\theta=\wt Q_\theta^{(b)}\otimes \wt Q_{S^n}^{(b)}$ follows a similar argument as the proof of Lemma~\ref{lem:expDacayverify} in Section~\ref{lem:KL_bound} by considering for each coordinate index $j=1,\ldots,d$ the same distribution family indexed by $\lambda\in[0,1]$ as
\begin{align*}
\wt Q_{\lambda}^{\dagger} = (1-\lambda)\,\wt Q^{(b)}|_{A_{n,j}^c} + \lambda \,\wt Q^{(b)}|_{A_{n,j}},
\end{align*}
where $A_{n,j}=\{\theta:\,|\theta_j-\theta^\ast_j|\geq D\varepsilon\}\times \m S^n$. In particular, by the optimality, $\wt \lambda= \wt Q^{(b)}(A_{n,j})=\wt Q^{(b)}_\theta(|\theta_j-\theta^\ast_j|\geq D\varepsilon)$ minimizes the weighted variational objective function $\wt D^{(b)}\big(\wt Q_\lambda^{\dagger}\,||\,\wt P^{(b)}(\cdot\,|\,X^n)\big)$ defined in~\eqref{eqn:vwlbObj} as a function of $\lambda$. In addition, we have a similar decomposition as
\begin{align*}
&\wt D^{(b)}\big(\wt Q_\lambda^{\dagger}\,||\,\wt P^{(b)}(\cdot\,|\,X^n)\big)=  (1 - \lambda)\, \wt d_{n2} + \lambda  \, \wt d_{n1}
+ D\big(\text{Ber}(\lambda)\, \| \, \text{Ber}(\wt\beta_n)\big),
\end{align*}
where $\wt \beta_n = \wt P^{(b)}(A_{n,j}\,|\,X^n) \leq \wt \Pi_n^{(b)}(\|\theta-\theta^\ast\|\geq D\varepsilon)$, 
$\text{Ber}(\wt \lambda)$ denotes a Bernoulli distribution with success probability $\wt \lambda$, and two constants $(d_{n1},d_{n2})$ independent of $\wt \lambda$ are
\begin{align*}
\wt d_{n1} = \wt D^{(b)}\big(\wt Q^{(b)}|_{A_{n,j}}\, \| \wt P^{(b)}|_{A_{n,j}}(\cdot\,|\,X^n)\big) \quad \mx{and} \quad \wt d_{n2} =  \wt D^{(b)}\big(\wt Q^{(b)}|_{A_{n,j}^c}\, \| \wt P^{(b)}|_{A_{n,j}^c}(\cdot\,|\,X^n)\big),
\end{align*}
where for any set $A\subset \Theta\times\m S^n$, $\wt p^{(b)}|_A$ (not a valid density function) is defined as 
\begin{align*}
\wt p^{(b)}|_A(\theta,\,S^n|\,X^n) =\frac{\wt p^{(b)}(X^n,S^n|\,\theta)\,\pi(\theta)}{\int_{A}\wt L^{(b)}(\theta;\,X^n) \,\pi(\theta)\,\dd\theta},\quad \mx{for all $(\theta, S^n)\in A$.}
\end{align*}
Due to the decomposition~\eqref{Eqn:WKL_decom}, each term in the decomposition remains nonnegative. Consequently, the rest steps in the proof of Lemma~\ref{lem:expDacayverify} still apply and it remains to 
prove a weighted version of Lemma~\ref{lem:expDacayverify} as $\wt D^{(b)}\big(\wt Q^{(b)}\,||\,\wt P^{(b)}(\cdot\,|\,X^n)\big)\leq CnM\wt \varepsilon_n^2$ holds with probability at least $1-CM^{-2}$. A proof for this again is almost the same as that of Lemma~\ref{lem:expDacayverify} by utilizing Lemma~\ref{lem:decomposeWKL}, the first inequality~\eqref{Eqn:contraction_WP} Theorem~\ref{thm:WMF_contraction} and the fact that the random weights $W_i^{(b)}$ has unit mean (when bounding the expectation as in~\eqref{Eqn:KL_mean}).


\subsection{Proof of Theorem~\ref{thm:BvM_wp}}\label{app:Proof_thm:BvM_wp}
We just sketch the proof about the bound on the KL divergence $D(\wt\Pi_n^{(b)}\,||\,\wt Q^{\ast(b)})$, and the bound on the total variational distance can be proceeded in a similar way (for a proof on the total variational distance in the unweighted case, which is the classical BvM theorem, c.f.~Chapter 1.4 of~\cite{ghoshbayesian}). In fact, the proof is simply a weighted extension of the proof of Lemma~\ref{lem:KL_normal_approx} with $Q_\theta$ chosen as the $\wt\Pi_n^{(b)}$ (by Theorem~\ref{thm:WMF_contraction} it satisfies the sub-Gaussian tail condition therein) and $\Pi_n$ being replaced with the weighted posterior $\wt\Pi_n^{(b)}$. The only difference in the proof is the Taylor expansion as~\eqref{eq:TylorLR} and~\eqref{Eqn:taylor_ex}, which in the weighted case they should be expanded at the weighed MLE $\wt\theta^{(b)}$ rather than the MLE $\mle$, as now $\wt\theta^{(b)}$ maximizes the log-weighted likelihood function $\wt l^{(b)}(\theta;\,X^n) = \log \wt L^{(b)}(\theta;\,X^n)$, whose gradient vanishes at $\wt\theta^{(b)}$, and the rest of the proofs are the same.


\subsection{Proof of Theorem~\ref{thm:3}}\label{app:Proof_thm:3}
With Theorem~\ref{thm:WMF_contraction}, the proof of this theorem is again similar to that of Theorem~\ref{thm:1} by adapting to the weighted case.
We just need to point out the difference. Lemma~\ref{lem:KL_identity} in step one remains unchanged. In step two, based on almost same lines of the proof (using the fact that the random weights have unit expectation in the last step of applying Markov inequality), the conclusion of Lemma~\ref{lem:LQA_VB} becomes
\begin{align*}
\Big|\wt F_n^{(b)}(Q_\theta) - D\big(Q_\theta\,||\,\wt \Pi^{(b)}_n\big) - \frac{n}{2}\,\text{tr}\big(\Sigma_{Q_\theta}\, I_s(\theta^\ast)\big)\Big| \leq \frac{CM^3(\log n)^{3}}{\sqrt{n}},
\end{align*}
where $\wt F_n^{(b)}(Q_\theta)$ denotes the weighted ``profile divergence" defined as as the right hand side of the identity in Lemma~\ref{lem:decomposeWKL}, and the exponent $3$ of the $\log n$ term is due to the fact that $\wt \varepsilon_n$ in Theorem~\ref{thm:WMF_contraction} is up to a constant $\sqrt{\log n}$ times larger than $\varepsilon_n$ in Lemma~\ref{lem:expDacayverify}. In step three, the approximation result in Lemma~\ref{lem:KL_normal_approx} becomes
\begin{align*}
\Big|D\big(Q_\theta\,||\,\wt \Pi^{(b)}_n\big) -D\big(Q_\theta\,||\,N\big(\wt\theta^{(b)}, \, [n I(\theta^\ast)]^{-1}\big)\big)\Big| \leq \frac{CM^3(\log n)^{d+4}}{\sqrt{n}},
\end{align*}
since now $\wt\theta^{(b)}$ plays the role of the MLE $\mle$ in the unweighted posterior case. Its proof is also almost the same as the proof of Lemma~\ref{lem:KL_normal_approx} plus the proof of Lemma~\ref{lem:deIntApprox} (a similar situation as the proof of Lemma~\ref{app:Proof_thm:BvM_wp}). Steps four and five remain valid for the weighted case as well, and we omit the details.


\subsection{Proof of Corollary~\ref{cor:WVB_center}}\label{app:Proof_cor:WVB_center}
The claimed bound is due to Theorem~\ref{thm:3} and the last identity in the proof of Lemma~\ref{lem:KL_identity}.


\subsection{Proof of Theorem~\ref{thm:VWLB_consistency}}\label{app:Proof_thm:VWLB_consistency}
We only need to prove the first inequality, since the second can be obtained by combining the first with the classical BvM theorem under the total variation metric (c.f.~Chapter 1.4 of~\cite{ghoshbayesian}).
An asymptotic version of the first inequality was proved in~\cite{wlb} (Theorem 2), we provide a proof for our non-asymptotic version, which proceeds as follows. 

Let $s(\theta;X^n) := \nabla l(\theta;\, X^n)=\sum_{i=1}^n \nabla \log p(X_i\,|\,\theta)$ be the (unweighted) score function. By the optimality of the MLE $\mle$, we have $s(\mle; X^n) = 0$. The classical analysis of the MLE (for example,~Chapter 1.4 of~\cite{ghoshbayesian}) implies that under Assumption A2
\begin{equation}\label{eqn:mleexpansion}
\sqrt{n}(\mle - \theta^\ast) = [I(\theta^\ast)]^{-1}\frac{1}{\sqrt{n}}\,s(\theta^\ast; X^n) + R_n,
\end{equation}
where the remainder term $R_n$ satisfies $P(|R_n|\geq C_0M(\log n)^{3/2}/\sqrt{n}) \leq 1-CM^{-2}$ for any $M\geq 1$.
A similar argument based on Taylor expansion can be applied to the weighted MLE, yielding that under Assumptions A2 and W,
\begin{equation}\label{eqn:wexpansion}
\sqrt{n}(\wt{\theta}^{(b)} - \theta^\ast) = [I(\theta^\ast)]^{-1}\frac{1}{\sqrt{n}}\,\wt{s}(\theta^\ast; X^n) + \wt{R}_n.
\end{equation}
where $\wt{s}(\theta;X^n) := \nabla \wt l^{(b)}(\theta;\, X^n) = \sum_{i=1}^n W_i^{(b)} \nabla \log p(X_i\,|\,\theta)$ is the weighted score function and $\wt R_n$ is a remainder term also satisfies $P(|\wt R_n|\geq C_0M(\log n)^{3/2}/\sqrt{n}) \leq 1-CM^{-2}$ for any $M\geq 1$. By taking the difference between these two, we reach
\begin{align*}
\sqrt{n}(\wt{\theta}^{(b)} -\mle) = [I(\theta^\ast)]^{-1}\frac{1}{\sqrt{n}}\,\sum_{i=1}^n (W_i^{(b)}-1)\, \nabla \log p(X_i\,|\,\theta) + R_n-\wt{R}_n.
\end{align*}
By combining this with Corollary~\ref{cor:WVB_center}, we can get
\begin{align}\label{Eqn:weighted_mle_bound}
\sqrt{n}(\wt{\theta}^{(b)}_{VB} -\mle) = [I(\theta^\ast)]^{-1}\frac{1}{\sqrt{n}}\,\sum_{i=1}^n (W_i^{(b)}-1)\, \nabla \log p(X_i\,|\,\theta) + \wt R_{VB},
\end{align}
where the remainder term $\wt R_{VB}$ satisfies $P(|\wt R_{VB}| \geq C'M^{3/2} (\log n)^{d/2+2}n^{-1/4})\leq 
1-CM^{-2}$. We use the shorthand $\zeta_n = [I(\theta^\ast)]^{-1}\frac{1}{\sqrt{n}}\,\sum_{i=1}^n (W_i^{(b)}-1)\, \nabla \log p(X_i\,|\,\theta)\in\mathbb R^d$. Under Assumptions A2 and W, by elementwisely applying the Markov inequality and a union bound argument, we obtain
\begin{align}\label{Eqn:cov_bound}
P_{\theta^\ast}\Big(\opnorm{\mbox{Cov}_W[\zeta_n\,|\,X^n] - [I(\theta^\ast)]^{-1}}\leq \frac{CM}{\sqrt{n}}\Big) \geq 1-CM^{-2},
\end{align}
where $\mbox{Cov}_W[\zeta_n\,|\,X^n]$ denotes the conditional covariance matrix of $\zeta_n$ given $X^n$, and we have used the fact that the matrix Frobenius norm is at most $\sqrt{d}$ times larger than the matrix operator norm (the constant $C$ may depend on the dimension $d$). Therefore, under this high probability event, $\mbox{Cov}_W[\zeta_n\,|\,X^n]$ is non-singular. In the rest of the proof we always work under this event. By combining~\eqref{Eqn:weighted_mle_bound}, \eqref{Eqn:cov_bound}, Assumption A2 and a Markov inequality for the sum $\zeta_n$, we can get that
\begin{align*}
\big\|\sqrt{n}[I(\theta^\ast)]^{1/2}(\wt{\theta}^{(b)}_{VB} -\mle) - \big\{\mbox{Cov}_W[\zeta_n\,|\,X^n]\big\}^{-1/2} \zeta_n\big\|_\infty \leq CM^{3/2} (\log n)^{d/2+2}n^{-1/4}=:\kappa_n,
\end{align*}
holds with probability at least $1-CM^{-2}$. This inequality can also be written as (recall that the notation $a\leq b$ for two vectors $a$ and $b$ means element-wise less than or equal to)
\begin{align*}
 \big\{\mbox{Cov}_W[\zeta_n\,|\,X^n]\big\}^{-1/2} \zeta_n - \kappa_n e_d\leq \sqrt{n}[I(\theta^\ast)]^{1/2}(\wt{\theta}^{(b)}_{VB} -\mle) \leq \big\{\mbox{Cov}_W[\zeta_n\,|\,X^n]\big\}^{-1/2} \zeta_n + \kappa_n e_d
\end{align*}
where $e_d$ denotes the all one vector in $\mb R^d$.

By applying the multivariate Berry-Esseen theorem (for example, Theorem 1 in~\cite{raivc2019multivariate}) with set $A$ therein restricted to all convex sets of the form $\{x\in\mb R^d:\,x\leq u\}$ over all $u\in \mb R^d$, we obtain
\begin{align*}
\sup_{u\in\mb R^d}\Big|P \big(\sqrt{n}\,\big\{\mbox{Cov}_W[\zeta_n\,|\,X^n]\big\}^{-1/2} \zeta_n\leq u\,\big|\,X^n\big) - P(Z \leq u) \Big|& \leq \frac{C}{n^{3/2}}\sum_{i=1}^n \|\nabla \log p(X_i\,|\,\theta) \|^3,
\end{align*}
where $Z\sim N(0,I_d)$.  Note that by applying Markov inequality and Assumption A2, the right hand side is upper bounded by $CM/\sqrt{n}$ with probability at least $1-CM^{-2}$. Since the pdf of $Z\sim N(0, I_d)$ is uniformly bounded from above, we have 
\begin{align*}
\sup_{u\in\mb R^d}\big|P(Z\leq u + \kappa_n e) - P(Z\leq u)\big|\leq C\kappa_n, \quad \sup_{u\in\mb R^d}\big|P(Z\leq u - \kappa_n e) - P(Z\leq u)\big|\leq C\kappa_n.
\end{align*}
 Finally, by combining the last three displays and the fact that $P(Z\leq u + te_d)$ is monotonically increasing in $t>0$, we obtain that
 \begin{align*}
 \sup_{u\in\mb R^d}\Big|P \big(\sqrt{n}[I(\theta^\ast)]^{1/2}(\wt{\theta}^{(b)}_{VB} -\mle)\leq u\,\big|\,X^n\big) - P(Z \leq u) \Big|& \leq C(\kappa_n +M/\sqrt{n})
 \end{align*}
holds with probability at least $1-CM^{-2}$, which completes the proof.


\subsection{Proof of Lemma~\ref{lem:KL_bound}}\label{app:lem:KL_bound_proof}
To provide the claimed bound, we will use the optimality of $\wht Q_\theta$ for minimzing the profile divergence $F_n(Q_\theta)$ defined in step two in the proof of Theorem~\ref{thm:1} in Section~\ref{Sec:proof_thm1}, and compare it with $\widebar Q_\theta$, the minimizer of the KL-divergence to the marginal posterior $\Pi_n$ in the mean-field family, or
\begin{align*}
\widebar Q_\theta = \argmin_{Q_\theta=\bigotimes_{j=1}^k Q_{\theta_j}} D\big(Q_\theta\,||\,\Pi_n\big).
\end{align*}
More specifically, by the definition of $F_n$ and the optimality of $\wht Q_\theta$, we have
\begin{align*}
D\big(\, \wht Q \,||\, P(\cdot\, ||\,X^n)\big) = F_n(\wht Q_\theta) \leq F_n(\widebar Q_\theta).
\end{align*}

\vspace{0.5em}
\noindent{\bf Step one:} First, we show that $\widebar Q_\theta$ satisfies the sub-Gaussian tail decay as in Lemma~\ref{lem:expDacayverify}, so that we can apply Lemma~\ref{lem:LQA_VB} to bound $F_n(\widebar Q_\theta)$ in step two below. In fact, the same variational argument as in the proof of Lemma~\ref{lem:expDacayverify} in Section~\ref{Sec:Proof_lem:expDacayverify} leads to (by using the optimality of $\bar Q_\theta$ as the minimizer of $D\big(Q_\theta\,||\,\Pi_n\big)$)
\begin{align*}
\widebar Q_\theta (\|\theta-\theta^\ast\| \geq D\varepsilon) =  \widebar {\lambda} \leq 2\beta_n \exp(\frac{D\big(\widebar Q_\theta\, \| \,\Pi_n)}{1 - \widebar{\lambda}}) \quad \mx{and}\quad  D\big(\text{Ber}(\widebar \lambda)\, \| \, \text{Ber}(\beta_n)\big) \leq D\big(\widebar Q_\theta\, \| \,\Pi_n),
\end{align*}
where recall that $\beta_n=\Pi_n(\|\theta-\theta^\ast\|\geq D\varepsilon) \leq e^{-CD^2n\varepsilon^2}\leq 1/2$. Similar to the argument in the last paragraph in Section~\ref{Sec:Proof_lem:expDacayverify}, in order to prove the sub-Gaussian tail bound as $\widebar Q_\theta(\|\theta-\theta^\ast\|\geq D\varepsilon) \leq \exp\big(-CD^2n\varepsilon^2/2\big)$, it suffices to show that $D\big(\widebar Q_\theta\, \| \,\Pi_n)\leq CnM^2\varepsilon_n^2$. The rest of step one devotes to the proof of this bound. We make use of the following identity for any $Q_\theta$ from~\cite{zhang2006,Chao2018},
\begin{align*}
E_{\theta^\ast}\big[ D\big(Q_\theta\,||\,\Pi_n\big)\big]& = D\big(Q_\theta\,||\,\Pi\big) + E_{Q_\theta}\Big[ D_{X^n}\big(p(X^n|\,\theta^\ast)\,||\,p(X^n|\,\theta)\big) - D_{X^n}\big(p(X^n|\,\theta^\ast)\,||\,p(X^n)\big) \Big]\\
&\leq D\big(Q_\theta\,||\,\Pi\big) + E_{Q_\theta}\Big[ D_{X^n}\big(p(X^n|\,\theta^\ast)\,||\,p(X^n|\,\theta)\big)\Big],
\end{align*}
where we use $D_{X^n}$ to denote (indicate) the KL divergence whose integral is with respect to $X^n$, not $\theta$, and recall $p(X^n) = \int_{\Theta} p(X^n|\,\theta)\,\pi(\theta)\,\dd\theta$ is the normalization constant (viewed as a function of $X^n$) in the posterior $\Pi_n$. Due to the optimality of $\widebar Q_\theta$, we can bound $E_{\theta^\ast}\big[D\big(\widebar Q_\theta\, \| \,\Pi_n)\big]$ by the above with any $Q_\theta$ in the mean-field family. In particular, we choose $Q_\theta = Q^\diamond_\theta$, where $Q^\diamond_\theta$ is the uniform distribution over the $\varepsilon_n$-$\ell_\infty$ ball $B_\infty(\theta^\ast;\,\varepsilon_n)$ at $\theta^\ast$, whose density function is
\begin{align*}
q^\diamond_\theta=\begin{cases}
\big[\mx{Vol}(B_\infty(\theta^\ast;\,\varepsilon_n))\big]^{-1} & \mx{if $\|\theta-\theta^\ast\|_\infty \leq \varepsilon_n$},\\
0 & \mx{otherwise}.
\end{cases}
\end{align*}
It is straightforward to see that $Q^\diamond_\theta$ belongs to the mean-field family.
Due to the continuity of $\log \pi(\theta)$ around $\theta^\ast$ by Assumption A1, we have
\begin{align*}
D\big(Q_\theta^\diamond\,||\,\Pi\big) &= \big[\mx{Vol}(B_\infty(\theta^\ast;\,\varepsilon_n))\big]^{-1}\int_{B_\infty(\theta^\ast;\,\varepsilon_n)} \big( -d \log (2\varepsilon_n) - \log \pi(\theta)\big)\,\dd\theta\\
& \leq \sup_{\theta\in B_\infty(\theta^\ast;\,\varepsilon_n)} \big|\,d \log (2\varepsilon_n) + \log \pi(\theta)\big|\leq C \log n.
\end{align*}
By Assumption A2 on the derivatives of log-likelihood function, we have 
\begin{align*}
 E_{Q_\theta^\diamond}\Big[ D_{X^n}\big(p(X^n|\,\theta^\ast)\,||\,p(X^n|\,\theta)\big)\Big]& 
 = n\,E_{Q_\theta^\diamond}\Big[ D_{X_1}\big(p(X_1|\,\theta^\ast)\,||\,p(X_1|\,\theta)\big)\Big] \\
&
 \leq Cn\, E_{Q_\theta^\diamond}\big[\|\theta-\theta^\ast\|^2\big] \leq C d\,n\,\varepsilon_n^2,
\end{align*}
where we have used the conditional independence of $\{X_i\}_{i=1}^n$ given $\theta$ and the last step is due to the inequality $\|\theta-\theta^\ast\|\leq \sqrt{d}\, \|\theta-\theta^\ast\|_\infty$.
Putting pieces together, we can reach 
\begin{align}\label{Eqn:KL_mean}
E_{\theta^\ast}\big[D\big(\widebar Q_\theta\, \| \,\Pi_n)\big] \leq D\big(Q_\theta^\diamond\,||\,\Pi\big) + 
E_{Q_\theta^\diamond}\Big[ D_{X^n}\big(p(X^n|\,\theta^\ast)\,||\,p(X^n|\,\theta)\big)\Big] \leq C\log n + C\,n\varepsilon_n^2 = C' n\varepsilon_n^2.
\end{align}
Finally, an application of the Markov inequality implies the following to hold with probability at least $1-CM^{-2}$,
\begin{align*}
D\big(\widebar Q_\theta\, \| \,\Pi_n) \leq Cn M^2 \varepsilon_n^2,
\end{align*}
which finished the proof of the sub-Gaussian tail bound for $\widebar Q_\theta$.

\vspace{0.5em}
\noindent{\bf Step two:} Due to the sub-Gaussian tail bound of $\widebar Q_\theta$ and the preceding display, we may apply Lemma~\ref{lem:LQA_VB} to conclude that
\begin{align*}
&D\big(\, \wht Q \,||\, P(\cdot\, ||\,X^n)\big) = F_n(\wht Q_\theta) \leq F_n(\widebar Q_\theta)\\
&\qquad
\leq  D\big(\widebar Q_\theta\, \| \,\Pi_n) + \frac{n}{2}\mx{tr}\big(\Sigma_{\widebar Q_\theta} I_s(\theta^\ast)\big) + \frac{CM^3(\log n)^{3/2}}{\sqrt n }
\leq Cn M^2 \varepsilon_n^2
\end{align*}
where the last step is due to the bound as $E_{\widebar Q_\theta}\big[\|\theta-\theta^\ast\|^2\big] \leq CM^2\varepsilon_n^2$ (due to the sub-Gaussian tail, c.f.~inequality~\eqref{Eqn:moment_bound} with $k=2$ for a proof) for the middle term.


\subsection{Proof of Lemma~\ref{lem:KL_identity}}\label{app:Proof_KL_identity}

Without loss of generality, we may assume $Q$ to admit a density function, denoted by $q$, with respect to the $d$-dim Lebesgue measure (otherwise both sides are equal to infinity and the claimed results hold).
By direct calculation and using the pdf of a multivariate normal distribution, we have
\begin{equation}\label{Eqn:KL_expression}
\begin{aligned}
&D\big(Q\,||\,N(\mu,\,\Gamma^{-1})\big) = \int_{\mb R^d} q(\theta)\log q(\theta)\, \dd\theta 
+\frac{1}{2} \log\big((2\pi)^{d} |\Gamma|^{-1}\big) \\
&\qquad\qquad\qquad\qquad\qquad\qquad + \frac{1}{2}\int_{\mb R^d} q(\theta) \Big(\sum_{j=1}^d \Gamma_{jj} (\theta_{j}-\mu_j)^2 + \sum_{j\neq k} \Gamma_{jk} (\theta_j-\mu_j)(\theta_k-\mu_k)\Big)\,\dd\theta\\
&= \sum_{j=1}^d\int_{\mb R} q_j(\theta_j)\log q_j(\theta_j)\, \dd\theta_j 
+\frac{1}{2} \log\big((2\pi)^{d} |\Gamma|^{-1}\big) + \frac{1}{2}\sum_{j=1}^d \Gamma_{jj} \text{Var}_{Q_j}[\theta_{j}] + \frac{1}{2}(\mu_{Q_\theta}-\mu)^T \Gamma\,(\mu_{Q_\theta}-\mu).
\end{aligned}
\end{equation}
By substituting $Q$ with $Q^\ast=N\big(\mu,(\text{diag}(\Gamma))^{-1}\big)$ in this display, we obtain
\begin{align*}
D\big(Q^\ast\,||\,N(\mu,\,\Gamma^{-1})\big) &= \sum_{j=1}^d\int_{\mb R} q^\ast_j(\theta_j)\log q^\ast_j(\theta_j)\, \dd\theta_j 
+\frac{1}{2} \log\big((2\pi)^{d} |\Gamma|^{-1}\big) + \frac{1}{2}\sum_{j=1}^d \Gamma_{jj} \text{Var}_{Q^\ast_j}[\theta_{j}]\\
&=\frac{1}{2} \log\big((2\pi)^{d} |\Gamma|^{-1}\big) - \frac{1}{2} \log\big((2\pi)^{d} |\text{diag}(\Gamma)|^{-1}\big)
\end{align*}
Taking the difference between the two preceding displays and using the translation invariance of $\int q_j(\theta)\log q_j(\theta_j)\,\dd\theta_j$ for each $j$, we can reach
\begin{align*}
&D\big(Q\,||\,N(\mu,\,\Gamma^{-1})\big)-D\big(Q^\ast\,||\,N(\mu,\,\Gamma^{-1})\big)=
\sum_{j=1}^d\int_{\mb R} q_{(\mu_j),j}(\theta_j)\log q_{(\mu_j),j}(\theta_j)\, \dd\theta_j \\
&\qquad\qquad\qquad\qquad\qquad
+\frac{1}{2}\sum_{j=1}^d \log(2\pi \Gamma_{jj}^{-1})+  \frac{1}{2}\sum_{j=1}^d \Gamma_{jj} \text{Var}_{Q_j}[\theta_{j}] + \frac{1}{2}(\mu_{Q_\theta}-\mu)^T \Gamma\,(\mu_{Q_\theta}-\mu)\\
&\qquad\qquad\qquad\qquad\qquad\qquad\qquad
= \sum_{j=1}^d
D\big(Q_{(\mu_j),j}\,||\,Q^\ast_j\big) + \frac{1}{2}(\mu_{Q_\theta} - \mu)^T \Gamma\,(\mu_{Q_\theta} - \mu),
\end{align*}
where the last step is due to an application of \eqref{Eqn:KL_expression} to each $D\big(Q_{(\mu_j),j}\,||\,N(\mu_j,\Gamma_{jj}^{-1})\big)$ (they have the same expectation $\mu_j$). This completes the proof of the first identity.

Now we apply the first identity in the lemma with $\Gamma$ being $\text{diag}(\Gamma)$ and $Q^\ast=N(\mu,\text{diag}(\Gamma))$ so that the second KL-divergence term becomes zero,
\begin{align*}
D\big(Q\,||\,Q^\ast\big)= \sum_{j=1}^d
D\big(Q_{(\mu_j),j}\,||\,Q^\ast_j\big) + \frac{1}{2}(\mu_{Q_\theta} - \mu)^T \text{diag}(\Gamma)\,(\mu_{Q_\theta} - \mu), 
\end{align*}
The second desired inequality follows by comparing this identity with the first identity in the lemma and using the fact that KL divergence is nonnegative.


\subsection{Proof of Lemma~\ref{lem:LQA_VB}}\label{app:Proof_LQA_VB}
For any fixed $Q_\theta$, due to the sub-Gaussian tail bound $Q_\theta(\|\theta-\theta^\ast\|\geq C\varepsilon) \leq e^{-Cn\varepsilon^2}$ for all $\varepsilon\geq M\varepsilon_n$, we can apply the identity $\int_{0}^\infty x \,P(\dd x) = \int_{0}^\infty P(X\geq x) \dd x$ for any probability measure over $[0,\infty)$ to bound the $k$th ($k\geq 1$) order moment of $\|\theta-\theta^\ast\|$ as
\begin{equation}\label{Eqn:moment_bound}
\begin{aligned}
m_{Q_\theta,k}&:\,=E_{Q_\theta}\|\theta-\theta^\ast\|^k = \int_{0}^\infty kt^{k-1} Q_\theta(\|\theta-\theta^\ast\|\geq t) \,\dd t \leq C D^kM^k\varepsilon_n^k + \int_{D\varepsilon_n}^\infty t^{k-1} e^{-Cnt^2} \,\dd t\\
&\leq CD^kM^k\varepsilon_n^k + Ce^{-CD^2n\varepsilon_n^2} \leq  CD^kM^k\varepsilon_n^k +Cn^{-D^2C}\leq CD^kM^k \varepsilon_n^k,
\end{aligned}
\end{equation}
by choosing a sufficiently large constant $D$. We will repeatedly use this moment bound for $k=1,2,3$ throughout the proof.

Due to Lemma~\ref{lem:decomposeKL}, it suffices to show 
\begin{align*}
\Big|\sum_{i=1}^n \log \Big(\int_{\m S} \exp\{r_i(s_i)\}\,\dd s_i\Big) + \frac{n}{2}\,\text{tr}\big(\Sigma_{Q_\theta}\, I_s(\theta^\ast)\big)\Big|\leq  \frac{CM^3(\log n)^{3/2}}{\sqrt{n}}.
\end{align*}
The proof of this result is based on tedious calculations via Taylor expansions.
To simplify the presentation of the proof, we assume $S_i$ to be discrete so that all integrals over $\m S$ reduce to summations, and $\theta$ is one-dimensional. According to the definition of $r_i(s)$, we can rewrite
\begin{equation*}
\sum_{s \in \mathcal{S}} \exp(r_i(s))
= \sum_{s \in \mathcal{S}} p(s_i = s\,| \,X_i, \mu_{Q_\theta})\, \exp\Big\{ \int q_\theta(\theta) \log \frac{p(s_i = s\,|\, X_i, \theta)}{p(s_i = s\,|\, X_i,  \mu_{Q_\theta})} \,\dd \theta\Big\},
\end{equation*}
where recall that $\mu_{Q_\theta}$ denotes the expectation of a probability measure $Q$. In the following we consider $d = 1$ case for ease of notation, while the proof can be generalized to other fixed dimensions trivially. We use the shorthand $p_{is}(\theta)$ to denote $p(s_i=s\,|\,X_i,\theta)$ for $i\in[n]$ and $s\in \m S$, $l_{is}(\theta)$ to denote $\log p_{is}(\theta)$, and $\big(l'_{is}(\theta),l''_{is}(\theta),l^{(3)}_{is}(\theta)\big)$ to denote its derivatives up to order three.

Consider a intermediate term 
\begin{equation*}
R_{sn}=
-\sum_{i = 1}^n \log \Big\{\sum_{s \in \mathcal{S}} \exp(r_i(s))\Big\} +\frac{1}{2} \Sigma_{Q_\theta}
\sum_{i = 1}^n  \sum_{s \in \mathcal{S}} p_{is}(\mu_{Q_\theta}) \, l''_{is}(\mu_{Q_\theta}).
\end{equation*}
By definition $-\sum_{s \in \mathcal{S}} p_{is}(\theta^\ast)\, l''_{is}(\theta^\ast)$ has expectation $I_s(\theta^\ast)$. Therefore, by using the Taylor expansion, Assumption A3, the moment bound~\eqref{Eqn:moment_bound} with $k=2,3$ and a Markov inequality for sum of i.i.d.~random variables, we obtain that with probability at least $1-CM^{-2}$ that
\begin{align*}
&\Big|\, \frac{1}{2}\, 
 \Sigma_{Q_\theta} \sum_{i = 1}^n  \sum_{s \in \mathcal{S}} p_{is}(\mu_{Q_\theta})\, l''_{is}(\mu_{Q_\theta}) + \frac{n}{2}\, \text{tr}\big[\Sigma_{Q_\theta}\,I_s(\theta^\ast)\big]\Big| \\
 &\qquad\qquad\qquad\qquad\qquad\qquad
 \leq C\,n\,\big(m_{Q_\theta,2}\,M\sqrt{\frac{\log n}{n}} + m_{Q_\theta,3}\big) \leq  \frac{CM^3(\log n)^{3/2}}{\sqrt{n}}.
\end{align*}
Therefore, it remains to show $|R_{sn}| \leq \frac{CM^3(\log n)^{3/2}}{\sqrt{n}}$.

We apply the Taylor expansion up to the third order to $l_{is}(\theta)$ at $\theta=\mu_{Q_\theta}$ in the preceding display,
\begin{align*}
& \sum_{s \in \mathcal{S}} \exp(r_i(s)) =\sum_{s \in \mathcal{S}} p_{is}(\mu_{Q_\theta}) \exp\Big\{ \int q_\theta(\theta) \log \frac{p_{is}(\theta)}{p_{is}(\mu_{Q_\theta})} \dd \theta\Big\}  \\
= & \sum_{s \in \mathcal{S}} p_{is}(\mu_{Q_\theta}) \exp\Big\{ \int q_\theta(\theta) \bigl\{ 
(\theta - \mu_{Q_\theta}) \, l'_{is}(\mu_{Q_\theta}) +
 (\theta - \mu_{Q_\theta})^2 \,l''_{is}(\mu_{Q_\theta})/2 + 
(\theta - \mu_{Q_\theta})^3\, l^{(3)}_{is}(\theta^\dagger)/6
 \bigr\} \dd \theta\Big\} \\
= & \sum_{s \in \mathcal{S}} p_{is}(\mu_{Q_\theta}) \exp \Big\{
\Sigma_{Q_\theta} \, l''_{is}\,(\mu_{Q_\theta}) /2 +
m_{Q_\theta,3}\,  R_i \Big\},
\end{align*}
where $\theta^\dagger$ is some point between $\theta$ and $\mu_{Q_\theta}$, $R_i$ is a remainder term, and recall that $m_{Q_\theta,3}=E_{Q_\theta}\|\theta-\theta^\ast\|^3$.
Here, the last step is due to the fact that $\mu_{Q_\theta}$ is the expectation of $Q_\theta$ so that the integral of the first linear term in the second line vanishes. Due to Assumption A3, the remainder term satisfies $|R_i|\leq M(s,X_i)/6$.
As a consequence, we obtain
\begin{align*}
&|R_{sn}| = 
\biggl|- \sum_{i = 1}^n \log \Big\{\sum_{s \in \mathcal{S}} \exp(r_i(s))\Big\} + \frac{1}{2}\,\Sigma_{Q_\theta}
\sum_{i = 1}^n  \sum_{s \in \mathcal{S}} p_{is}(\mu_{Q_\theta})\, l''_{is}(\mu_{Q_\theta}) \biggr|\\
= &\, \biggl| -\sum_{i = 1}^n \log \Big\{\sum_{s \in \mathcal{S}} p_{is}(\mu_{Q_\theta}) \exp \Big\{
\frac{1}{2}\,\Sigma_{Q_\theta} \, l''_{is}\,(\mu_{Q_\theta}) +
m_{Q_\theta,3}\,  R_i \Big\}\Big\}
+\frac{1}{2}\Sigma_{Q_\theta} \sum_{i = 1}^n  \sum_{s \in \mathcal{S}} p_{is}(\mu_{Q_\theta})\, l''_{is}(\mu_{Q_\theta})\biggr|.
\end{align*}
Now we apply inequalities $|\log(1 + \delta)| \leq |\delta| + 2\delta^2$, $|e^\delta - 1| \leq |\delta| + 2\delta^2$ for $|\delta|\in[0,1/2]$, and the identity $\sum_{s\in \m S} p_{is}(s)=1$ for any $i\in[n]$ to obtain 
\begin{align*}
|R_{sn}| \leq &\, \biggl|
\sum_{i=1}^n \big\{ 1 - \sum_{s \in \mathcal{S}} 
p_{is}(\mu_{Q_\theta}) \exp \bigl\{ \frac{1}{2}\,\Sigma_{Q_\theta} \, l''_{is}\,(\mu_{Q_\theta}) +
m_{Q_\theta,3}\,  R_i  \bigr\} \bigr\}
+ \frac{1}{2}\Sigma_{Q_\theta} \sum_{i = 1}^n  \sum_{s \in \mathcal{S}} p_{is}(\mu_{Q_\theta})\, l''_{is}(\mu_{Q_\theta})
\biggr| \\
&\, + 
2 \sum_{i=1}^n \biggl| 1 - \sum_{s \in \mathcal{S}} 
p_{is}(\mu_{Q_\theta}) \exp \bigl\{ \frac{1}{2}\,\Sigma_{Q_\theta} \, l''_{is}\,(\mu_{Q_\theta}) +
m_{Q_\theta,3}\,  R_i \bigr\}
\biggr|^2 \\
\leq &\, \biggl|
- \sum_{i=1}^n \sum_{s \in \mathcal{S}} 
p_{is}(\mu_{Q_\theta}) \bigl\{ \frac{1}{2}\,\Sigma_{Q_\theta} \, l''_{is}\,(\mu_{Q_\theta}) +
m_{Q_\theta,3}\,  R_i  \bigr\} 
+\frac{1}{2}\Sigma_{Q_\theta} \sum_{i = 1}^n  \sum_{s \in \mathcal{S}} p_{is}(\mu_{Q_\theta})\, l''_{is}(\mu_{Q_\theta})
\biggr| \\
\,& + 
2 \sum_{i=1}^n \biggl| 1 - \sum_{s \in \mathcal{S}} 
p_{is}(\mu_{Q_\theta}) \exp \bigl\{ \frac{1}{2}\,\Sigma_{Q_\theta} \, l''_{is}\,(\mu_{Q_\theta}) +
m_{Q_\theta,3}\,  R_i  \bigr\}
\biggr|^2 \\
\,& + 2\sum_{i=1}^n  \sum_{s \in \mathcal{S}} 
p_{is}(\mu_{Q_\theta}) \biggl|\frac{1}{2}\,\Sigma_{Q_\theta} \, l''_{is}\,(\mu_{Q_\theta}) +
m_{Q_\theta,3}\,  R_i 
\biggr|^2 \\
=\, &  \bigl| \sum_{i=1}^n \sum_{s \in \mathcal{S}} p_{is}(\mu_{Q_\theta})\, m_{Q_\theta,3}\,  R_i \bigr| + 
2 \sum_{i=1}^n \biggl| 1 - \sum_{s \in \mathcal{S}} 
p_{is}(\mu_{Q_\theta}) \exp \bigl\{\frac{1}{2}\,\Sigma_{Q_\theta} \, l''_{is}\,(\mu_{Q_\theta}) +
m_{Q_\theta,3}\,  R_i  \bigr\}
\biggr|^2 \\
\,& + 2\sum_{i=1}^n  \sum_{s \in \mathcal{S}} 
p_{is}(\mu_{Q_\theta}) \biggl|\frac{1}{2}\,\Sigma_{Q_\theta} \, l''_{is}\,(\mu_{Q_\theta}) +
m_{Q_\theta,3}\,  R_i  
\biggr|^2 \\
\leq \,& \bigl| \sum_{i=1}^n \sum_{s \in \mathcal{S}} p_{is}(\mu_{Q_\theta})\,m_{Q_\theta,3}\,  R_i\bigr| + 
6\sum_{i=1}^n  \sum_{s \in \mathcal{S}} p_{is}(\mu_{Q_\theta}) \biggl|
\frac{1}{2}\,\Sigma_{Q_\theta} \, l''_{is}\,(\mu_{Q_\theta}) +
m_{Q_\theta,3}\,  R_i  
\biggr|^2.
\end{align*}
Using the last bound, we can now invoke the moment bound~\eqref{Eqn:moment_bound} with $k=2,3$ again and a Markov inequality on $\sum_i|R_i|\leq \sum M(s,X_i)/6$ to obtain that with probability at least $1 - CM^{-2}$, 
$|R_{sn}| \leq \frac{CM^3(\log n)^{3/2}}{\sqrt{n}}$. This completes the proof.


\subsection{Proof of Lemma~\ref{lem:KL_normal_approx}}\label{app:proof_KL_normal_approx}
In the proof we omit the $\|\theta-\theta^\ast\|^L$ terms in Assumption A1 and A2 for ease of read, but the proof can be generalized trivially. Use $l(\theta:\, X^n)=p(X^n\,|\,\theta)$ to denote the log-likelihood function, and $\phi_n$ to denote the density function of $N\big(\mle,[nI(\theta^\ast)]^{-1}\big)$. Then the posterior density function of $\theta$ can be expressed as
\begin{align}\label{Eqn:posterior_form}
\pi_n(\theta) = \frac{\pi(\theta)\,\exp\big\{l(\theta; X^n) - l(\mle; X^n)\big\}}{\int \pi(\theta)\,\exp\big\{l(\theta; X^n) - l(\mle; X^n)\big\}\,\dd\theta},\quad \forall \theta\in\mb R^d,
\end{align}
where we have divided both the numerator and denominator by $\exp\big\{l(\mle; X^n)\big\}$ for technical convenience. We invoke the following lemma that provides an approximation to the integral in the denominator. A proof is provided in Appendix~\ref{app:proof_deIntApprox}.

\begin{lemma}\label{lem:deIntApprox}
Under Assumption A1 and A2, for any $M\geq 1$, the denominator in equation~\eqref{Eqn:posterior_form} satisfies
\begin{equation}\label{Eqn:normalization}
\bigg|\frac{ \int \pi(\theta)\,\exp\big\{l(\theta; X^n) - l(\mle; X^n)\big\}\,\dd\theta}{\displaystyle\pi(\theta^\ast) \big(2\pi/n\big)^{d/2}\,|I(\theta^\ast)|^{-1/2}} - 1 \bigg| \leq \frac{CM^3\,(\log n)^{d+3}}{\sqrt{n}},
\end{equation}
with probability at least $1-CM^{-2}$.
\end{lemma}

Denote the quantity on the right hand side of~\eqref{Eqn:normalization} by $\alpha_n\in(0,1/2]$.
Now we can express the log ratio between $\pi_n$ and $\phi_n$ as
\begin{equation}\label{Eqn:taylor_ex}
\begin{aligned}
\log(\frac{\pi_n}{\phi_n})(\theta) =&\, \log(\pi(\theta)\exp(l(\theta; X^n) - l(\mle; X^n))) \\
&\, -\log(\int \pi(\theta)\exp(l(\theta; X^n) - l(\mle; X^n)) \dd \theta) \\
&\, +\log(\frac{|nI(\theta^\ast)|^{-1/2}}{(2\pi)^{d/2}}) - \frac{n}{2}(\theta - \mle)^TI(\theta^\ast)(\theta - \mle)\\
=&\,l(\theta; X^n) - l(\mle; X^n) - \frac{n}{2}(\theta - \mle)^TI(\theta^\ast)(\theta - \mle)\\
&\, +  \log(\pi(\theta)) - \log(\pi(\theta^\ast)) + \log\big(1 + \alpha_n|I(\theta^\ast)|^{1/2}\big).
\end{aligned}
\end{equation}
Now we divide the space into two regions $K_n=\{\theta:\,\|\theta-\theta^\ast\| \leq c_0 M\varepsilon_n\}$ and $K_n^c=\mb R^d\setminus K_n$ for some sufficiently large constant $c_0$. For the region $K_n$, similar to the Taylor expansion in equation~\eqref{eq:TylorLR}, we may apply Assumption A1 and A2, and the preceding display to obtain that
\begin{align}\label{Eqn:bound_K_n}
\sup_{\theta\in K_n}\Big|\log(\frac{\pi_n}{\phi_n})(\theta)\Big| \leq \frac{Cc_0^3M^3 (\log n)^3}{\sqrt{n}} + C\alpha_n \leq \frac{Cc_0^3M^3\,(\log n)^{d+3}}{\sqrt{n}},
\end{align}
holds with probability at least $1-CM^{-2}$.

For the region $K_n^c$, we have by Assumption A1 and A2, and the expression of $\log\big(\pi_n/\phi_n\big)(\theta)$ that
\begin{align*}
\int_{K_n^c} \log(\frac{\pi_n}{\phi_n})(\theta)\, \dd Q_\theta(\theta)\leq  C\int_{K_n^c} \big(1 + n\|\theta - \theta^\ast\|^3\big)\, \dd Q_\theta(\theta) + C\alpha_n.
\end{align*}
Now since by the condition of the lemma, $Q_\theta$ satisfies $Q_\theta(\|\theta-\theta^\ast\|\geq C\varepsilon)\leq e^{-Cn\varepsilon^2}$ for all $\varepsilon\geq M\varepsilon_n$, we can bound the preceding display as
\begin{align*}
&\int_{K_n^c} \log(\frac{\pi_n}{\phi_n})(\theta)\, \dd Q_\theta(\theta)
\leq n\,\int_{c_0M\varepsilon_n}^\infty 3u^2 \,Q_\theta(\|\theta - \theta^\ast\| \geq u)\, \dd u  + e^{-Cc_0^2M^2\log n}+ C\alpha_n\\
&\qquad\qquad\qquad\leq  n\,\int_{c_0M\varepsilon_n}^\infty 3u^2 \,e^{-Cc_0^2M^2n u^2}\, \dd u  + e^{-Cc_0^2M^2\log n}+ C\alpha_n\\
&\qquad\qquad\qquad \leq   e^{-Cc_0^2M^2\log n}+ C\alpha_n \leq \frac{CM^3\,(\log n)^{d+3}}{\sqrt{n}},
\end{align*}
for sufficiently large constant $c_0$.

Putting all pieces together, we have
\begin{align*}
&\Big|D\big(Q_\theta\,||\,\Pi_n\big) -D\big(Q_\theta\,||\,N\big(\mle, \, [n I(\theta^\ast)]^{-1}\big)\big)\Big|=\Big| \int \log(\frac{\pi_n}{\phi_n})(\theta)\, \dd Q_\theta(\theta)\Big|\\
\leq &\,\sup_{\theta\in K_n}\Big|\log(\frac{\pi_n}{\phi_n})(\theta)\Big| 
+\Big| \int_{K_n^c} \log(\frac{\pi_n}{\phi_n})(\theta)\, \dd Q_\theta(\theta)\Big|\leq \frac{CM^3\,(\log n)^{d+3}}{\sqrt{n}},
\end{align*}
holds with probability at least $1-CM^{-2}$, which completes the proof.


\section{Proofs of technical lemmas}
In this appendix, we collect proofs of all technical lemmas.


\subsection{Proof of Lemma~\ref{lem:deIntApprox}}\label{app:proof_deIntApprox}
Denoting $\sqrt{n}(\theta - \mle)$ by $s_n$, we will approximate the integration of
$$i_1(\theta) =\pi(\theta) \exp(l(\theta; X^n) - l(\mle; X^n))$$
by
$i_2(\theta)= \pi(\theta^\ast) \exp(-\frac{1}{2}\,s_n^TI(\theta^\ast)\,s_n)$, where
\begin{equation*}
\int i_2(\theta)\, \dd \theta = \pi(\theta^\ast) (\frac{2\pi}{n})^{d/2}|I(\theta^\ast)|^{-1/2}.
\end{equation*}
Therefore, it suffices to show
\begin{equation}\label{Eqn:posterior_goal}
\Big|\int \big[i_1(\theta) - i_2(\theta)\big]\,\dd\theta \Big| \leq \frac{CM^3\,(\log n)^{d+3}}{\sqrt{n}}\, n^{-d/2}.
\end{equation}
We split the space $\mathbb R^d$ into two regions:
\begin{align*}
A_1 = \{ |s_n| \leq c_1 \log n\}\quad\mx{and}\quad \quad A_2 = \{ |s_n|> c_1 \log n\},
\end{align*}
for some sufficiently large constant $c_1$.

\vspace{0.5em}
\noindent {\bf Integral over $A_1$:} We start with $A_1$, by using the $\ell_\infty$-$\ell_1$ type Holder's inequality:
\begin{equation*}
\Big|\int_{A_1} \big[i_1(\theta) - i_2(\theta)\big]\,\dd\theta\Big| \leq \text{Vol}(A_1) \,\sup_{\theta\in A_1} \big|i_1(\theta - i_2(\theta)\big|,
\end{equation*}
where the volume $\text{Vol}(A_1)$ is $C\big(\frac{c_1 \log n}{\sqrt{n}}\big)^{d}$. Applying the Taylor expansion to $l(\theta; X^n)$ at $\theta=\mle$, we obtain
\begin{align}
\exp(l(\theta; X^n) - l(\mle; X^n)) = \exp(\frac{1}{2}(\theta - \mle)^T\nabla^2l(\mle; X^n)(\theta - \mle)   + R_n), \label{eq:TylorLR}
\end{align}
where we have used the fact that $\mle$ maximizes $l(\theta; X^n)$ so that $\nabla l(\mle; X^n)=0$, and the remainder term $R_n$ satisfies
\begin{equation*}
|R_n| \leq \sum_i |M(X_i)| (\frac{c_1 \log n}{\sqrt{n}})^3 \leq C\frac{c_1^3M^2(\log n)^3}{\sqrt{n}},
\end{equation*}
with probability at least $1-CM^{-2}$ by using Assumption A2 and the Markov inequality to the sum $\sum_{i=1}^n|M(X_i)|$ of i.i.d.~random variables.
As a consequence, we can bound the difference
\begin{align*}
&\Big|\exp(l(\theta; X^n) - l(\mle; X^n)) - \exp(-\frac{1}{2}\,s_n^TI(\theta^\ast)\,s_n)\Big| \\
=&\, \exp(-\frac{1}{2}\,s_n^TI(\theta^\ast)\,s_n)\, \Big|\exp(\frac{1}{2}\,s_n^T\big(I(\theta^\ast)+\frac{1}{n}\nabla^2l(\mle; X^n)\big)\,s_n + R_n) - 1\Big|.
\end{align*}
By the definition of $I(\theta^\ast)=n^{-1}E_{\theta^\ast}[\nabla^2 l(\theta^\ast;\, X^n)]$ in Assumption A2, the fact that $\|\mle -\theta^\ast\|\leq \frac{CM(\log n)^{1/2}}{\sqrt n}$ with probability at least $1-CM^{-2}$, and the inequality $|e^x-1|\leq 2x$ for $x\in[0,1/2]$, we can bound the above by using the Markov inequality on the i.i.d.~sum $\nabla^2 l(\theta^\ast;\, X^n)$ to obtain that
\begin{align*}
\Big|\exp(l(\theta; X^n) - l(\mle; X^n)) - \exp(-\frac{1}{2}\,s_n^TI(\theta^\ast)\,s_n)\Big|
\leq&\, \frac{Cc_1^2M^2(\log n)^3}{\sqrt{n}},\quad\mx{for all }s_n\in A_1,
\end{align*}
holds with probability at least $1-CM^{-2}$.
Using Assumption A1, we have $|\pi(\theta) - \pi(\theta^\ast)| \leq C c_1 \frac{\log n}{\sqrt{n}}$ under the aforementioned high probability event. Finally, by putting pieces together and using the triangle inequality, we have that with probability at least $1-CM^{-2}$,
\begin{align}\label{eqn:A_1_bound}
\Big|\int_{A_1} \big[i_1(\theta) - i_2(\theta)\big]\,\dd\theta\Big| \leq \frac{CM^3\,(\log n)^{d+3}}{\sqrt{n}}\, n^{-d/2}.
\end{align}

\vspace{0.5em}
\noindent {\bf Integral over $A_2$:}
We use the inequality
\begin{align*}
\Big|\int_{A_2} \big[i_1(\theta) - i_2(\theta)\big]\,\dd\theta\Big| \leq\Big|\int_{A_2} i_1(\theta)\,\dd\theta\Big| + \Big|\int_{A_2} i_2(\theta)\,\dd\theta\Big|,
\end{align*}
and bound the two terms separately. For the second term, we have
\begin{align}\label{eqn:A_2_bound}
\Big|\int_{A_2} i_2(\theta)\, \dd \theta\Big| &\leq C\, \exp(-C\,n\,\Big(c_1\frac{\log n}{\sqrt{n}}\Big)^2) \leq C\,n^{-Cc_1}\leq \frac{CM^3\,(\log n)^{d+3}}{\sqrt{n}}\, n^{-d/2},
\end{align}
where we choose $c_1$ large enough so that the last step is true.

For the first term, use inequality~\eqref{eqn:A_1_bound} and $\int_{A_2} i_2(\theta)\,\dd\theta \geq Cn^{-d/2}$, we have by the triangle inequality that 
\begin{align*}
\int_{A_1} i_1(\theta)\,\dd\theta \geq Cn^{-d/2}.
\end{align*}
From the form~\eqref{Eqn:posterior_form} of the posterior density $\pi_n$,the inequality $\|\mle-\theta^\ast\|\leq \frac{CM(\log n)^{1/2}}{\sqrt n}$ and the posterior tail probability bound~\eqref{Eqn:post_contraction}, we have that for sufficiently large constant $c_1$,
\begin{align*}
\frac{ \int_{A_2} i_1(\theta)\,\dd\theta}{\int_{A_1} i_1(\theta)\,\dd\theta+\int_{A_2} i_1(\theta)\,\dd\theta}& = \Pi_n\Big(\|\theta-\mle\|\geq \frac{c_1\log n}{\sqrt{n}}\Big) \\
&\leq \Pi_n\Big(\|\theta-\theta^\ast\|\geq \frac{c_1\log n}{2\sqrt{n}}\Big)\leq e^{-Cc_1^2 (\log n)^2}\leq \frac{C}{n},
\end{align*}
holds with probability at least $1-CM^{-2}$. The last two displays together imply
\begin{align*}
\int_{A_2} i_1(\theta)\,\dd\theta \leq \frac{Cn^{-1}}{1-Cn^{-1}} n^{-d/2} \leq Cn^{-d/2-1}.
\end{align*}
Now, we can combine the above and inequality~\eqref{eqn:A_2_bound} to obtain 
\begin{align}\label{eqn:A_2_bound}
\Big|\int_{A_2} \big[i_1(\theta) - i_2(\theta)\big]\,\dd\theta\Big| \leq  \frac{CM^3\,(\log n)^{d+3}}{\sqrt{n}}\, n^{-d/2}.
\end{align}

Finally, the desired bound~\eqref{Eqn:posterior_goal} follows by combining inequalities~\eqref{eqn:A_1_bound} and~\eqref{eqn:A_2_bound}.

\subsection{Proof of Lemma~\ref{lem:WPdenominator}}\label{app:proof_WPdenominator}
The proof of this lemma follows almost same lines as Lemma 8.1 in~\cite{ghosal2000convergence} that generalizes from the usual posterior to the weighted posterior.
Recall the definition of $\wt m_n(A)$ as:
\begin{equation*}
m_n(A) = \int_A \frac{\widetilde L^{(b)}(\theta;\, X^n)}{\widetilde L^{(b)}(\theta_0;\, X^n)} \, \dd \Pi(\theta).
\end{equation*}
To simplify the notation, we use shorthand $\mathrm{PB}(\delta)$ for the set $\{\theta\in \Theta: \|\theta - \theta^\ast \|^2 \leq C_0\delta^2 \}$.
Since the integrand is non-negative, we have
\begin{equation*}
m_n(\Theta) \geq \Pi(\mathrm{PB}(\delta)) \int_{\mathrm{PB}(\delta)} \frac{\widetilde L^{(b)}(\theta;\, X^n)}{\widetilde L^{(b)}(\theta_0;\, X^n)} \, \dd \Pi|_{\mathrm{PB}(\delta)}(\theta),
\end{equation*}
where recall that $\Pi|_A$ denotes the restriction of $\Pi$ on $A$.
By Jensen's inequality, we have further:
\begin{align*}
m_n(\Theta) \geq & \Pi(\mathrm{PB}(\delta)) \exp( \int_{\mathrm{PB}(\delta)} \log \frac{\widetilde L^{(b)}(\theta;\, X^n)}{\widetilde L^{(b)}(\theta_0;\, X^n)} \, \dd \Pi|_{\mathrm{PB}(\delta)}(\theta)) \\
= &  \Pi(\mathrm{PB}(\delta)) \exp(\int_{\mathrm{PB}(\delta)} \big(\widetilde l^{(b)}(\theta;\, X^n) - \widetilde l^{(b)}(\theta^\ast;\, X^n) \big)\, \dd \Pi|_{\mathrm{PB}(\delta)}(\theta))
\end{align*}
where $\widetilde l^{(b)}(\theta;\, X^n)$ denotes the logarithm of the weighted likelihood. Therefore, it remains to prove 
\begin{align*}
 \int_{\mathrm{PB}(\delta)} \big( \widetilde l^{(b)}(\theta;\, X^n) - \widetilde l^{(b)}(\theta^\ast;\, X^n)\big) \, \dd \Pi|_{\mathrm{PB}(\delta)}(\theta) \leq -(1 + C)n \delta^2.
\end{align*}
Assumption A2 implies that any $\theta$ in $\mathrm{PB}(\delta)$ satisfies
\begin{align*}
E_{\theta^\ast} \big[l(\theta^\ast; X) - l(\theta; X)\big] \leq C'\delta^2, \quad\mx{and}\quad E_{\theta^\ast} \big[l(\theta^\ast; X) - l(\theta; X)^2\big] \leq C'\delta^2.
\end{align*}
Since the weights $\{W_i^{(b)}\}_{i=1}^n$ have unit mean and variance, we have
\begin{equation*}
E_{\theta^\ast} E_W\bigg[ \int_{\mathrm{PB}(\delta)} \big( \widetilde l^{(b)}(\theta;\, X^n) - \widetilde l^{(b)}(\theta^\ast;\, X^n)\big) \, \dd \Pi|_{\mathrm{PB}(\delta)}(\theta)\bigg] \geq -C'n \delta^2,
\end{equation*}
and
\begin{align*}
& \mathrm{var}_{\theta^\ast, W}\bigg[ \int_{\mathrm{PB}(\delta)} \big( \widetilde l^{(b)}(\theta;\, X^n) - \widetilde l^{(b)}(\theta^\ast;\, X^n) \big)\, \dd \Pi|_{\mathrm{PB}(\delta)}(\theta)\bigg] \\
&= n \mathrm{var}_{\theta^\ast, w}\bigg[ W_1^{(b)} \int_{\mathrm{PB}(\delta)} \big( (( l(\theta;\, X_1) - l(\theta^\ast;\, X_1))\big) \, \dd \Pi|_{\mathrm{PB}(\delta)}(\theta)\bigg]  \\
&\leq n E_{\theta^\ast, w} \bigg[ (W_1^{(b)})^2 \int_{\mathrm{PB}(\delta)} \big(  l(\theta;\, X_1) - l(\theta^\ast;\, X_1)\big)^2 \, \dd \Pi|_{\mathrm{PB}(\delta)}(\theta)\bigg] 
\leq  2C'n\delta^2.
\end{align*}
by the independence.
We complete the proof by applying Chebyshev’s inequality, where the probability that the random variable $\int_{\mathrm{PB}(\delta)} \big( \widetilde l^{(b)}(\theta;\, X^n) - \widetilde l^{(b)}(\theta^\ast;\, X^n)\big) \, \dd \Pi|_{\mathrm{PB}(\delta)}(\theta)$ with variance less than $2C'n\delta^2$ deviates from its mean at least $C'n\delta^2$ is upper bounded by
$\frac{2C'n\delta^2}{(C'n\delta^2)^2} = \frac{2}{C'n\delta^2}$.

\subsection{Proof of Lemma~\ref{lem:WPnumerator}}\label{app:proof_WPnumerator}
By the sub-exponential condition in Assumption W, the distribution of the random weight has tail probability bounded by $P(W > t) \leq  C_0e^{-C_1t}$, for any $t >C_2$. This implies by a simple union bound argument that
\begin{align*}
& P\biggl( \bigcup_{i = 1}^n \{ W_i^{(b)} \geq C \log n \} \biggr)\leq C_0ne^{-C_1 C \log n} \leq C'n^{-1}
\end{align*}
by choosing a sufficiently large constant $C$. 
Thus, in the rest of the proof we may simply assume without loss of generality that the weights are uniformly bounded by $M_W := C \log n$.

The rest of the proof follows closely the steps in the proof of Theorem 1 in \cite{RatesPosterior}. Many intermediate lemmas therein can be reused for our proof except for two that need to be adapted to the weight posterior. First recall Bernstein's inequality for sum of i.i.d.~sub-exponential random variables: let $Z_1, \ldots, Z_n$ be i.i.d. r.v.'s satisfying:
\begin{equation*}
E[|Z|^j] \leq j!\, b^{j-2} v/2, \quad \mx{for some constant $b,v>0$ and for all $j \geq 2$}.
\end{equation*}
Let $\nu_n$ be the empirical process functional, i.e. $\nu_n(Z) = \sum_{i=1}^n (Z_i - E[Z])/\sqrt{n}$, then
\begin{equation*}
P(\nu_n(Z) > t) \leq \exp(-\frac{t^2}{4(2v + bt/\sqrt{n})}),\quad\forall t>0.
\end{equation*}
A key ingredient of the proof is a modification of Lemma 6 in \cite{RatesPosterior} to the weighted likelihood ratio process. Following their notation, for any density function we define the following lower-truncated log-likelihood ratio at some level $\tau$ to be determined later:
\begin{equation*}
Z_f := \max (\log \frac{f}{p^\ast}, -\tau),
\end{equation*}
where $p^\ast$ is the true density function $p_{\theta^\ast}$ in our case, and the truncation is needed since the log-likelihood ratio can be ill-behaved in the lower end. The truncated one can be uniformly controlled with exponentially decay tail probability. In our case, we introduce a version of the likelihood ratio empirical process as $\wt{\nu}_n(Z_f) = \sum_{i=1}^n (W_i^{(b)} Z_f(X_i) - E[W_i^{(b)}Z_f(X_i)])$ as analogous to $\nu_n$ in \cite{RatesPosterior}. 
\begin{lemma}\label{lem:technicalBernstein}
Under Assumption W, for any fixed density function $f$  with a finite Hellinger distance to $p^\ast$, we have
\begin{equation*}
P(\wt{v}_n(Z_f) \geq t) \leq \exp(-\psi_2(t, H^2(f, p_{\theta^\ast}), n)),
\end{equation*}
where function $\psi_2$ is defined as:
\begin{equation*}
\psi_2 (t, u, n) = \frac{Ct^2}{u + M_Wt/\sqrt{n}}, \quad\mx{for some constant $C>0$.}
\end{equation*}
\end{lemma}
\noindent The proof is deferred to Section~\ref{Sec:proof_lem:technicalBernstein}. The next lemma is a weighted version of Lemma 7 in~\cite{RatesPosterior} that provides a uniform control on the weighted likelihood ratio empirical process $\wt{\nu}_n(Z)$ via the chaining technique.
\begin{lemma}\label{lem:locent}
Suppose Assumption W holds. Consider any $t > 0$, $0 < k < 1$ and $R > 0$ such that
$R \leq k \sqrt{n} t^2 /4$ and 
\begin{equation*}
\int_{kR/(32\sqrt{n})}^t \mathrm{H}_B^{1/2} \Big( \frac{u}{2\exp(\tau/2)}, \big\{p_\theta,\,\theta\in\Theta:\, H(P_\theta,\, P_{\theta^\ast})\leq t\big\}  \Big) \dd u \leq \frac{Rk^{3/2}}{2^{11}(c_0 + 1/8)M_W},
\end{equation*}
where $c_0 = (\exp(\tau/2) - 1 -\tau/2)/(1 - \exp(-\tau/2))^2$ is a constant. Then
\begin{equation*}
P_{\theta^\ast}\Big(\sup_{\theta\in\Theta:\, H(P_\theta,\, P_{\theta^\ast})\leq t} \wt{v}_n(Z_{p_\theta}) \geq R\Big) \leq 3\exp(-(1-k)\psi_2(R, t^2, n)).
\end{equation*}
\end{lemma}
\noindent A proof of this lemma is provided in Section~\ref{Sec:proof_lem_locent} below.
Now we can continue the proof of the lemma, which resemble the proof of Theorem 1 in~\cite{wong1995}. We first verify the condition of Lemma~\ref{lem:locent} with $s=C_2'c_2{-1/2}\varepsilon\geq C_2'M\wt\varepsilon_n^2$, $t = \sqrt{2}s$, $1/2 < k < 1$, $\exp(-\tau/2)=1/5$ and $R = \frac{k\sqrt{n}s^2}{2}$ under Assumption A2 (where $c_2$ is the constant in Assumption A2). Specifically, for a sufficiently large constant $C_2'$, we have
\begin{align*}
&\int_{s^2/2^8}^{\sqrt{2}s} \mathrm{H}_B^{1/2}\Big(\frac{u}{2\exp(\tau/2)}, \big\{p_\theta,\,\theta\in\Theta:\, H(P_\theta,\, P_{\theta^\ast})\leq \sqrt{2s}\big\} \Big) \dd u\\
&\qquad\qquad \leq  C\,s \int_{0}^1 \sqrt{\log\frac{\exp(\tau/2)}{u}}\,\dd u \leq \frac{C'c_2^{1/2}\sqrt{n}s^2}{C_2'M_W}
\leq  \frac{Rk^{3/2}}{2^{11}(c_0 + 1/8)M_W}.
\end{align*}
Moreover, for the above inequality, since the left hand side has a linear growth in $s$ while the right hand side has a quadratic growth, it also holds for all $s \geq C_2'c_2^{-1/2} \varepsilon$.

The remaining of the proof is almost identical to the proof of Theorem 1 in \cite{RatesPosterior}. As a consequence of Lemma~\ref{lem:locent} and the above argument, we have:
\begin{equation*}
P\bigg(\sup_{\theta\in\Theta:\, H(P_\theta,\, P_{\theta^\ast})\leq \sqrt{2}s} \wt{\nu}_n (Z_{p_\theta}) \geq \frac{k}{2}\sqrt{n}s^2\bigg) 
\leq 3\exp(-\frac{Cns^2}{1+M_W}),\quad\forall s\geq C_2'c_2^{-1/2}\varepsilon.
\end{equation*}

By the equivalence between the Hellinger distance and the $\ell_2$ distance in Assumption A2, we can decompose $\{\theta:\,\|\theta-\theta^\ast\| \geq C_2'\varepsilon\} = \bigcup_{i=0}^\infty \big(\Theta^{(i+1)}\setminus \Theta^{(i)}\big)$ with $\Theta^{(i)}:=\{\theta\in\Theta:\, H(P_\theta,\, P_{\theta^\ast})\leq 2^{i/2}C_2'c_2^{-1/2}\varepsilon\}$. Lemma 4 of~\cite{wong1995} implies $E\wt Z_{p_\theta} \leq (1-\delta_0) H^2(P_\theta,P_{\theta^\ast}) \leq (1-\delta_0)2^{i}(C_2')^2c_2^{-1}\varepsilon^2$ for each $\theta\in\Theta^{(i)}$, where  $\delta_0 = 2\exp(-\tau/2)/(1 - \exp(-\tau/2))^2$.
Therefore, a simple union bound with the previous decomposition and the preceding display (also with the fact that $E W_i^{(b)}=1$) renders:
\begin{align*}
P\bigg(\sup_{\theta:\,\|\theta-\theta^\ast\| \geq C_2'\varepsilon} \prod \Big(\frac{p_\theta(X_i)}{p_{\theta^\ast}(X_i)}\Big)^{W^{(b)}_i} \geq &\, \exp\Big(-\big((1 - k/2)(C_2')^2c_2{-1}  - \delta_0\big)\,n\varepsilon^2\Big) \bigg)\\
&\qquad\qquad\leq 6\exp(- \frac{C(1 - k)k^2 (C_2')^2c_2^{-1}n\varepsilon^2}{1+M_W}).
\end{align*}
As a consequence, by choosing $\exp(-\tau/2) = 1/5$ and $k = 2/3$, and using the fact that $n\varepsilon^2/\log n\geq C'M^2$ for some constant $C'$ and any $\varepsilon\geq M\wt \varepsilon=MC_0' \log n/\sqrt{n}$, we obtain that by choosing a sufficiently large $C_2'$, the following holds with probability at least $1-e^{-C''M^2}\geq 1- CM^{-2}$,
\begin{equation*}
\sup_{\theta:\,\|\theta-\theta^\ast\| \geq C_2'\varepsilon} \prod \Big(\frac{p_\theta(X_i)}{p_{\theta^\ast}(X_i)}\Big)^{W^{(b)}_i} \leq \exp(-n\varepsilon^2/24).
\end{equation*}


\subsection{Proof of Lemma~\ref{lem:technicalBernstein}}\label{Sec:proof_lem:technicalBernstein}
The proof is a direct application of Bernstein's inequality. We only need to verify the Bernstein condition for the random variable $\wt Z_f$. In particular, similar to Lemma 5 of~\cite{wong1995}, it is easy to verify that
\begin{equation*}
E[|W_i^{(b)}Z_f(X_i)|^j] \leq (j!)\, 2^jc_0 H^2(f, p_{\theta^\ast}), \quad \forall j \geq 2,
\end{equation*}
with $c_0 = \frac{\exp(\tau/2) - 1 - \frac{\tau}{2}}{(1 - \exp(-\tau/2))^2}$.
Since the weights $W_i^{(b)}$'s are uniformly bounded by $M_W = C\log n$ and have second moment $E[W_i^{(b)}]^2=2$, we obtain $E[|w Z_f|^j] \leq (j!)2^{j+1} M_W^{j-2}c_0 H^2(f, p_{\theta^\ast})$ for any $j\geq 2$. Therefore, we can then directly apply Bernstein's inequality with $b = 2M_W$, and $v = 8c_0H^2(f, p_{\theta^\ast})$ to obtain the claimed inequality.


\subsection{Proof of Lemma~\ref{lem:locent}}\label{Sec:proof_lem_locent}
In the following, we provide a sketched proof as an extension of the proofs of Lemma 7 in~\cite{RatesPosterior} and Theorem 3 in~\cite{Wong1994} to the weighted likelihoods. We also keep their notation, and use 
the shorthand $\mathrm{HB}(R) := \{p_\theta,\,\theta\in\Theta:\, H^2(p_{\theta}, p_{\theta^\ast}) \leq R^2 \}$ to denote  a Hellinger ball centered at $p_{\theta^\ast}$ with radius $R > 0$. It is proved in lemma 3 of~\cite{RatesPosterior} that 
\begin{equation*}
\| Z_{f_1} - Z_{f_2}\|_2^2 \leq 4\exp(\tau) H^2(f_1, f_2)
\end{equation*}
which will be repeatedly used throughout the proof.

The two bracketing sequence $\{ \mathscr{Z}_j \}_{j \leq N}$ is constructed in the same way as in Section 2.1 of \cite{RatesPosterior}, whose members have pairwise $L_2$ distance upper bounded by a sequence $\{ 2\exp(\tau/2)\delta_j \}_{j \leq N}$  defined as:
\begin{align*}
\delta_0 =& \inf_{u} \{ \mathrm{H}_B(u, \mathrm{HB}(t)) \geq \frac{k}{4} \psi_2(R, t^2, n)\},\\
\delta_{j+1} =& \max\biggl( \frac{kR}{8\sqrt{n}}, \sup \bigl\{ x\leq \frac{\delta_j}{2}: \mathrm{H}_B(x, \mathrm{HB}(t)) \geq 4 \mathrm{H}_B(\delta_j, \mathrm{HB}(t)) \bigr\} \biggr),\\
N =& \min \{ j: \delta_j \leq \frac{kR}{8\sqrt{n}} \}.
\end{align*}
If $\frac{kR}{8\sqrt{n}} \leq t$, then similar to the proof of Theorem 3 in~\cite{Wong1994} the probability on the left hand side of the desired inequality can be decomposed into four terms, 
\begin{equation*}
P_{\theta^\ast}\Big(\sup_{\mathrm{HB}(t)} \wt{v}_n(Z_{p_\theta}) \geq R\Big) \leq P_1 + P_2 + P_3 + P_4,
\end{equation*}
where $P_i$, $i=1,2,3,4$, is (the definition of $u_j$ below is on page 601 of \cite{Wong1994})
\begin{align*}
P_1 &= |\mathscr{Z}_0|\sup_{u_0} P(\wt{\nu}_n (Z_{u_0}) > (1 - \frac{k}{4}) R), \\
P_2 &= \sum_{j=1}^N \prod_{l=0}^{j-1} |\mathscr{Z}_l| \prod_{m=0}^{j-1} |\mathscr{Z}_m| \sup_{u_j, u_{j-1}} P(\wt{v}_n(Z_{u_j} - Z_{u_{j-1}})I_{\bigcup_{k\geq j}B_k} > M_W \eta_j),\\
P_3 &= \sum_{j = 0}^{N - 1} P^\ast( \sup_{\mathrm{HB(t)}} (\wt{\nu}_n (Z_{p_\theta} - Z_{u_j}) I_{B_j}) > M_W\eta_{j + 1}),\\
P_4 &= P^\ast(\sup_{\mathrm{HB}(t)}\wt{v}_n(Z_{p_\theta}-Z_{u_N})I_{B_N} > \frac{kR}{8}),
\end{align*}
with the following changed definitions for $a_j$ and $\eta_j$:
\begin{align*}
\eta_j =& 16\exp(\tau/2) \delta_{j-1} \sqrt{\frac{\sum_{l\leq j} \mathrm{H}_B(\delta_l, \mathrm{HB}(t))}{k}},\\
a_j =& \frac{32\sqrt{n}\exp(\tau) \delta_{j-1}^2}{\eta_j}
\end{align*}
and the same $B_j$ on page 602 of \cite{Wong1994} except for the underlying family of function to be $Z_f$. Following the steps in their proof, it can be verified that
\begin{equation*}
\sum_{j = 0}^N \eta_j \leq \frac{2^7}{\sqrt{k}} \int_{kR/(32\sqrt{n})}^t \mathrm{H}_B^{1/2} \bigl( \frac{u}{2\exp(\tau/2)}, \mathrm{HB}(t) \bigr) \dd u
\leq \frac{kR}{8M_W}.
\end{equation*}

Now it remains to bound $P_1$--$P_4$ respectively. First, by applying the inequality in our Lemma~\ref{lem:technicalBernstein} to replace their unweighted version and using the same argument, we have
\begin{equation*}
P_1 \leq \exp(-(1 - k)\psi_2(R, t^2, n)).
\end{equation*}
Similarly, by applying a one-sided Bernstein's inequality with their argument, we have
\begin{equation*}
P_2  \leq 2\exp(-(1 - k)\psi_2(R, t^2, n)),
\end{equation*}
where the only difference is in replacing the $L_2$ bracketing entropy with the Hellinger bracketing entropy (they are equivalent up to a constant). Similarly to their argument on page 605-606, we also have $P_3 = P_4 = 0$ when $\frac{kR}{8\sqrt{n}} \leq t$.

The case of $\frac{kR}{8\sqrt{n}} > t$ can be similarly proceeded as in \cite{Wong1994} and we omit the details.

\end{document}